\documentclass[11pt]{imsart}
\RequirePackage[OT1]{fontenc}
\RequirePackage{amsthm,amsmath}
\RequirePackage[numbers]{natbib}
\RequirePackage[colorlinks,citecolor=blue,urlcolor=blue]{hyperref}

\usepackage{graphicx,amsmath,amsthm,amsfonts,amssymb}
\usepackage{booktabs, rotating, float}
\usepackage[hmargin=3cm,vmargin=3cm]{geometry}
\usepackage{mathtools}
\usepackage{bbm}
\mathtoolsset{showonlyrefs}
\usepackage{enumerate}
\usepackage{color}
\usepackage{listings}
\usepackage{epstopdf}
\usepackage{makecell}
\usepackage{tikz}
\usetikzlibrary{decorations.pathreplacing}
\usepackage{array}
\usepackage{caption}
\usepackage{subcaption}
\usepackage{soul}
\usepackage{hyperref}
\usepackage{verbatim}
\usepackage{siunitx}
\usetikzlibrary{calc} \tikzset{>=latex}
\usetikzlibrary{calc,decorations.markings}
\usetikzlibrary{positioning,shapes,decorations.text}
\theoremstyle{plain}
\newtheorem{thm}{Theorem}[section]
\newtheorem{cor}[thm]{Corollary}

\newtheorem{prop}[thm]{Proposition}
\newtheorem{lemma}[thm]{Lemma}
\newtheorem{assumption}[thm]{Assumption}
\theoremstyle{definition}
\newtheorem{defn}[thm]{Definition}

\newtheorem{eg}[thm]{Example}
\newtheorem{fact}[thm]{Fact}
\newtheorem{observe}[thm]{Observation}

\numberwithin{equation}{section}

\newcommand{\rpm}{\sbox0{$1$}\sbox2{$\scriptstyle\pm$}
  \raise\dimexpr(\ht0-\ht2)/2\relax\box2 }

\tikzstyle{nd} = [anchor=base, inner sep=0pt]
\tikzstyle{ndpic} = [remember picture, baseline, every node/.style={nd}]

\def\beq{\begin{equation}}
\def\eeq{\end{equation}}
\def\ba{\begin{enumerate}[(a)]}
\def\bei{\begin{enumerate}[(i)]}
\def\be{\begin{enumerate}[(1)]}
\def\ee{\end{enumerate}}
\def\bi{\begin{itemize}}
\def\ei{\end{itemize}}
\def\beg{\begin{eg}}
\def\eeg{\end{eg}}
\def\bd{\begin{defn}}
\def\ed{\end{defn}}
\def\bt{\begin{thm}}
\def\et{\end{thm}}
\def\bl{\begin{lemma}}
\def\el{\end{lemma}}
\def\bfac{\begin{fact}}
\def\efac{\end{fact}}

\def\bc{\begin{cor}}
\def\ec{\end{cor}}
\def\bp{\begin{prop}}
\def\ep{\end{prop}}
\def\bo{\begin{observe}}
\def\eo{\end{observe}}
\def\bas{\begin{assumption}}
\def\eas{\end{assumption}}
\def\RR{\mathbb{R}}
\def\CC{\mathbb{C}}

\def\ZZ{\mathbb{Z}}
\def\NN{\mathbb{N}}

\def\ii{\item}

\def\beg{\begin{eg}}
\def\eeg{\end{eg}}

\def\p{\thinspace}

\makeindex

\numberwithin{equation}{section}
\numberwithin{table}{section}
\newtheorem{remark}{Remark}

\startlocaldefs
\endlocaldefs

\begin{document}

\begin{frontmatter}

\title{Hall-Littlewood-PushTASEP and its KPZ limit}
\runtitle{HL-PushTASEP \& its KPZ limit}


\begin{aug}
  \author{\fnms{Promit} \snm{Ghosal}\ead[label=e1]{pg2475@columbia.edu}\vspace{0.1in} \\ {\it Columbia University}}

   \affiliation{Columbia University}


  \runauthor{Ghosal}


 \address{ Department of Statistics, 1255 Amsterdam Avenue, New York, NY 10027\\ \printead{e1}}

\end{aug}

\begin{abstract}
\noindent We study a new model of interactive particle systems which we call the randomly activated cascading exclusion process (RACEP). Particles wake up according to exponential clocks and then take a geometric number of steps. If another particle is encountered during these steps, the first particle goes to sleep at that location and the second is activated and proceeds accordingly. We consider a totally asymmetric version of this model which we refer as Hall-Littlewood-PushTASEP (HL-PushTASEP) on  $\ZZ_{\geq  0}$ lattice where particles only move right and where initially particles are distributed according to Bernoulli product measure on $\ZZ_{\geq 0}$. We prove KPZ-class limit theorems for the height function fluctuations. Under a particular weak scaling, we also prove convergence to the solution of the KPZ equation.
\end{abstract}




\end{frontmatter}

\tableofcontents

\section{Introduction}
%
%
%
%
%

We introduce and study a special class of interacting particle systems which we call {\it randomly activated cascading exclusion processes} (RACEP). Fix a directed graph $G=(V,E)$ with conductances $c_e$ (i.e., non-negative weights) along directed edges $e\in E$ and define a random walk measure as the Markov chain on vertices with transition probability from $v$ to $v'$ given by $c_{v\to v'}$ normalized by the sum of all conductances out of $v$. The state space for RACEP is $\{0,1\}^{V}$ where $1$ denotes a particle and $0$ denotes a hole. Each particle is activated according to independent Poisson clocks. Once active, a particle moves randomly to one its adjacent sites. Therein, it chooses an independent random number of steps according to a geometric distribution and then performs an independent random walk (according to the random walk measure just described) of that many steps. However, if along that random walk trajectory, the first particle arrives at a site occupied by a second particle, then, the first particle stays at that site (and goes back to sleep) and the second particle becomes active and proceeds as if its Poisson clock had rung (according to the above rules). The Poisson clock activation processes are supposed to model a much slower process than the random walk cascades, so for simplicity we assume that these random walk cascades occur instantaneously. There is some work necessary to prove well-definedness of this process in general, and we do not pursue that here since we will mainly focus on one concrete case.

RACEP is a special type of exclusion process (see, e.g. \cite{Liggett85,Liggett99}) in which the jump distribution depends on the state around the particle in a rather non-trivial way. Exclusion processes are important models of lattice gases, transport (such as in various ecological or biological contexts), traffic, and queues in series (see for example \cite{tagkey2011iv,Corwin12} and references therein). The prototypical example of an exclusion process is the asymmetric simple exclusion process (ASEP) which was introduced in the biology literature in 1968 \cite{MGP68} (see also more recent applications such as in \cite{Chowdhury05,Greulich10}) and two years later in the probability literature \cite{Spitzer70}.

RACEP can also be thought of as a variant of a {\it frog model}  (see, e.g. \cite{AMP02a,AMP01}) in which particles do not fall asleep once active but continue to move and wake up other particles. A natural interpolation between RACEP and the frog model is to have the probability that a particle goes back to sleep be in $(0,1)$. These relations to exclusion processes and frog models suggest a number of natural probabilistic questions for RACEP, such as understanding its hydrodynamic and fluctuation limit theorems for various families of infinite graphs with simple conductances. We attack these questions exactly for the special case of RACEP in which the underlying graph is $\ZZ_{\geq 0}$ and the random walk is totally asymmetric (in the positive direction). From here on out, we will only focus on this one-dimensional case of RACEP.

ASEP and RACEP on $\ZZ$ are siblings in that they both arise as (different) special limits of the stochastic six vertex model \cite{GS92,BCG2016} -- see Section~\ref{SSVMtoHL-PushTASEP}. They are, in fact, part of a broader hierarchy of integrable probability particle systems which are solvable due to connections to quantum integrable systems -- see \cite{Corwin2016, BP16a, BP16b}. A totally asymmetric version of RACEP also comes up as a marginal of certain continuous time RSK-type dynamics which preserve the class of Hall-Littlewood processes -- see \cite{BBW16} wherein they refer to RACEP as the $t$-pushTASEP due to its similarity with the model of $q$-PushTASEP \cite{BP16a,CorwinPetrov2015}. Here, we are mostly interested in this one dimensional totally asymmetric specialization of RACEP. To avoid the confusion with the time parameter $t$, we refer RACEP as Hall-Littlewood-PushTASEP (HL-PushTASEP) throughout the rest of the paper.

These systems enjoy many concise and exact formulas from which one can readily perform asymptotics. Such results for ASEP go back to the now seminal papers \cite{Tracy2008,TracyWidom11} (and earlier to \cite{Kurt00} for TASEP), and results for the stochastic six vertex model were worked out in \cite{BCG2016,AB16}. In our present paper, we work out the analogous asymptotics for the HL-PushTASEP on $\ZZ_{\geq 0}$ using the approach of \cite{BCG2016}. (We also provide some additional details in the asymptotics compared to the previous works.) Owing to the recent work \cite{B16c,BO16} such asymptotics could alternatively be performed via reduction to Schur process asymptotics. We do not pursue that route here and instead opt for the more direct (albeit technically demanding) route.

The asymptotic results we show for HL-PushTASEP on $\ZZ_{\geq 0}$ (as well as those previously shown for ASEP and the stochastic six vertex model) demonstrate its membership in the Kardar-Parisi-Zhang universality class (see, e.g. the review \cite{Corwin12} and references therein). In particular, in Theorem \ref{Asymptheo1Mainstate} we show that for step Bernoulli initial data, HL-PushTASEP has fluctuations in its rarefaction fan of order $t^{1/3}$ and with statistics determined by the GUE Tracy-Widom distribution. At the edge of the rarefaction fan, we also demonstrate the occurrence the Baik-Ben Arous-P\'{e}ch\'{e} crossover distributions in \citep{BBP05} which arises also in ASEP \cite{Tracy2009}, TASEP \cite{BBP05, BenCor11}, and the stochastic six vertex model \cite{AB16}. See \cite{CorwinICM2014} for further references to asymptotics of other integrable probabilistic systems. The centering and scaling that our result demonstrates agrees with the predictions of KPZ scaling theory from the physics literature -- see Section \ref{KPZScaling}.

The KPZ scaling theory also suggests that under certain {\it weak scalings} a broad class of interacting particle systems should converge to the Cole-Hopf solution to the KPZ stochastic PDE (i.e., the KPZ equation). There are, in fact, many different choices of weak scalings under which the KPZ equation remains statistically invariant (see, e.g.\cite[Section~2.4]{QS15}). The main motivation for our present investigation into the KPZ equation limit for HL-PushTASEP was our desire to solve the analogous question for the stochastic six vertex model. Let us briefly explain why that question is difficult and what our HL-PushTASEP results suggest regarding it. 

For ASEP under weak asymmetry scaling, the KPZ equation limit was proved in \cite{Bertini1997} (see also \cite{Amir11,Dembo2016,CSL16}). That work relies upon two main identities. The first is the  G\"{a}rtner (or microscopic Cole-Hopf) transform which turns ASEP into a discrete stochastic heat equation (SHE). The second is a non-trivial key identity (Proposition~4.8 and Lemma~A.1 in \cite{Bertini1997}) which allows one to identify white-noise as the limit of the martingale part of the discrete SHE. For the stochastic six vertex model, one still has an analogous G\"{a}rtner transform (in fact, this holds for all higher-spin vertex models in the hierarchy of \cite{Corwin2016} due to the duality shown therein). The discrete time nature of the stochastic six vertex model, however, renders the identification of the white-noise much more complicated and presently it is unclear how to proceed. 

For higher-spin vertex models with unbounded particle occupation capacity, \cite{CL15} proved convergence to the KPZ equation under certain weak scalings. Even though these models are still discrete time, the scaling was such that the overall particle density goes to zero with the scaling parameter $\epsilon$. Owing to that fact, there was no need for an analog of the key estimate used in the case of ASEP.

Our present analysis of HL-PushTASEP also involves a weak scaling under which the local density of particles goes to zero. Consequently, we are able to prove the KPZ limit without the key estimate. This suggests that it may be simpler to derive the KPZ equation directly from the stochastic six vertex model (or other finite-spin integrable stochastic vertex models) when the weak scaling facilitates local density decay to zero. Indeed, there are many different weak scalings which should all lead to the KPZ equation. This can be anticipated, for instance, from analyzing moment or one-point distribution formulas, see e.g. \cite[Section~12]{BO16}. We intend on investigating such zero-density KPZ equation limits of the stochastic six vertex model in subsequent work.

There is another approach for proving KPZ limit of particle system (when started from their invariant measure) via energy solutions. This approach was introduced by Assing \cite{Assing2002} and then significantly developed in \cite{GJ2014a,GP2015a}. The invariant measure is  Bernoulli product measure for HL-PushTASEP (as well as for the stochastic six vertex model), so it would be interesting to see if these methods apply. Presently, the energy solution approach is only developed for continuous time systems, so an analysis of the stochastic six vertex model results would require developing a discrete time variant.

\subsection*{Outline}

We introduced above RACEP in an arbitrary infinite directed graph. In Section~\ref{ModelMainResult}, we describe the dynamics of HL-PushTASEP in one dimension and its connection with stochastic six vertex model. Further, we state our results on asymptotic fluctuation and KPZ limit of HL-PushTASEP in Section~\ref{ModelMainResult}.  At the end of this section, we explain the KPZ scaling theory for HL-PushTASEP in brief. In Section~\ref{TransitionMatrixLaplaceTransfom}, we first derive the eigenfunctions of the transition matrix of HL-PushTASEP and then use it to get Fredholm determinant formulas. Section~\ref{secAsymptotics} contains the proof of the limiting fluctuation results under the assumption of step Bernoulli initial data. Interestingly, we see a phase transition in the limiting behavior of the fluctuation. In Section~\ref{WeakScaleLim}, we prove the KPZ limit theory for HL-PushTASEP. First, we start with the construction of a discrete SHE. Next, we prove few estimates on the discrete heat kernel and subsequently, show the tightness of the rescaled SHE. At the end, we elaborate on the equivalence of all the limit points by solving the martingale problem for HL-PushTASEP.          

\subsection*{Acknowledgements}
Special thanks are due to Ivan Corwin for suggesting the problem and guiding throughout this work. We thank Guillaume Barraquand and Peter Nejjar for several constructive comments on the first draft of this paper. We would also like to thank Li-Cheng Tsai and Rajat Subhra Hazra for helpful discussions.

\section{Model and main results}\label{ModelMainResult}

We now define the HL-PushTASEP on $\ZZ_{\geq 0}$, describe its connection to the stochastic six vertex model, and then state our main results. Theorems~\ref{Asymptheo1Mainstate} and ~\ref{Asymptheo2Mainstate} give the fluctuations limits for HL-PushTASEP started with step Bernoulli initial data. Theorem~\ref{KPZLimitTheo} and ~\ref{KPZLimitTheoAtStepInit}  contain the KPZ equation limits for HL-PushTASEP under weak scaling given in Definition~\ref{GartnerTransform}.

\subsection{HL-PushTASEP on $\ZZ_{\geq 0}$}
We will consider HL-PushTASEP with a finite number of particles supported on $\ZZ_{\geq 0}$. The dynamics of HL-PushTASEP are such that the behavior of particles restricted to the interval $[0,L]$ is Markov. Thus, even if we are interested in HL-PushTASEP with an infinite number of particles, if we only care about events involving the restriction to finite intervals, then it suffices to consider the finite particle number version we define here.

The $N$-particle HL-PushTASEP on $\ZZ_{\geq 0}$ has state at time $t$ given by $(x_1(t),\ldots ,x_N(t))$ where $x_i(t)\in \ZZ_{\geq 0}$ for all $t\in \RR^+$ and $x_1(t)<\ldots <x_N(t)$. For notational simplicity we add two virtual particles fixed for all time so that $x_0=-1$ and $x_{N+1}=+\infty$. For later reference, we denote the state space of $N$-particle HL-PushTASEP on $\ZZ_{\geq 0}$ by
\begin{align}\label{eq:SateOfSystem}
X^N:=\big\{\vec{x}=(-1=x_0<x_1<x_2<\ldots <x_N<x_{N+1}=\infty):\forall i, x_i\in \ZZ_{\geq 0}\big\}.
\end{align}

\begin{figure}
\begin{tikzpicture}
\shade[ball color=white] (-2,0) circle (1ex);
\shade[ball color= white] (0,0) circle (1ex);
\shade[ball color= green] (7,0) circle (1ex);
\shade[ball color= red] (7,0.5) circle (1ex);
\shade[ball color= black ] (11,0) circle (1ex);
\draw[-] (-3,0) -- (12,0) ;
\draw[-] (-2,-0.05) -- (-2,0.05) node[below] {$x_{m-2}(t)$};
\draw[-] (-1,-0.05) -- (-1,0.05);
\draw[-] (0,-0.05) -- (0,0.05) node[below]  {$x_{m-1}(t)$};
\draw[-] (1,-0.05) -- (1,0.05);
\draw[-] (2,-0.05) -- (2,0.05);
\draw[-] (3,-0.05) -- (3,0.05);
\draw[-] (4,-0.05) -- (4,0.05) node[below] {$x_{m}(t)$};
\draw[-] (5,-0.05) -- (5,0.05) ;
\draw[-] (6,-0.05) -- (6,0.05) ;
\draw[-] (7,-0.05) -- (7,0.05) node[below] {$x_{m+1}(t)$};
\draw[-] (8,-0.05) -- (8,0.05) ;
\draw[-] (9,-0.,0.05) -- (9,0.05) ;
\draw[-] (10,-0.05) -- (10,0.05) ;
\draw[-] (11,-0.05) -- (11,0.05) node[below] {$x_{m+2}(t)$};
\node[above] at (7,0) (A) {};
\node[below] at (11,0.2) (D) {};
\node[above] at (7,0.3) (B) {};
\node[below] at (10,0.2) (C) {};
\draw[->] (A)to [out=145,in=158,looseness=6] (B)  ;
\draw[->] (B)to [out=30,in=150,looseness=1] (C) ;
\draw[->] (B)to [out=60,in=150,looseness=1] (D) ;
\end{tikzpicture}
\caption{Once red particle has been activated, it moves by one to the right, and then proceeds according to the same rules which had governed the green particle. For instance, it can jump by $j$ steps with probability $(1-b)b^{j-1}$ for $1\leq j<x_{m+2}(t)-x_{m+1}(t)$ and with probability $b^{j-1}$ if $j=x_{m+2}(t)-x_{m+1}(t)$. We illustrated here the cases when $j$ equals $x_{m+2}(t)-x_{m+1}(t)$ and $x_{m+2}(t)-x_{m+1}(t)-1$.}
\label{fig:2}
\end{figure}
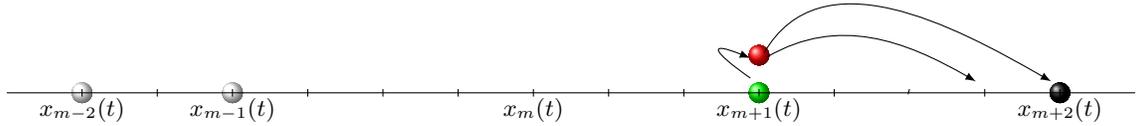

We illustrate typical moves in HL-PushTASEP on $\ZZ_{\geq 0}$ in Fig~\ref{fig:1} and Fig~\ref{fig:2}.
Each of the particles in HL-PushTASEP carries an exponential clock of rate $1$. When the clock in a particle rings, it jumps to the right. Unlike the simple exclusion processes, if the particle at $x_m(t)$ gets excited at time $t$, then it can jump $j$ steps to the right with a geometric probability $(1-b)b^{j-1}$ for $1\leq j<x_{m+1}(t)-x_m(t)$. It can also knock its neighbor on the right side at $x_{m+1}(t)$ out from its position with probability $b^{x_{m+1}(t)-x_m(t)-1}$. In that case, particle at position $x_{m+1}(t)$ will also get excited and jumps in the same way independently to others. To be precise, the HL-PushTASEP shares a large extent of resemblance with the particle dynamics in the case of stochastic six vertex model (SSVM) (see, \cite[Section~2.2]{BCG2016}). Although the dynamics in latter case is governed by a discrete time process, but pushing effect of the particles comes into play in the same way. We embark more upon the connection between these two model later in the following section. One can see similar pushing mechanism in the case of Push-TASEP \cite{BorFer08}, q-PushASEP\cite{CorwinPetrov2015} which are also continuous time countable state space Markov processes. But in contrast with the HL-PushTASEP, each of the particles in those cases can only jump by one step when their clocks ring.

\begin{figure}
\begin{tikzpicture}
\shade[ball color=white] (-2,0) circle (1ex);
\shade[ball color= white] (0,0) circle (1ex);
\shade[ball color= green] (4,0) circle (1ex);
\shade[ball color= red] (7,0) circle (1ex);
\shade[ball color= black ] (11,0) circle (1ex);
\draw[-] (-3,0) -- (12,0) ;
\draw[-] (-2,-0.05) -- (-2,0.05) node[below] {$x_{m-2}(t)$};
\draw[-] (-1,-0.05) -- (-1,0.05);
\draw[-] (0,-0.05) -- (0,0.05) node[below]  {$x_{m-1}(t)$};
\draw[-] (1,-0.05) -- (1,0.05);
\draw[-] (2,-0.05) -- (2,0.05);
\draw[-] (3,-0.05) -- (3,0.05);
\draw[-] (4,-0.05) -- (4,0.05) node[below] {$x_{m}(t)$};
\draw[-] (5,-0.05) -- (5,0.05) ;
\draw[-] (6,-0.05) -- (6,0.05) ;
\draw[-] (7,-0.05) -- (7,0.05) node[below] {$x_{m+1}(t)$};
\draw[-] (8,-0.05) -- (8,0.05) ;
\draw[-] (9,-0.,0.05) -- (9,0.05) ;
\draw[-] (10,-0.05) -- (10,0.05) ;
\draw[-] (11,-0.05) -- (11,0.05) node[below] {$x_{m+2}(t)$};
\node[above] at (4,0.1) (A) {};
\node[below] at (4,-0.2) (D) {};
\node[above] at (7,0.1) (B) {};
\node[below] at (6,0.1) (C) {};
\draw[->] (A)to [out=22,in=158,looseness=1] (B)  ;
\draw[->] (A)to [out=15,in=155,looseness=1] (C) ;
\end{tikzpicture}
\caption{Assume that particle at $x_m(t)$ (colored green) becomes active at time $t$. It immediately takes  a step to the right and then continues to jump according to a geometric distribution with parameter $b$. Thus, the probability that it takes a jump of exactly size $2$ is $b(1-b)$. With probability $b^2$ it reaches the position of $x_{m+1}(t)$ (colored red). In that scenario, the green particle occupies the position $x_{m+1}(t)$ and activates the red particle as shown in Fig. \ref{fig:2}.}

\label{fig:1}
\end{figure}
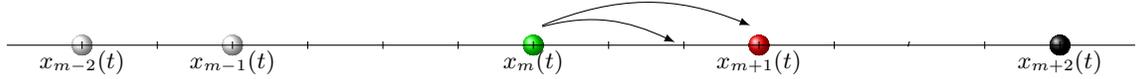

For HL-PushTASEP on $\ZZ_{\geq 0}$, we are interested in studying the evolution of the empirical particles counts, i.e, for any $\nu\in \RR$, what is the limit shape and fluctuation of the number of particles in an interval $[0,\lfloor \nu t\rfloor]$ as the time evolves to infinity. One can note that empirical particle counts fits into the stereotype of the height function for the exclusion processes. In what follows, we provide the explicit definitions of the height function and the step Bernoulli initial condition.

  \bd\label{eq:ParticleCountDefinition}
  Fix a positive integer $N$. Define a function $N_x(t;.)$ and $\eta_t(x;.)$ of $\vec{x}=(x_1<x_2<\ldots<x_N)\in X^{N}$  via 
\[N_x(t;\vec{x})=\#\{i;x_i\leq x\}, \quad \eta_x(t;\vec{x})=\left\{\begin{matrix}
1 &  x_i=x \quad \text{for some }i.\\
0 & \text{otherwise}.
\end{matrix}\right.\]
For the rest of the paper, we prefer to use $N_x(t)$ and $\eta_t(x)$ instead of $N_x(t;\vec{x})$ and $\eta_t(x;\vec{x})$ with the understanding that those are evaluated at $\vec{x}$. We call $N_x(t)$ as the height function of HL-PushTASEP and $\eta_t(x)$ as the occupation variable.
 
  \ed

  \bd[Step Bernoulli Initial Data]\label{eq:InitialData} 
  Fix any positive integers $L$. In step Bernoulli initial condition with parameter $\rho$ over $[0,L]$, each integer lattice site is occupied by a particle with probability $\rho$ independently of others. Whenever $\rho=1$, we refer it as step initial condition.       
  \ed

\subsection{Connection with Stochastic Six Vertex Model}\label{SSVMtoHL-PushTASEP}
Six vertex model was first studied in \cite{Lieb67} to compute the residual entropy of the ice structure. One can think of six vertex model as a lattice model where the configurations are the assignments of six different types of $H_2O$ molecular structure to the vertices of a square grid so that $O$ atoms are at the vertices of the grid. To each $O$ atoms, there are two $H$ atoms attached so that they are at the angles $90^\circ$ or $180^\circ$ to each other. Later, a stochastic version of the six vertex model was introduced by Gwa and Spohn in \cite{GS92}. Recently, several new discoveries related to the stochastic six vertex model came into light in works including \cite{BCG2016,BP16a,BP16b,AB16,A16b}. Furthermore, there are three equivalent ways of describing stochastic six vertex model: $(i)$ as a ferro-electric asymmetric six vertex model on a long rectangle with specific vertex weights (see, \cite{GS92} or \cite[Section~2.2]{BCG2016}); $(ii)$ as an interacting particle system subject to the pushing effect and exclusion of mass (see, \cite[Section~2.2]{BCG2016}, \cite{GS92}); or $(iii)$ as a probability measure on directed path ensemble (see, \cite[Section~1]{BP16a}, \cite[Section~1.1]{AB16}). In what follows, we elaborate in details on the second description of stochastic six vertex model.

Particle dynamics of stochastic six vertex model is governed by a discrete time Markov chain with local interactions.  Fix $b_1,b_2\in [0,1]$. Consider a particle configuration $\mathbb{Y}^{b_1,b_2}(t)=(y_1(t)<y_2(t)<\ldots )$ at time $t$ where $0<b_1,b_2<1$. Then, at time $t+1$, each particle decides with probability $b_1$ whether it will stay in its position, or not. If not, then it moves to the right for $j$ steps with probability $(1-b_2)b_2^{j-1}$ where $j$ can be any integer in between $1$ and the distance of the right neighbor from the particle. Further, the particle at $y_i(t)$ can jump to $y_{i+1}(t)$ with probability $b^{y_{i+1}(t)-y_{i}(t)-1}$. In that case, the particle which was initially at $y_{i+1}(t)$ is pushed towards right by one step and starts to move in the same way.

It has been noted in \cite[Section~2.2]{BCG2016} that the dynamics of the limit
\begin{align}\label{eq:SixVertexToHL-PushTASEP}
Y^{b}(t):=\lim_{\epsilon\to 0}\mathbb{Y}^{1-\epsilon,b}(t\epsilon^{-1})
\end{align}
are those of HL-PushTASEP on $\ZZ_{\geq 0}$. Similarly, they show that the dynamics of \[Z^{p,q}(t):=\lim_{\epsilon\to 0}\{\mathbb{Y}^{\epsilon p,\epsilon q}(t\epsilon^{-1})- t\epsilon^{-1}\}\] are those of ASEP with jump parameters $p$ and $q$ (see also \cite{A16b}). For any finite particle configuration, it is easy to show that the limiting dynamics in \eqref{eq:SixVertexToHL-PushTASEP} indeed matches with HL-PushTASEP on $\ZZ_{\geq 0}$. Heuristically, this is attributed to the facts
that in an time interval $[t,t+\Delta t]$, $(i)$ the expected number of jumps  for each of the particles of the limiting process $Y^b(t)$ varies linearly with $\Delta t$ and $(ii)$ the probability of more than one jumps decays quadratically with $\Delta t$. In particular, $(i)$ shows the number of jumps of each of the particles follows a Poisson point process with intensity $1$ and $(ii)$ shows that the Poisson clocks in all particles are independent of each other. In the case of countable particle configurations, one can complete the proof of the convergence of stochastic six vertex model to HL-PushTASEP on $\ZZ_{\geq 0}$ using the similar arguments as in \cite{A16b}.

\subsection{Main Results}



In the following theorems, we state the asymptotic phase transition result for the height function $N_{\nu t}(t)$ of HL-PushTASEP on $\ZZ_{\geq 0}$ under the step Bernoulli initial data. Let us explain the scalings used in the paper. Assume that initial condition for the HL-PushTASEP dynamics is step Bernoulli with parameter $\rho$. Thus, it is expected that the macroscopic density of the particles varies between $0$ to $\rho$. In fact, it will be shown that whenever $\nu \leq (1-b)^{-1}$, then macroscopic density around $\nu t$ will be $0$. Furthermore, in the case when $(1-b)^{-1}< \nu <(1-\rho)^{-2}(1-b)^{-1}$, macroscopic density varies as $1-(\nu(1-b))^{-1/2}$ around the position $\nu t$. Consequently, one can expect KPZ scaling theory to hold in such scenario. On the contrary, when $\nu>(1-\rho)^{-2}(1-b)^{-1}$, density profile around the position $\nu t$ turns out to be flat, thus yielding the gaussian fluctuation. 
Before stating the theorems, we must define three distributions, namely the GUE Tracy-Widom distribution, square of GOE Tracy-Widom distribution and Gaussian distribution which will arise as limits in three different scenario.

\bd\label{KernelIntroduction}
\be
\ii
Consider a piecewise linear curve $\Gamma^{(1)}$ consists of two linear segments: from $\infty e^{-i\pi/3}$ to $0$ and from $0$ to $\infty e^{i\pi/3}$.
The distribution function $F_{GUE}(x)$ of the GUE Tracy-widom distributions is defined by $F_{GUE}(x)=\mbox{det}(I+K_{Ai})_{L^2(\Gamma^{(1)})}$ where $K_{Ai}$ is the \emph{Airy kernel}. For any $w,w^\prime \in \Gamma^{(1)}$, Airy kernel is expressed as
\begin{align}\label{eq:KAiryDef}
K_{Ai}(w,w^\prime) = \frac{1}{2\pi i}\int_{e^{-2\pi i/3 }\infty}^{e^{2\pi i/3}\infty}  \frac{e^{w^3/3-sw}}{e^{v^3/3-sv}}\frac{1}{(v-w)(v-w^\prime)}dv
\end{align}
where the contour of $v$ is piecewise linear from $-1+e^{-2\pi i/3 }\infty$ to $-1$ to $-1+e^{2\pi i/3 }\infty$. It is important to note that the real part of $v-w$ is always negative.

\ii Let $\delta>0$. Consider a piecewise linear curve $\Gamma^{(2)}$ which is composed of two linear segments: from $-\delta+\infty e^{-i\pi/3}$ to $-\delta$ and from $-\delta$ to $-\delta+\infty e^{i\pi/3}$. Then, the $F^2_{GOE}$ is defined by $F^2_{GOE}(xs)=\mbox{det}(I+K_{GOE,2})_{L^2(\Gamma^{(2)})}$ where $K_{GOE,2}$ is given by
\begin{align}\label{eq:KBBPDef}
K_{GOE,2}(w,w^\prime)=\frac{1}{2\pi i}\int_{-2\delta+\infty e^{-2\pi i/3 }}^{-2\delta+\infty e^{2\pi i/3}}  \frac{e^{w^3/3-sw}}{e^{v^3/3-sv}}\frac{1}{(v-w)(v-w^\prime)}\frac{v}{w} dv
\end{align}
 for aany $w,w^\prime\in \Gamma^{(2)}$. The contour of $v$ is piecewise linear from  $-2\delta+\infty e^{-2\pi i/3 }$ to $-2\delta$ to $-2\delta+\infty e^{2\pi i/3 }$. 

 \ii For completeness, we also add here a very brief description of the Gaussian distribution. Density of the standard Gaussian distribution function $\Phi(x)$ is given by
 \[\frac{1}{\sqrt{2\pi}}\exp\left(-\frac{x^2}{2}\right).\]
\ee
\ed

\bt\label{Asymptheo1Mainstate}
Consider HL-PushTASEP on $\ZZ_{\geq 0}$. Assume initial data is step Bernoulli data. Set $\alpha=\frac{1-\rho}{\rho}$. Let us also fix a real number $\nu>(1-b)^{-1}$. Define the parameters,
\[m_\nu:=\frac{(\sqrt{\nu(1-b)}-1)^2}{1-b},\quad \sigma_\nu:=\frac{b^{1/3}(\sqrt{\nu(1-b)}-1)^{2/3}}{(1-b)^{1/2}\nu^{1/6}}, \quad\varrho := b(\sqrt{\nu(1-b)}-1),\]
\[\tilde{m}_\nu := \frac{\nu}{1+\alpha}-\frac{1}{\alpha(1-b)},\quad \quad \tilde{\sigma}_\nu:=\sqrt{\frac{\nu(1-b)(1-\rho)^2-1}{\alpha^2(1-b)}}.\]

\ba
 \ii If $\rho >1-\left(\nu(1-b)\right)^{-1/2}$, then
 \begin{align}\label{eq:HighDenLimit}
 \lim_{t\to \infty}\mathbb{P}\left(\frac{m_{\nu}t-N_{\nu t}(t)}{\sigma_\nu t^{1/3}}\leq s\right)=F_{GUE}(s).
 \end{align}
 \ii If $\rho= 1-\left(\nu(1-b)\right)^{-1/2}$, then 
 \begin{align}\label{eq:CriticalDenLimit}
 \lim_{t\to \infty}\mathbb{P}\left(\frac{m_{\nu}t-N_{\nu t}(t)}{\sigma_\nu t^{1/3}}\leq s\right)=F^2_{GOE}(s).
 \end{align}
 \ii If $\rho <1-\left(\nu(1-b)\right)^{-1/2}$, then
 \begin{align}\label{eq:LowDenLimit}
\lim_{t\to \infty}\mathbb{P}\left(\frac{\tilde{m}_{\nu}t - N_{\nu t}(t)}{\tilde{\sigma}_\nu t^{1/2}}\leq s\right)= \Phi(s).
 \end{align}
\ee
\et
Let us denote the position of $m$-th particle at time $t$ by $x_{m}(t)$. To this end, one can note that for any $y,t\in \RR^+$ and $m\in \ZZ_{\geq 0}$, the event $\{N_{y}(t)\geq m\}$ is same as the event $\{x_m(t)\leq y\}$. Thus, Theorem~\ref{Asymptheo1Mainstate} translates to the following results (noted down below) in terms of the particle position $x_m(t)$.
\bt\label{Asymptheo2Mainstate}
 Continuing with all the notations introduced in Theorem~\ref{Asymptheo1Mainstate} above, we have
 \ba
 \ii if $\rho >1-\left(\nu(1-b)\right)^{-1/2}$, then 
 \begin{align}\label{eq:HighDenLimit2}
 \lim_{t\to \infty}\mathbb{P}\left(x_{tm_\nu}(t)\leq \nu t+ \left(\frac{\partial m_\nu}{\partial \nu}\right)^{-1}\sigma_\nu s t^{1/3}\right)&=F_{GUE}(s),
 \end{align}
 \ii if $\rho= 1-\left(\nu(1-b)\right)^{-1/2}$, then
 \begin{align}\label{eq:CriticalDenLimit2}
 \lim_{t\to \infty}\mathbb{P}\left(x_{tm_\nu}(t)\leq \nu t+ \left(\frac{\partial m_\nu}{\partial \nu}\right)^{-1}\sigma_\nu s t^{1/3}\right)&=F^2_{GOE}(s),
 \end{align}
 \ii if $\rho <1-\left(\nu(1-b)\right)^{-1/2}$, then
 \begin{align}\label{eq:LowDenLimit2}
\lim_{t\to\infty}\mathbb{P}\left(x_{t\tilde{m}_\nu}(t)\leq t\nu +\left(\frac{\partial \tilde{m}_\nu}{\partial \nu}\right)^{-1}\tilde{\sigma}_\nu st^{1/2}\right)&=\Phi(s).
 \end{align}
\ee
\et

 Using the theorems above, we derive the following law of large numbers for HL-PushTASEP.
 \bc[Law of large numbers]\label{LimitShapesCorollary}
 Consider HL-PushTASEP on $\ZZ_{\geq 0}$ with step initial condition. Denote $\nu_0:=(1-b)^{-1}$. Then,
 \begin{align}\label{eq:Limitshapescorollary}
 \lim_{t\to\infty }\frac{N_{\nu t}(t)}{t}=\left\{\begin{matrix}
 \frac{(\sqrt{\nu(1-b)}-1)^2}{(1-b)} & \quad \text{when }\nu>\nu_0,\\
 0 & \quad \text{otherwise}.
 \end{matrix}\right.
 \end{align}
 \ec

Interactive particle systems which obey KPZ class conjecture (see, Section~\ref{KPZScaling}) form the KPZ universality class. Models is the KPZ universality class are in a sense characterized by the KPZ equation which is written as 
 \begin{align}\label{eq:KPZequation}
 \partial_t\mathcal{H}(t,x)=\frac{1}{2}\delta \partial^2_x\mathcal{H}(t,x)+\frac{1}{2}\kappa(\partial_x\mathcal{H}(t,x))^2+\sqrt{D}\xi(t,x),
\end{align} 
for $\delta,\kappa\in \RR$ and $D>0$ where $\xi(.,.)$ denotes the space time white noise. In a seminal work \cite{KPZ86}, Karder, Parisi and Zhang introduced the KPZ equation as a continuum model for the randomly growing interface which can be subjected to the local dynamics like smoothing, lateral growth and space-time noise. Most prominent feature of the models in KPZ universality class is that the scaling exponents of the fluctuation, the space and the time follow a ratio $1:2:3$. Further, it is expected that the long time asymptotics of any model in the KPZ universality class coincides with the corresponding statistics of the KPZ equation when initial condition of the model corresponds with the initial condition in KPZ equation. To illustrate, consider ASEP on $\ZZ$. Initially, if all the sites in $\ZZ_{\leq 0}$ are occupied and the rest are left empty, then \cite{Tracy2009} shows that the height fluctuation in $t^{1/3}$ scale converges to the Tracy-Widom GUE distribution. Later, it is verified in the work of \cite{Amir11}, \cite{DOT03}, \cite{DR02} and \cite{SS10} that  under narrow wedge initial condition the height fluctuation of the KPZ equation in $t^{1/3}$ scale has the same asymptotic distribution.

 It is now natural to ask what are the other distributions which arise as the long time fluctuation limits of stochastic models in KPZ universality class. There are only very few instances where this question has been addressed. In the case of ASEP under the step Bernoulli initial condition (sites in $\ZZ_{\leq 0}$ are occupied w.p $\rho$ and other places are empty), Tracy and Widom \cite{Tracy2009} demonstrates that the limiting behavior of the height fluctuation undergoes a phase transition. In particular, whenever $\rho$ is close to $1$, scaling exponent of the fluctuations is still $t^{1/3}$ and limiting distribution is Tracy-Widom GUE. On the contrary, when $\rho$ nears $0$, the fluctuation converges to the Gaussian distribution under $t^{1/2}$ scaling. In between those two regime, there lies a critical point where the limiting distribution of the fluctuation is $F^2_1$. Nevertheless, the scaling exponent of the fluctuation at the critical point remains $1/3$. It is shown in \cite{CJ13} that the asymptotic limit of the fluctuation (in ASEP) at the critical point corresponds to the Half-Brownian initial data of the KPZ equation. 
 
 Similar phase transition was in fact first unearthed in the work of Baik, Ben Arous and P\'ech\'e \cite{BBP05}, where they observed a transition phenomenon (now known as BBP transition) in the limiting law of the largest eigenvalue $\lambda_1$ of large spiked covariance matrices. In a follow up paper, \cite{P06} proved the same for the finite rank perturbation of the GUE matrices. When the perturbation has rank one, then the asymptotic distribution of $\lambda_1$ at criticality coincides with $F^2_1$. In the context of interacting particle systems, Barraquand \cite{B15} established the BBP transition in the case of $q-TASEP$ with few slower particles. More recently, Aggarwal and Borodin \cite{AB16} showed that the same phenomenon also holds in the case of stochastic six vertex model when subjected to generalized Bernoulli initial data.

 Our next result is on the convergence of an exponential transform of the height function $N_x(t)$ to a mild solution of the stochastic heat equation (SHE). The stochastic heat equation is written as
\begin{align}\label{eq:SHE}
\partial_t \mathfrak{Z}(t,x)=\frac{1}{2}\delta \partial^2_x\mathfrak{Z}(t,x)+\frac{\kappa}{\delta}\sqrt{D}\mathfrak{Z}(t,x)\xi(t,x)
\end{align}
$\xi(\cdot, \cdot )$ denotes the space time white noise.
 Let us call $\mathfrak{Z}(t,r)\in C([0,\infty ), C(\RR))$ a mild solution of the SHE starting from the initial condition $\mathfrak{Z}^{in}$ if
 \begin{align}\label{eq:SHEDef}
 \mathfrak{Z}(t,r)= \int_{\RR} P_{s}(r-r^\prime)\mathfrak{Z}^{in}(s,r^\prime)dr^\prime+ \int_0^t
 \int_{\RR} P_{t-s}(r-r^\prime)\mathfrak{Z}^{in}(s,r-r^\prime)\xi(ds,dr^\prime)
 \end{align}
 where $P_t(r):=\frac{1}{\sqrt{2\pi t}}\exp(-r^2/2t)$ is the standard heat kernel. For existence, uniqueness, continuity and positivity of the solutions to \eqref{eq:SHEDef}, see \cite[Proposition~2.5]{Corwin12} and \cite{Bertini1997}, \cite{Bertini1995}, \cite{Mueller91}. Interestingly, the Cole - Hopf solution $\mathcal{H}(t,x)$ of the KPZ equation is defined by $\mathcal{H}(t,x)=\frac{\delta}{\kappa}\log \mathfrak{Z}(t,x)$ whenever $\mathfrak{Z}$ is a mild solution (see, \eqref{eq:SHEDef}) of the SHE.
  Bertini and Giacomin \cite{Bertini1997} proved that ASEP under weak asymmetric scaling and near equilibrium initial condition converges to the Cole - Hopf solution of the KPZ equation. Later, the same result has been shown for the narrow wedge initial data in \cite{Amir11}. In both of these works, G\"artner transformation which is a discrete analogue of Cole - Hopf transform plays crucial roles. For  weakly asymmetric exclusion process with hopping range $m$ more than $1$, similar result has first appeared in the work of Dembo and Tsai \cite{Dembo2016}. But, their proof breaks down when $m\geq 4$. In the following discussion, we show that under a particular weak scaling, HL-PushTASEP on $\ZZ_{\geq 0}$ started from the step initial condition converges to the KPZ equation.

 \bd[G\"artner Tranform \& Weak Noise Scaling]\label{GartnerTransform}
  Here, we turn to define explicitly the exponential transformation (or, G\"artner transformation)
 \begin{align}
 Z(t,x)=\left\{
 \begin{matrix}
 b^{N_{ x+\lfloor t/(1-b^{\nu_\epsilon})\rfloor}(t)-(1-{\nu_\epsilon})\left( x+\lfloor t/(1-b^{\nu_\epsilon})\rfloor\right)}\exp(-\mu_\epsilon t)  & \quad \text{when } x+\lfloor t/(1-b^{\nu_\epsilon})\rfloor \in \ZZ_{\geq 0}\nonumber\\
 0 & \quad \text{o.w}
 \end{matrix}\right. \\\label{eq:HopfColeTransform}
 \end{align}
where $t\in \RR_{+}$, $\nu_\epsilon=(\sqrt{5}-1)/2$  and
\begin{align}\label{eq:MuEpsilon}
\mu_\epsilon=\frac{b^{-1}-1}{b^{-\nu_\epsilon}-1}+\frac{(1-\nu_\epsilon)\log b}{1-b^{\nu_\epsilon}}.
\end{align}
 In the rest of our discussion, we will keep on making an abuse of notation by writing $t/(1-b^{\nu_\epsilon})$ instead of $\lfloor t/(1-b^{\nu_\epsilon})\rfloor$. Let us mention that we use a weak noise scaling $b=b_\epsilon:=e^{-\lambda_\epsilon\epsilon^{1/2}}$, where $\epsilon\to 0$ and $\lambda_\epsilon = \nu^{-3/2}_\epsilon$.
  \ed

  We extend the process $Z(t,x)$, defined for $t\in \RR_{+}$ and $x\in \ZZ$, to a continuous process in $\RR_{+}\times \RR$ by linearly interpolating in $x$. Let us introduce the scaled process
\begin{align}\label{eq:ScaledProcess}
Z_\epsilon(t,x):=Z(\epsilon^{-1}t,\epsilon^{-1}x).
\end{align}
Furthermore, we endow the space $C(\RR_+)$, and $C(\RR_+\times \RR)$ with the topology of uniform convergence and use $\Rightarrow$ to denote the weak convergence.
\bd[Near-Equilibrium Condition]\label{NearEquiCon}
Discrete SHE $Z_\epsilon(.,.)$ is defined to have near equilibrium initial condition when the initial data $Z_\epsilon(0,x)$ satisfies the following two conditions:
\begin{align}
\|Z_\epsilon(0,x)\|_{2k}&\leq Ce^{\epsilon\tau |x|}, \label{eq:EquiInitCond1}\\
 \|Z_\epsilon(0,x)-Z_\epsilon(0,x^\prime)\|_{2k}& \leq C  \left(\epsilon |x-x^\prime|\right)^{v}e^{\epsilon\tau(|x|+|x^\prime|)},\label{eq:EquiInitCond2}\\
\end{align}
 for some $v\in (0,1)$ and $\tau>0$.
 \ed
\bt\label{KPZLimitTheo}
Let $\mathfrak{Z}$ be the unique $C(\RR_+\times \RR)$- valued solution of the SHE starting from some initial condition $\Delta$ and let $Z_\epsilon(t,x)\in C(\RR_+\times \RR)$ as defined in \eqref{eq:ScaledProcess} evolved from some near equilibrium initial condition. If $Z_\epsilon(0,.)$ converges to $\mathfrak{Z}(0,.)$ whenever $\epsilon$ tends to $0$, then
\begin{align}\label{eq:WeakConvergence}
Z_\epsilon(.,.)\Rightarrow \mathfrak{Z}(.,.)\quad \text{on }C(\RR_+\times \RR)
\end{align}
as $\epsilon\to 0$.
\et
 It can be noted that the step initial condition, i.e, $x_n(0)=n$ for $n\geq \ZZ_{+}$ and $x_n(0)=-\infty$ for $n\in \ZZ_{-}$, doesn't belong to the family of near equilibrium initial conditions in Definition~~\ref{NearEquiCon}. Nevertheless, the rescaled height function of HL-PushTASEP started with step initial condition also has KPZ limit. We illustrate this further in the following theorem.

 \bt\label{KPZLimitTheoAtStepInit}
  Let $\mathfrak{Z}$ be the unique $C(\RR_+\times \RR)$- valued solution of the SHE starting from delta initial measure. Consider $\tilde{Z}_\epsilon(t,x):= \epsilon^{-1}(1-\exp(-\lambda_\epsilon\nu_\epsilon))Z_\epsilon(t,x)$ in $C(\RR_+\times \RR)$ evolving from the step initial condition. Then,
  \begin{align}
  \tilde{Z}_\epsilon(.,.) &\Rightarrow \mathfrak{Z}(.,.) \quad \text{on }C(\RR_+\times \RR),\label{eq:SHEUnderStep}
  \end{align}
  as $\epsilon$ tends to $0$.
  \et

\bc\label{KPZLimitCoro}
Let $\mathfrak{Z}$ be as in the Theorem \ref{KPZLimitTheo} so that $\mathcal{H}(t,x):=\log \mathfrak{Z}(t,x)$ be the unique Cole-Hopf solution of the KPZ equation starting from the delta initial measure. Then,
\begin{align}
 \log \tilde{Z}_\epsilon(.,.)&\Rightarrow \mathcal{H}(.,.) \quad \text{on } C(\RR_+\times \RR)\label{eq:KPZUnderStep}
\end{align}
as $\epsilon\to 0$.
\ec

 \textsc{Proof of Theorem~\ref{KPZLimitTheoAtStepInit}.}
  To begin with, note that $\epsilon\sum_{\zeta\in \Xi(0)}\tilde{Z}_\epsilon(0,\zeta)=1$. Thus, $\tilde{Z}_\epsilon(0,.)$ converges to the delta initial measure as $\epsilon$ goes to $0$. Below, in Proposition~\ref{MomentEstimateProp} which we prove in Section~\ref{WeakScaleLim}, we show the moment estimates for $\tilde{Z}_\epsilon$ at any future time point $t$. Thus, using Theorem~\ref{KPZLimitTheo}, one can complete the proof in the same way as in \citep[Section~2]{Amir11}.

 \bp\label{MomentEstimateProp}
 Fix any $T>0$. Consider $\tilde{Z}(t,x)=\epsilon^{-1}(1-\exp(-\lambda_\epsilon\nu_\epsilon))Z(t,x)$ evolving from the step initial condition. Then, for any $t\in (0,\epsilon^{-1}T]$, $\zeta_1,\zeta_2\in \Xi(t)$ and $v\in (0,1/2)$, we have
 \begin{align}
 \Vert \tilde{Z}(t,\zeta)\Vert &\leq C\min\{\epsilon^{-1/2},(\epsilon t)^{-1/2}\},\label{eq:MomemtEstimateFromStep1}
\\ \Vert \tilde{Z}(t,\zeta)-\tilde{Z}(t,\zeta^\prime)\Vert &\leq (\epsilon|\zeta-\zeta^\prime|)^v (\epsilon t)^{-(1+v)/2}.\label{eq:MomemtEstimateFromStep2}
 \end{align}
where the constant $C$ depend only on $T$ and $v$.
 \ep

\subsection{KPZ Scaling Theory}\label{KPZScaling}
  
  Here, we explain how the asymptotic fluctuation results in Theorem~\ref{Asymptheo1Mainstate} confirms the KPZ scaling theory (see, \cite{KPH92}, \cite{Spohn12}) of the physics literature. To start with, we present the predictions of KPZ scaling theory in the context of exclusion process following \cite{Spohn12}. For that, one need to assume that spatially ergodic and time stationary measures of the exclusion process on $\ZZ$ are precisely labelled by the density of the particles $\rho$, where
\begin{align}\label{eq}
\rho:= \lim_{n\to \infty}\frac{1}{2n+1}\sum_{k=-n}^n \eta(n).
\end{align}
We denote the stationary measure arising out of the density $\rho$ by $\nu_\rho$. Consider the average steady state current $j(\rho)$ which counts the average number of particles transported from $0$ to $1$ in unit time where the system is distributed as $\nu_\rho$. Further, denote the integrated covariance $\sum_{j\in \ZZ}\mbox{Cov}(\eta_0,\eta_j)$ by $A(\rho)$ where the covariance between the occupation variables $\eta_0$ and $\eta_j$ is computed under the stationary measure $\nu_\rho$. One expects that on the microscopic scale the rescaled particle density $\varrho(x,t)$ given as
\begin{align}\label{eq:RescaledDensity}
\varrho(x,\tau)= \lim_{t\to \infty}\mathbb{P}\left(\text{there is a particle at site }\lfloor xt\rfloor\text{ at time }\lfloor \tau t\rfloor\right)
\end{align}
satisfies the conservation equation
\begin{align}\label{eq:ConservPartialDiff}
\frac{\partial }{\partial \tau}\varrho(x,\tau)+\frac{\partial }{\partial x} j(\varrho(x,\tau))=0
\end{align}
when the system is started with the step initial condition, i.e,
\begin{align}\label{eq:StepInitCond}
\varrho(x,0)=\left\{\begin{matrix}
1 & x\geq 0\\
0 & x< 0.
\end{matrix}\right.
\end{align}
This result can also be phrased as a hydrodynamic limit theory. For any exclusion process, consider the height function $h(.,.):\ZZ\times \RR_{+}\to \ZZ$ given by
\begin{align}\label{eq:HeightFunctionDef}
h(j,t):=\left\{\begin{matrix}
N_t +\sum_{i=1}^j \eta_t(j) & \text{ if }j> 0\\
N_t & \text{ if } j=0 \\
N_t - \sum_{i=1}^{-j} \eta_t(-i) & \text{ if } j<0.
\end{matrix}\right.
\end{align}
Then, using the arguments in \cite[Part~II, Section~3.3]{spohn2012large}, one can prove that \eqref{eq:ConservPartialDiff} implies the following law of large numbers
\begin{align}\label{eq:LawOfLaregNumbersExclussionProcess}
\lim_{t\to \infty}\frac{h(\nu t,t)}{t} = \phi(\nu)
\end{align}
where the limit shape $\phi(.)$ is given by
\begin{align}\label{eq:GeneralLimitShapesForm}
\phi(y)= \sup_{\rho\in[0, 1]}\{y\rho -j(\rho)\}.
\end{align}
Furthermore, Let us denote $\lambda(\rho)= -j^{\prime\prime}(\rho)$. Under such parametrization, we have the following KPZ class conjecture.

\textbf{KPZ Class Conjecture:}(\cite{Spohn12}) Let $y$ be such that $\phi$ is twice differentiable at $y$ with $\phi^{\prime\prime}\neq 0$. Set $\rho_0=\phi^\prime(y)$. If $0\leq \rho_0<1$, $A(\rho_0)<\infty$ and $\lambda(\rho_0)\neq 0$, then under the step initial condition in \eqref{eq:StepInitCond}, one has
\begin{align}\label{eq:KPZclassConjecture}
\lim_{t\to \infty}\mathbb{P}\left(\frac{t\phi(y)-h(\lfloor yt\rfloor ,t)}{(-\frac{1}{2}\lambda A^2)^{1/3}t^{1/3}}\leq s\right)= F_{GUE}(s).
\end{align}
In the case when $\phi(y)$ is linear or cusp like, then KPZ conjecture doesn't hold. It is likely that the scaling of the fluctuations would be different in those cases. Historically, the conjecture on the scaling exponent $t^{1/3}$ of the fluctuation above dates back to the the work \cite{KPZ86}. To that end, the magnitude of the fluctuation was obtained in \cite{KPH92} and the limit distribution was first surfaced in
\cite{Kurt00}.

\bl\label{StationaryMeasureLemma}
 Fix $0\leq \rho\leq 1$. Consider the HL-PushTASEP model on $\ZZ_{\geq 0}$. Further, assume that a jet of particles is entering into the region $\ZZ_{\geq 0}$ through $0$ at rate $\tau$ where $\tau:= \rho(1-\rho)^{-1}(1-b)^{-1}$. Then, the product Bernoulli measure $\nu_\rho$ will be the invariant measure for the HL-PushTASEP model. In particular, each of the occupation variables $\{\eta_t(x)|x\in \ZZ_{\geq 0}\}$ follows $\mbox{Ber}(\rho)$ independent of others.
\el
\begin{proof}
Assume, at time $t=0$, each site on $\ZZ_{+}$ is occupied by at most one particle with probability $\rho$ independently of others. First, we show $\nu_\rho(\eta_t(0)=1)=\rho$ for any $t\geq 0$. Note that influx rate of particles at site $x=0$ is equal to $\tau$. Then, rate at which particle leaves site $0$ is given by $\chi_1+\chi_2+\chi_3$
\begin{align}\label{eq:OutfluxRate}
\chi_1:= \tau \rho, \quad \chi_2:=\tau b(1-\rho),\quad \text{and }, \quad \chi_3:= \rho.
\end{align}
To begin with, consider the case when initially the site $0$ was occupied. Probability of this event is $\rho$. In such scenario, $\chi_1$ captures the rate of out flux of any particle to the right from the site $0$ after a particle arrives there. On the contrary, probability that the site $0$ was initially unoccupied is $1-\rho$. Thus, if a particle arrives when the site $0$ is empty, it moves further one step towards the right with probability $b$. Henceforth, the value of $\chi_2$ synchronizes perfectly with the rate of outflux from the site $0$ whenever initially it was unoccupied. Finally, the particle which was initially at the site $0$ jumps to the right at rate $1$ and this contribution has added up in the total rate out flux through $\chi_3$. Due to the specific choice of $\tau$, we get $\tau= \chi_1+\chi_2+\chi_3$. Consequently, the rate of in flux matches with the rate of out flux of the particles from the site $0$. Therefore, probability that the site $0$ is occupied doesn't change over time. More importantly, the rate at which particles will move into the site $1$ is also equal to $\tau$. This further implies $\mathbb{P}(\eta_t(1)=1)=\rho$ for any $t\geq 0$. Using induction, one can now prove that the occupation variable $\eta_t(x)$ for any $x\in \ZZ_{+}$ and $t\geq 0$ follows $\mbox{Ber}(\rho)$. In fact, using similar argument, one can show that at any time $t$, the rate of influx of the particles into any interval $[x_1,x_2]$ of finite length over $\ZZ_{+}$ is equal to the rate at which the particles will come out of the interval $[x_1,x_2]$.

Now, we turn to prove the independence of any two occupation variables $\eta_t(x_1)$ and $\eta_t(x_2)$ for $x_1<x_2\in \ZZ_{+}$. Let us consider $(\epsilon_{x_1},\epsilon_{x_{1}+1},\ldots ,\epsilon_{x_2})\in \{0,1\}^{x_2-x_1+1}$. Further, denote
\begin{align}\label{eq:EpsilonEta}
\eta^{\epsilon_a}_t(a):=\left\{\begin{matrix}
\eta_t(a) & \quad \text{if }\epsilon_a=1,\\
1-\eta_t(a) & \quad \text{if }\epsilon_a=0.
\end{matrix}\right.
\end{align}
To this end, one can note that
\begin{align}
\mathbb{E}\left(\prod_{a=x_1}^{x_2}\eta^{\epsilon_a}_t(a)\right)=\mathbb{E}\left(\prod_{a=x_1}^{x_2}\eta^{\epsilon_a}_0(a)\right)
\end{align}
for any choice of the variables $(\epsilon_{x_1},\epsilon_{x_{1}+1},\ldots ,\epsilon_{x_2})$ in the space $\{0,1\}^{x_2-x_1+1}$. This is again due to the principle of conservation of mass, i.e, the influx rate into the interval $[x_1,x_2]$ is same as the rate of out flux of the particles. Thus, we get
\begin{align}\label{eq:IndependenceAssertion}
\mathbb{E}(\eta_t(x_1)\eta_t(x_2))&=\mathbb{E}\left(\sum_{\epsilon_{x_1+1},\ldots ,\epsilon_{x_2-1}}\eta_t(x_1)\eta_t(x_2)\prod_{a=x_1+1}^{x_2-1}\eta^{\epsilon_a}_{t}(a)\right) \\
&=\mathbb{E}\left(\sum_{\epsilon_{x_1+1},\ldots ,\epsilon_{x_2-1}}\eta_0(x_1)\eta_0(x_2)\prod_{a=x_1+1}^{x_2-1}\eta^{\epsilon_a}_{0}(a)\right)
 &=\mathbb{E}\left(\eta_0(x_1)\eta_0(x_2)\right)
 &=\rho^2.
\end{align}
To that effect, covariance of $\eta_t(x_1)$ and $\eta_t(x_2)$ becomes $0$. As these are Bernoulli random variables, thus the independence is proved.
\end{proof}

Now, we turn to showing how Theorem~\ref{Asymptheo1Mainstate} verifies the KPZ class conjecture in the context of HL-PushTASEP. In Lemma~\ref{StationaryMeasureLemma}, it is proved that all the stationary translation invariant measure of HL-PushTASEP are given by product measure $\nu_\rho$ where each site follows $\mbox{Ber}(\rho)$ distribution. Further, it has been also shown that under stationary measure $\nu_\rho$ steady state current will be $\rho(1-\rho)^{-1}(1-b)^{-1}$. To that effect, we have
\begin{align}\label{eq:LimitShapesInHL-PushTASEP}
\phi(y)=\sup_{\rho\in[0,1]}\left\{\rho y-\frac{\rho}{(1-\rho)(1-b)}\right\}.
\end{align}
If $y(1-b)>1$, then the function $g(\rho)=y\rho -\rho(1-\rho)^{-1}(1-b)^{-1}$ is maximized when $\rho=1-(y(1-b))^{-1/2}$. In this scenario, we get
\begin{align}\label{eq:LimitShapesSpecification}
\phi(y)=\frac{\left(\sqrt{y(1-b)}-1\right)^2}{(1-b)}.
\end{align}
One can further derive that the equilibrium density $\rho$ corresponding to $y$ is $1-(y(1-b))^{-1/2}$. To this end, we have $\lambda(\rho)=-j^{\prime\prime}(\rho)=  -2y^{3/2}b(1-b)^{1/2}$. Further, stationarity of the product Bernoulli measure implies $A(\rho)=\rho(1-\rho)=(y(1-b))^{-1/2}(1-(y(1-b))^{-1/2})$. Finally, it can be noted that
\begin{align}\label{eq:SigmaDeriv}
\left(-\frac{1}{2}\lambda(\rho)A^2(\rho)\right)^{1/3}=\frac{b^{1/3}\left(\sqrt{y(1-b)}-1\right)^{2/3}}{y^{1/6}(1-b)^{1/2}}=\sigma_y.
\end{align}
   Thus, KPZ class conjecture implies that the height function $N_{\lfloor yt\rfloor }(t)$ of HL-PushTASEP started from the step initial condition has a limit shape given by $t\phi(y)$ in \eqref{eq:LimitShapesInHL-PushTASEP} and more importantly, the fluctuation around the limit shape under the scaling of $(-2^{-1}\lambda A^2 t)^{1/3}$ converges to the Tracy-Widom GUE distribution. As one can see, Theorem~\ref{Asymptheo1Mainstate} verifies the conjecture and at the same time, identifies the form of the limit shape and the correct scaling of the fluctuation.

\section{Transition Matrix \& Laplace Transform}\label{TransitionMatrixLaplaceTransfom}

For proving the asymptotic results, we need to have concrete knowledge on the distribution of the positions of the particles at any fixed time. In most of the cases of interactive particle systems (see, \citep{TracyWidom11} for ASEP, \citep{BCG2016} for the stochastic six vertex model), transition matrices for the underlying Markov processes play important roles towards that goal. In particular, transition matrices assists in getting suitable Fredholm determinant formulas for the Laplace transforms of some of the important observables like in the case of several quantum integrable interacting particle systems (see, \cite{BorodinCorwinMac} for a survey). In the following subsections, we elaborate on the derivation of the transition matrix of the HL-PushTASEP and its use to
proclaim the formulas for the Laplace transforms.


\subsection{Transition Matrix \& Transition Probabilities}

Primary goal of this section is to compute the transition matrix in the case of finite particle configurations of the HL-PushTASEP. Later in this section, we use the transition matrix to specialize on the transition probabilities of a single particle in a system of finite particle configuration. At first, fix some positive integer $N$. For any $\mathcal{A}\subseteq X^N$ define  \[P_{\vec{x}}(\mathcal{A};t):=\sum_{\vec{x}\in \mathcal{A}}\mathcal{T}^{(N)}_t(\vec{x}\to \vec{y})\]
where $\mathcal{T}^{(N)}_t(\vec{x}\to \vec{y})$ denotes the probability of transition to $\vec{y}\in X^N$ from $\vec{x}\in X^N$. In the following proposition, we derive the eigenfunctions of the transfer matrix $\mathcal{T}^{(N)}_t$.
\bp\label{EigenFunctionLemma}
For $N>0$, fix $N$ small complex numbers $z_1,z_2,\ldots ,z_N$ such that $|bz_i|<1$ for $1\leq i\leq N$ and $1-(2+b)z_j+bz_iz_j\neq 0$ for all $1\leq i\neq j\leq N$. Fix a permutation $\sigma\in \mathfrak{S}(N)$. Define,
\[A_\sigma=(-1)^\sigma\prod_{i<j}\frac{1-(1+b)z_{\sigma(i)}+bz_{\sigma(i)}z_{\sigma(j)}}{1-(1+b)z_i+bz_iz_j}.\]
Then the function
\[\Phi(x_1,x_2,\ldots ,x_N;z_1,z_2,\ldots ,z_N)=\sum_{\sigma\in \mathfrak{S}(n)}A_\sigma\prod_{i=1}^N z_{\sigma(i)}^{x_i}\]
is an eigenfunction of the transfer matrix $\mathcal{T}^{(N)}$, i.e
 \begin{align}\label{eq:eigenfunc2}
 \sum_{\vec{y}\in X^N}\mathcal{T}^{(N)}_t(\vec{x}\to \vec{y})&\Phi(\vec{y};z_1,z_2,\ldots ,z_n)=\exp\left(\sum_{i=1}^N -t\frac{1-z_i}{1-bz_i}\right)\Phi(\vec{y};z_1,z_2,\ldots ,z_n).
\end{align}
\ep

As pointed out in Section~\ref{SSVMtoHL-PushTASEP}, HL-PushTASEP can be considered as a continuum version of the stochastic six vertex model for which the transfer matrix has been worked in full details in the past. See \cite{Lieb67}, \cite{Nolden92} for the derivation in the case of the stochastic six vertex model under periodic boundary condition and \cite[Theorem~3.4]{BCG2016} for the derivation on $\ZZ$.  
Let us note that one can derive the relation \eqref{eq:eigenfunc2} by taking appropriate limit (see, \eqref{eq:SixVertexToHL-PushTASEP}) of the transfer matrix of stochastic six vertex model in \cite[Eq.16]{BCG2016}. Although, one has to be careful because the limit might not exist. For completeness, we present a self contained proof of the result in Appendix~\ref{suppA}.

\bp\label{Ntransition}
For $N>0$, fix $N$ small complex numbers $z_1,z_2,\ldots ,z_N$ such that $|bz_i|<1$ for $1\leq i\leq N$ and $1-(2+b)z_j+bz_iz_j\neq 0$ for all $1\leq i\neq j\leq N$. For any permutation $\sigma\in \mathfrak{S}(N)$, define
\[A^\prime_\sigma=(-1)^\sigma\prod_{i<j}\frac{1-(1+b)z_{\sigma(j)}+bz_{\sigma(i)}z_{\sigma(j)}}{1-(1+b)z_j+bz_iz_j}\]
and
\[A^{\prime\prime}_\sigma=(-1)^\sigma\prod_{i<j}\frac{1-(1+b^{-1})z_{\sigma(i)}+b^{-1}z_{\sigma(i)}z_{\sigma(j)}}{1-(1+b^{-1})z_i+b^{-1}z_iz_j}.\]
Let us assume $\vec{x}=(x_1,\ldots ,x_N)$ and $\vec{y}=(y_1,\ldots ,y_N)$. Denote the transition matrix for $N$ particles in time $t$ by $\mathcal{T}^{(N)}_t$.
Then, we have
\begin{align}\label{eq:EigenInversion1}
\mathcal{T}^{(N)}_t(\vec{x}\to \vec{y})=\int_{(C_r)^N}\sum_{\sigma\in \mathfrak{S}(N)}A^\prime_\sigma \prod_{i=1}^n z_{\sigma(i)}^{-x_i}z_i^{y_i-1}\exp\left(-t\frac{1-z_i}{1-bz_i}\right)dz_1\ldots dz_N
\end{align}
and
\begin{align}\label{eq:EigenInversion2}
\mathcal{T}^{(N)}_t(\vec{x}\to \vec{y})=\int_{(C_R)^N}\sum_{\sigma\in \mathfrak{S}(N)}A^{\prime\prime}_\sigma \prod_{i=1}^n z_{\sigma(i)}^{x_i}z_i^{-y_i-1}\exp\left(-t\frac{1-z^{-1}_i}{1-bz^{-1}_i}\right)dz_1\ldots dz_N
\end{align}
where $C_r$ is small positively oriented circle leaving all the singularities outside and $C_R$ is large positively oriented circle containing all the singularities inside.
\ep
\begin{proof}
Using the relation \eqref{eq:eigenfunc2}, both the formulas in \eqref{eq:EigenInversion1} and \eqref{eq:EigenInversion2} can be proved in the same way as in \cite[Theorem 3.6]{BCG2016} (see also \cite[Section II.4.b]{TracyWidom11}). We omit here further details of the proof.
\end{proof}

\bt\label{TransitProb}
 Fix an integer $N>0$ and a positive real number $t$. Let us assume that we have $N$ particles which are initially at the position $\vec{y}=(y_1,y_2,\ldots ,y_N)\in X^N$ evolve according the dynamics of HL-PushTASEP. Then, for $1\leq m\leq N$, we have
 \begin{align}\label{eq:1}
 P_{\vec{y}}(x_m=x;t)&=(-1)^{m-1}b^{m(m-1)/2}\sum_{m\leq k\leq N}\sum_{|S|=k}b^{\kappa(S,\ZZ_{>0})-mk -k(k-1)/2}\binom{k-1}{m-1}_{b}\nonumber\\
 &\times \frac{1}{(2\pi i)^k}\oint\ldots \oint \prod_{i,j\in S, i<j}\frac{z_j-z_i}{1-(1+b^{-1})z_i+b^{-1}z_iz_j}\nonumber\\
 &\times \frac{1-\prod_{i\in S}z_i}{\prod_{i\in S}(1-z_i)}\prod_{i\in S}z_i^{x-y_i-1}\exp\left(-t\frac{1-z_i^{-1}}{1-bz_i^{-1}}\right)dz_i
 \end{align}
 where the summation goes over $S\subset \{1,2,\ldots ,N\}$ of size $k$, $\kappa(S,\ZZ_{>0})$ is the sum of the elements in $S$, and contours are positively oriented large circles of equal radius which contains all the singularities.
\et
\begin{proof}
 This result is in the same spirit of \cite[Theorem 4.9]{BCG2016} where they found out explicitly the probability of the similar event in the case of the stochastic six vertex model.  They closely followed the techniques used in \cite[Section 6]{Tracy2008} in the context of N-particle ASEP. It begins with the realization that $\mathbb{P}(x_m(t)=x;1,b)$ is the sum of the probabilities $\mathcal{T}^{(N)}_t(\vec{y}\to \vec{x})$ over all $\vec{x}=(x_1,x_2,\ldots ,x_N)\in X^N$ such that $-\infty <x_1<x_2<\ldots <x_{m-1}<x_m=x$ and $x<x_{m+1}<\ldots <x_N<\infty$. Now, one have to make use of the formulas for $\mathcal{T}^{(N)}_t(\vec{y}\to \vec{x})$ mentioned in \eqref{eq:EigenInversion1} and \eqref{eq:EigenInversion2}. Notice that for the sum over all $\vec{x}\in X^N$ such that $-\infty <x_1<x_2<\ldots <x_{m-1}<x_m=x$ contours of integration is $C_r$ and for other part of the sum contours will be
 $C_R$. For brevity, we are skipping here the details of the proof.
 \end{proof}

\subsection{Moments Formulas}
In this subsection, we aim to present a series of results on moments of the height function of  HL-PushTASEP. One can see \cite[Section 4.4]{BCG2016} for similar results in the case of the stochastic six vertex model. In Section~\ref{secAsymptotics}, these formulas turns out to be instrumental for obtaining the asymptotics.

For a function $f:X^\NN\to \RR$, let $\mathbb{E}_{\vec{y}}(f;t)$ represent the expectation of $f$ at time $t$, i.e
\[\mathbb{E}_{\vec{y}}(f;t)=\sum_{\vec{x}\in X^\NN}\mathcal{T}^{(N)}_t(\vec{y}\to \vec{x}).\]
Also, let us introduce $(a;q)_k$ and $(a;q)_{\infty}$ by defining
\begin{align}\label{eq:QFac}
(a;q)_k:=\prod_{j=1}^k (1-aq^{j-1}) \quad \text{and }\quad (a;q)_k:=\prod_{j=1}^\infty (1-aq^{j-1}).
\end{align}
Furthermore, there is a $q$ - analogue of binomials defined as follows
\begin{align}\label{eq:QBinom}
\binom{n}{k}_{q}:=\frac{(1-q^N)(1-q^{N-1})\ldots (1-q^{N-k+1})}{(1-q)(1-q^2)\ldots (1-q^k)}=\frac{(q^{N-k+1};q)_k}{(q;q)_k}.
\end{align}
In both of the above definitions, we assume $k$ is a non-negative integers.
\bp\label{MGF}
Consider HL-PushTASEP on $\ZZ_{\geq 0}$. Fix a positive real number $t$. Assume that the process has been started from the step Bernoulli initial condition with parameter $\rho$. Denote the initial data by $\mbox{stepb}$ which belongs to $X^\NN$. Then for $L=0,1,2,\ldots $ and any positive integer $x$, we have
\begin{align}\label{eq:5}
 \mathbb{E}_{\mbox{stepb}}(b^{LN_x};t)&=1+(b^{-L}-1)\sum_{k=1}^L\left(\prod_{i=1}^{k-2}(1-b^{1-k+L}b^i)\right)\frac{b^{-k(k-1)/2+L-k}\rho^k}{(b^{-1};b^{-1})_k}\nonumber\\
 &\times \frac{1}{(2\pi i)^k}\oint\ldots \oint \prod_{1\leq i<j\leq k}\frac{z_j-z_i}{1-(1+b^{-1})z_i+b^{-1}z_iz_j}\nonumber\\
 &\times \prod_{i=1}^k\frac{z_i^x(b^{-1}-1)}{(\rho+(1-\rho)b^{-1}-z_ib^{-1})(1-z_i)}\exp\left(-t\frac{1-z_i^{-1}}{1-bz_i^{-1}}\right)dz_i
 \end{align}
 where the contours of integrations are largely oriented circles with equal radius and contain all the singularities of the integrand. 
\ep
\begin{proof}
One can find a similar result in \cite[Proposition 4.11]{BCG2016} in the case of stochastic six vertex model with step initial data. For the sake of clarity, we present an independent proof in ~\ref{suppC}. 
\end{proof}

\bp\label{MGFSimplified}
Consider HL-PushTASEP on $\ZZ_{\geq 0}$. Fix a positive number $t$. As usual $\mbox{stepb}\in X^\NN$ denotes the step Bernoulli initial condition  with parameter $\rho$. Then for $L=0,1,2,\ldots$ and for any integer $x$, we have
\begin{align*}
\mathbb{E}_{\mbox{stepb}}(b^{LN_x};t)=\frac{b^{L(L-1)/2}}{(2\pi i)^L}\oint \ldots \oint \prod_{A<B}\frac{u_A-u_B}{u_A-bu_B}\prod_{i=1}^L\exp\left(\frac{tu_i}{b}\right)\left(\frac{1+u_i}{1+u_ib^{-1}}\right)^x\frac{du_i}{u_i(1-\alpha u_ib^{-1})},
\end{align*}
where the positively oriented contour for $u_A$ includes $0,-b$, but does not include $b\rho(1-\rho)^{-1}$ or $\{bu_B\}_{B>A}$. Note that $\alpha:=(1-\rho)\rho^{-1}$.
\ep
\begin{proof}
To prove this result, one needs to substitute $z_i=(1+u_i)/(1+u_ib^{-1})$ into the formula \eqref{eq:5}.  Consequently, one can notice
\[\frac{z_j-z_i}{1-(1+b^{-1})z_i+b^{-1}z_iz_j}=\frac{u_i-u_j}{u_i-bu_j}\quad ,\quad \frac{\rho(b^{-1}-1) dz_i}{(\rho+(1-\rho)b^{-1}-z_ib^{-1})(1-z_i)}=\frac{du_i}{u_i(1 - \alpha u_ib^{-1})}\]
where $\alpha=(1-\rho)/\rho$. Let us mention that a similar result is proved in the case of stochastic six vertex model in \cite[Theorem~4.12]{BCG2016}. See \cite[Theorem~4.20]{BCS14} for a different proof in the case of ASEP. Rest of the calculation after the substitution can be completed in the same way as in \cite[Theorem~4.12]{BCG2016}. We omit any further details from here.

\subsection{Fredholm Determinant Formula}
In this section, we discuss the representation of the moment formulas of HL-PushTASEP in terms of Fredholm determinant. In the last ten years, a large number of works had been put forward in the literature of KPZ universality class and Fredholm determinant formulas play a crucial role in reaching out such results. For similar exposition in stochastic six vertex model, see \cite[Section~4.5]{BCG2016}.
See \cite[Section 3.2]{BorodinCorwinMac} and \cite[Section 3 and 5]{BCS14} for corresponding formulas in the case of $q$-Whittaker process and $q$-TASEP. Moreover, various other instances can be found in the works like \citep[Theorem~4.2]{Corwin2016}, \cite[Section~3.3]{barraq2016}, \cite{BCF12}. Definition and some of the useful properties of the Fredholm determinants are discussed in Appendix~\ref{App:ApendixA}.

\bd
Here, we define the contours $D_{R,d,\delta}$ and $D_{R,d,\delta;\kappa}$. Let us fix positive numbers $R,\delta,d$ such that $d,\delta<1$. Then the contour $D_{R,d,\delta}$ is  composed of five linear sections:
 \begin{align}\label{eq:DRdel}
 D_{R,d,\delta}:= (R-i\infty,R-id]\cup [R-id,\delta -id]\cup [\delta-id,\delta+id]\cup [\delta+id,R+id]\cup[R+id,R+i\infty).
\end{align}
Let $\kappa$ be an integer which is greater than $R$. Denote two points at which circle of radius $\kappa+\delta$ and centre at $0$ crosses $(R-i\infty,R-id]$ and $(R-i\infty,R-id]$ by $\zeta$ and $\zeta^\prime$. Let us call that minor arc joining $\xi_\kappa$ and $\xi^\prime_\kappa$ to be $I_\kappa$. Then $D_{R,\delta,d;\kappa}$ is defined to be a positively oriented closed contour formed by the union
 \begin{align}\label{eq:DRdelk}
 D_{R,d,\delta}:= (\xi_\kappa,R-id]\cup [R-id,\delta -id]\cup [\delta-id,\delta+id]\cup [\delta+id,R+id]\cup[R+id,\xi^\prime_\kappa)\cup I_\kappa.
\end{align}
See Figure~\ref{fig:3} for details. 
\ed

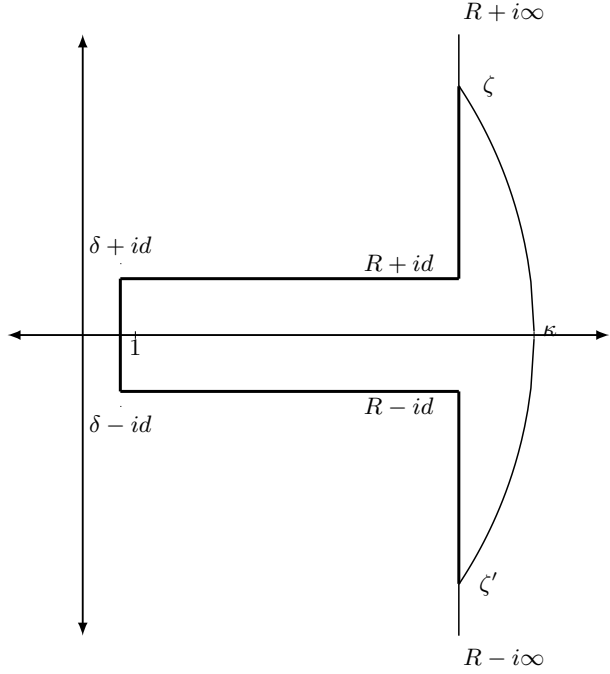
\begin{figure}
\begin{center}
\begin{tikzpicture}
  \draw[line width = 0.25mm, <->] (-1,0) -- (7,0) node {};
  \draw[line width = 0.25mm, <->] (0,-4) -- (0,4) node {};
  \draw[line width = 0.4mm, -] (0.5,-0.75) -- (0.5, 0.75) node {};
  \draw[line width = 0.4mm, -] (0.5, -0.75) -- (5,-0.75) node {};
  \draw[line width = 0.4mm, -] (0.5, 0.75) -- (5, 0.75) node{};
  \draw[line width = 0.4mm, -] (5,0.75) -- (5, 3.31) node {};
  \draw[line width = 0.4mm, -] (5,0.-0.75) -- (5, -3.31) node {};
  \draw[line width = 0.2mm, -] (5,-3.31) -- (5, -4) node {};
  \draw[line width = 0.2mm, -] (5,3.31) -- (5, 4) node {};
  \draw[] (0.5,0.95) -- (0.5,0.95) node[above] {$\delta+id$};
  \draw[] (0.5,-0.95) -- (0.5,-0.95) node[below] {$\delta - id$};
  \node[] at (4.2,-0.95) (a) {$R-id$};
  \node[] at (4.2,0.95) (b) {$R+id$};
  \node[] at (5.6,4.3) (c) {$R+i\infty$};
  \node[] at (5.6,-4.3) (d) {$R-i\infty$};
  \node[] at (5.4,3.31) (e) {$\zeta$};
  \node[] at (5.4,-3.31) (f) {$\zeta^\prime$};
  \draw[-] (6,-0.05) -- (6,0.05) node[right] {$\kappa$};
   \draw[-] (0.7,-0.05) -- (0.7,0.05) node[below] {$1$};
   \draw[domain=5:6, variable=\x, line width = 0.2mm] plot ({\x},{sqrt(36-\x^2)}) ;
   \draw[domain=5:6, variable=\x, line width = 0.2mm] plot ({\x},{-sqrt(36-\x^2)}) ;
\end{tikzpicture}
\caption{Plot of the contour $D_{R,d,\delta}$ and $D_{R,d,\delta;\kappa}$.}
\label{fig:3}
\end{center}
\end{figure}

\bt\label{FredRep}
Consider HL-PushTASEP on $\ZZ_{\geq 0}$. Fix a positive number $t$. Let $\mbox{stepb}\in \mathbb{W}_{\mathbb{N}}$ stands for the step Bernoulli initial data with parameter $\rho$. Let $C_r$ be a positively oriented circle centered at some negative real number $r$ such that $C_r$ include $0$ and $-b$, but does not contain $b\alpha^{-1}$. Furthermore, $C_r$ does  not intersect with $bC_r$, i.e, $bC_r$ is contained in the interior of $C_r$ . Then, there exist a sufficiently large real number $R>1$, and sufficiently small real numbers $\delta, d<1$ such that for all $\zeta \in \mathbb{C}\backslash \RR_{\geq 0}$, we have
\begin{align}\label{eq:UniformBoundCond}
(A)\p\p\inf_{\substack{w,w^\prime\in C_r\\ \kappa\in (2R,\infty)\cap \ZZ\\s\in D_{R,d,\delta;\kappa}}}|b^sw-w^\prime|>0, \quad\quad (B)\p\p\sup_{\substack{w\in C_r\\ \kappa\in (2R,\infty)\cap\ZZ\\s\in D_{R,d,\delta;\kappa}}}\frac{g(w)}{g(b^s w)}<\infty.
\end{align}
Moreover, for any $\zeta\in \mathbb{C}\backslash \RR_{\geq 0}$, we have
\begin{align}\label{eq:qLaplaceTransform}
\mathbb{E}_{\mbox{stepb}}\left(\frac{1}{(\zeta b^{N_x};b)_\infty};t\right)=\det\left(I+K^{b}_{\zeta}\right)_{L^2(C_r)},
\end{align}
where the kernel $K^{b}_\zeta$ is given by
\begin{align}\label{eq:KernelIntro}
K^{b}_\zeta(w,w^\prime)=\frac{1}{2\pi i}\int_{D_{R,d,\delta}}\frac{(-\zeta)^s}{\sin(\pi s)}\cdot\frac{g(w;1,b,x,t)}{g(b^sw;1,b,x,t)}\cdot\frac{ds}{b^sw-w^\prime},
\end{align}
 integration contour $D_{R,d,\delta}$ is oriented from bottom to top and
\begin{align}\label{eq:KernelPrinceComp}
g(z;1,b,x,t)=\left(\frac{1}{1+zb^{-1}}\right)^x\exp\left(\frac{tz}{b(1-b)}\right)\frac{1}{(\alpha b^{-1}z;b)_{\infty}}.
\end{align}
where $\alpha=\rho^{-1}(1-\rho)$.
\et
\begin{proof}
We first prove the existence of $R,d$ and $\delta$ such that $(A)$ and $(B)$ in \eqref{eq:UniformBoundCond} are satisfied. Note that $C_r$ doesn't contain $0$ on it and it is a bounded closed set in $\CC$, thus, compact. Hence, $|w^\prime w^{-1}|$ has a lower bound whenever $w,w^\prime\in C_r$. Let us choose $R$ in such a way that $b^R$ becomes less than the lower bound of $|w^\prime w^{-1}|$. This ensures that the infimum of $|b^sw-w^\prime|$ is bounded away from $0$ along $(\xi_k, R-id]\cup[R+id,\xi^\prime_k)\cup I_k$. Further, note that $bC_r$ never touches $C_r$. Thereafter, using boundedness of $C_r$, we can found $\delta$ sufficiently close to $1$ such that $b^sC_r$ remains disjoint from $C_r$ for all $s\in [\delta-id,\delta+id]$. However, it implies somewhat much stronger statement. On can say further that $b^sC_r$ escapes from $C_r$ for all $s\in [\delta-id,R-id]\cup[\delta+id,R+id]$. Hence, condition $(A)$ in \eqref{eq:UniformBoundCond} is satisfied. Furthermore, one can see
\begin{align}\label{eq:gratios}
\frac{g(w)}{g(b^s w)}=\left(\frac{1+b^{s-1}w}{1+b^{-1}w}\right)\exp\left(\frac{tw(1-b^s)}{b(1-b)}\right)\frac{(\alpha b^{s-1}w;b)_\infty}{(\alpha b^{-1}w;b)_\infty}.
\end{align}
Choice of $C_r$ and $R,\delta,d$ implies immediately that condition $(B)$ in \eqref{eq:UniformBoundCond} is also met.
 Now, we turn into proving \eqref{eq:qLaplaceTransform}. Techniques for similar determinantal formula is quite well known by now, thanks to the vast literature mentioned above. We will sketch here the basic outline of such technique. For brevity, we skip major portion of details in the proof. We first look into the expansion of $1/(\zeta b^{N_x};b)_\infty$. Using \cite[Corollary 10.2.2 a]{andrews1999}, one can say
\begin{align}\label{eq:Q-ExpoExpan}
\frac{1}{(\zeta b^{N_x};b)_\infty}=\sum_{L=0}^\infty \frac{\zeta^Lb^{LN_x}}{(b;b)_L}
\end{align}
where $(b;b)_L=\prod_{i=1}^L(1-b\cdot b^{i-1})$. Recall that in Lemma~\ref{MGFSimplified}, we have
\begin{align}\label{eq:LMoment}
\mu_L:=\mathbb{E}_{\mbox{stepb}}(b^{LN_x};t)=\frac{b^{L(L-1)/2}}{(2\pi i)^L}\oint \ldots \oint \prod_{A<B}\frac{u_A-u_B}{u_A-bu_B}\prod_{i=1}^L \frac{f(u_i)}{u_i} du_i,
\end{align}
where \[f(u):=\exp\left(\frac{tu}{b}\right)\left(\frac{1+u}{1+ub^{-1}}\right)^x\frac{1}{1-\alpha ub^{-1}}\]
and the integration contours do include $0,-b$ inside, but, not $\alpha^{-1}b$. Moreover, one can notice that $f$ doesn't have a pole in an open neighborhood of the line joining $0$ and $-b$. Let us denote any partition $\lambda=(\lambda_1,\lambda_2,\ldots)$ of any positive integer $k$ by $\lambda\vdash k$ where $\sum_{i=1}^\infty \lambda_i=k$. Further, denote number of positive $\lambda_i$'s in $\lambda$ by $l(\lambda)$. Thus, using \cite[Proposition~5.2]{BCF12}, one can conclude that
\begin{align}\label{eq:NewFormLMoment}
\mu_L=(b;b)_L &\sum_{\stackrel{\lambda\vdash L}{1^{m_1}2^{m_2}\ldots }}\frac{1}{m_1!m_2!\ldots }\frac{1}{(2\pi i)^{l(\lambda)}}\\&\times \int_C\ldots \int_C \mbox{det}\left(\frac{-1}{b^{\lambda_i}w_i-w_j}\right)_{i,j=1}^{l(\lambda)}\prod_{j=1}^{l(\lambda)}f(w_j)f(bw_j)\ldots f(b^{\lambda_j-1}w_j)dw_j
\end{align}
where the contour $C$ also avoids $0,-b$ and doesn't intersect with $b^nC$ for any $n\in\NN$. To this end, one can conform the contour $C$ to the the circle $C_r$ without passing any singularity of the integrand. Further, using the form of the function $g$, it is easy to see
\[\frac{g(w_j;1,b,x,t)}{g(b^{\lambda_j}w_j;1,b,x,t)}=f(w_j)f(bw_j)\ldots f(b^{\lambda_j-1}w_j).\]
Consequently, we can write
\begin{align}\label{eq:LaplaceExpan}
\sum_{L=0}^\infty \frac{\zeta^L\mu_L}{(b;b)_L}=\sum_{L=0}^\infty\sum_{\stackrel{\lambda\vdash L}{1^{m_1}2^{m_2}\ldots }}&\frac{\zeta^L}{m_1!m_2!\ldots }\frac{1}{(2\pi i)^{l(\lambda)}}\\&\times \int_{C_r}\ldots \int_{C_r} \mbox{det}\left(\frac{-1}{b^{\lambda_i}w_i-w_j}\frac{g(w_j;1,b,x,t)}{g(b^{\lambda_j}w_j;1,b,x,t)}\right)_{i,j=1}^{l(\lambda)}\prod_{j=1}^{l(\lambda)}dw_j.
\end{align}
At this point, we will make an interchange of two summations above. Let us fix a positive integer $k$ and first sum over all partitions $\lambda$ which has size $k$, i.e $l(\lambda)=k$. Interchange of summations is justified by the linearity of the integrals. Furthermore, one can write $\zeta^L=\zeta^{\lambda_1+\ldots +\lambda_{l(\lambda)}}$. To that effect, the right side of \eqref{eq:LaplaceExpan} simplifies to
\begin{align}
\sum_{k=0}^\infty \frac{1}{k!(2\pi i)^k}\int_C\ldots \int_C \mbox{det}\left(\sum_{L=1}^\infty\frac{(-\zeta)^s}{b^{L}w_i-w_j}\frac{g(w_j;1,b,x,t)}{g(b^{L} w_j;1,b,x,t)}\right)_{i,j=1}^{k}\prod_{j=1}^{k}dw_j.
\end{align}
  Now, it remains to show that
  \begin{align}\label{eq:SimpleLaplaceExpan}
  \sum_{L=1}^\infty\frac{(-\zeta)^s}{b^{L}w_i-w_j}\frac{g(w_j;1,b,x,t)}{g(b^{L} w_j;1,b,x,t)}=\frac{1}{2\pi i}\int_{D_{R,d,\delta}}\frac{(-\zeta)^s}{\sin(\pi s)}\cdot\frac{g(w;1,b,x,t)}{g(b^sw;1,b,x,t)}\cdot\frac{ds}{b^sw-w^\prime}.
  \end{align}
 Note that the function $1/\sin(\pi s)$ decays exponentially as $|s|\to \infty $ along the segments $(R-i\infty,R-id]$ and $[R+id,R+i\infty)$. On the flip side, $(b^s w-w^\prime)^{-1}g(w)/g(b^sw)$ is uniformly bounded on $D_{R,d,\delta,\kappa}$ for all $k\leq R$. Thus, integrand in \eqref{eq:SimpleLaplaceExpan} goes to zero at at both ends of the line $\mbox{Re}(z)=R$. Instead of integrating over $D_{R,d,\delta}$, if we choose to integrate over $D_{R,d,\delta; \kappa}$, then using residue theorem, we obtain
 \begin{align}\label{eq:ResidueCal}
 \frac{1}{2\pi i}\int_{D_{R,d,\delta; \kappa}}\frac{(-\zeta)^s}{\sin(\pi s)}\cdot\frac{g(w;1,b,x,t)}{g(b^sw;1,b,x,t)}\cdot &\frac{ds}{b^sw-w^\prime}\\
 &= \sum_{L=1}^{\kappa}\mbox{Res}_{s=L}\left[\frac{(-\zeta)^s}{\sin(\pi s)}\cdot\frac{g(w;1,b,x,t)}{g(b^sw;1,b,x,t)}\cdot\frac{1}{b^sw-w^\prime}\right].
 \end{align}
 Moreover, if we let $\kappa$ grows to infinity, then the integral along the semi-circular arc $I_\kappa$ tends $0$ thanks again to the exponential growth of $\sin(\pi s)$. Thus, integral over $D_{R,d,\delta;\kappa}$ converges to the right side in \eqref{eq:SimpleLaplaceExpan} whereas when $\kappa\to \infty$, the right side in \eqref{eq:ResidueCal} recovers the other side of the equality in \eqref{eq:SimpleLaplaceExpan}. This completes the proof.
\end{proof}

\section{Asymptotics}\label{secAsymptotics}

Main aim of this section is to prove the Theorem~\ref{Asymptheo1Mainstate} and Theorem~\ref{Asymptheo2Mainstate}. Asymptotic analysis that we present here adapts closely to the similar expositions in \cite[Section~5]{BCG2016}. The crux of the proof of Theorem~\ref{Asymptheo1Mainstate} essentially is governed by the the steepest descent analysis of the kernel of the Fredholm determinant in Theorem~\ref{FredRep}. In the following propositions, large time limit of the Fredholm determinant will be provided. We conclude the proof of Theorem~\ref{Asymptheo1Mainstate} after Lemma~\ref{GenLemma} modulo the proof of these propositions.

\bp\label{FinalAsympProp}
Recall the definition of the contour $C_r$ from Theorem~\ref{FredRep}. Adopt the notations used in Theorem~\ref{Asymptheo1Mainstate}.
\ba
\ii Assume $\rho >1-\left(\nu(1-b)\right)^{-1/2}$. Set $x(t)=\lfloor\nu t\rfloor$ and  $\zeta_t= -b^{-m_\nu  t+s\sigma_\nu t^{1/3}}$ . Then, we have
\begin{align}\label{eq:TW1limit}
\lim_{t\to \infty} \det\left(1+K^b_{\zeta_t}\right)_{L^2(C_r)}= F_{GUE}(s).
\end{align}
\ii Let $\rho= 1-\left(\nu(1-b)\right)^{-1/2}$. For the same choices of $x(t)$ and $\zeta_t$ as above, we have
\begin{align}\label{eq:TW2limit}
\lim_{t\to \infty} \det\left(1+K^b_{\zeta_t}\right)_{L^2(C_r)}= F^2_{GOE}(s).
\end{align}
\ii In the case when $\rho <1-\left(\nu(1-b)\right)^{-1/2}$, set $x(t)=\lfloor\nu^\prime t\rfloor$ and $\zeta_t= -b^{-m_\nu  t+s\sigma^\prime\nu t^{1/3}}$. Then, one can get
\begin{align}\label{eq:Gaussianlimit}
\lim_{t\to \infty} \det\left(1+K^b_{\zeta_t}\right)_{L^2(C_r)}= \Phi(s).
\end{align}
\ee
\ep

  \bl (\cite[Lemma~4.1.39]{BorodinCorwinMac})\label{GenLemma}
  Consider a set of functions $\{f_t\}_{t\geq 0}$ mapping $\RR\to [0,1]$ such that for each $t$, $f_t(x)$ is strictly decreasing in $x$ with a limit of $1$ at $x=-\infty$ and $0$ at $x=\infty$ and for each $\delta>0$, on $\RR\backslash [-\delta,\delta]$, $f_t$ converges uniformly to $\mathbf{1}(x\leq 0)$. Define the $r$-shift of $f_t$ as $f^{r}_n(x)=f_n(x-r)$. Consider a set of random variables $X_t$ indexed by $t\in \RR^+$ such that for each $r\in \RR$,
  \begin{align}\label{eq:CondOnLimit}
  \mathbb{E}[f^r_t(X_t)]\stackrel{t\to \infty}{\longrightarrow} p(r)
  \end{align}
  and assume that $p(r)$ is a continuous probability distribution. Then $X_t$ weakly converges in distribution to a random variable $X$ which is distributed according to $\mathbb{P}(X\leq r)=p(r)$.
  \el

  \textsc{Proof of Theorem~\ref{Asymptheo1Mainstate}}.  For part $(a)$ and $(b)$ of Theorem~\ref{Asymptheo1Mainstate}, we use Lemma~\ref{GenLemma} with the functions
  \[f_t(z):=\frac{1}{(-b^{-t^{1/3}z};b)_\infty}, t\in \RR^+.\]
  Note that $f_t(z)$ is a monotonously decreasing function of the real argument of $z$, yielding $\lim_{t\to\infty}f_t(z)=0$ if $z>0$ and $\lim_{t\to-\infty}f_t(Z)=1$ if $z\leq 0$. Set $X_t=\sigma^{-1}_{\nu}t^{-1/3}(m_\nu t-N_{\nu t}(t))$. Thus, for any $s\in\RR$, we can write
  \[f^s_t(X_t)= \frac{1}{(\zeta_tb^{N_{\nu t}(t)};b)_\infty}\]
  where $\zeta_t=b^{-m_\nu  t+s\sigma_\nu t^{1/3}}$. To this end, part $(a)$ and $(b)$ of Proposition~\ref{FinalAsympProp} shows that the limit of $\mathbb{E}\left(f^s_t(X_t)\right)$ are indeed probability distribution function. Hence, this proves part $(a)$ and $(b)$ of Theorem~\ref{Asymptheo1Mainstate}. For part $(c)$, we replace $f_t$ by $g_t$ which is defined as
  \[g_t(z):=\frac{1}{(-b^{-t^{1/2}z};b)_\infty}, t\in \RR^+.\]
  Setting $X_t=(\tilde{\sigma}_{\nu})^{-1}t^{-1/2}(\tilde{m}_\nu t-N_{\nu t}(t))$ and $\zeta_t=b^{-m_\nu  t+s\tilde{\sigma}_\nu t^{1/3}}$, one can note $g^s_t(X_t)= (\zeta_tb^{N_{\nu t}(t)};b)^{-1}_\infty$. Therefore, rest of the proof follows from the part $(c)$ of Proposition~\ref{FinalAsympProp}.

  \textsc{Proof of Proposition~\ref{FinalAsympProp}} We start with some necessary pre-processing of the kernel $K^b_{\zeta_t}$. Notice that in the statements of all three parts of Proposition~\ref{FinalAsympProp}, $\zeta_t$ is substituted with $-b^{\beta}$ and $x$ is substituted with $\lfloor\nu t\rfloor$ where $\beta$ is some function of the time variable $t$ and $\nu>(1-b)^{-1}$.  Under these substitutions, the form of the kernel $K^b_{\zeta_t}$ effectively reduces to
  \begin{align}\label{eq:KernelIntMedForm}
  K^{b}_{\zeta_t}(w,w^\prime)=\frac{1}{2i}\int_{D_{R,d,\delta}}\frac{b^{h\beta}}{\sin(\pi h)}\cdot\frac{g(w;1,b,\nu t,t)}{g(b^hw;1,b,\nu t,t)}\cdot\frac{dh}{b^hw-w^\prime}
  \end{align}
  where the contour of the integration is oriented from the bottom to top and
  $g(z;1,b,x,t)=\left(\frac{1}{1+z/b}\right)^x\exp\left(\frac{t}{b(1-b)}\right)(\alpha b^{-1}z;b)_\infty^{-1}$. Like in \cite{BCG2016}, the particular form of the integrand in \eqref{eq:KernelIntro} calls for the substitution $b^hw=v$. Unfortunately, the latter substitution is not injective. Thus, to study how this effect the integral, we have to divide the contour $D_{R,d,\delta}$ into countably many parts which have one-to-one images under the above substitution. For instance, each of the linear sections $\mathcal{I}_l:=[R+i(d+2\pi\log(b)^{-1}l ), R+i(d+2\pi (\log(b))^{-1}(l+1))]$ of the contour of $h$ for $l\in \ZZ\backslash g\{-1\}$ maps to the same image. Furthermore, other parts of $D_{R,d,\delta}$, i.e., the union
  \[[R-i(d+2\pi\log(b)^{-1}), R-id]\cup [R-id, \delta -id]\cup [\delta-id,\delta +id]\cup [\delta+id,R+id]\]
  maps to a closed contour $D_{R,d,\delta,w}$ composed of four different portions. For completeness, we present here a formal definition of the contour $D_{R,d,\delta,w}$.
\bd\label{eq:DefOfDeformedContour}
 For each $w$, define a contour $D_{R,d,\delta,w}$ as the union
 \[T_1\cup T_2\cup T_3\cup T_4.\]
 We describe $T_i$'s successively. To begin with, $T_1$ is a major arc of a small circle around origin with radius $b^R|w|$. To be exact, $T_1$ is the image of the segment $[R-i(d+2\pi\log(b)^{-1}), R-id]$ under the map $h\mapsto b^h w$. Similarly, $T_2$, image of $[\delta-id,\delta +id]$, is a minor arc of a larger circle of radius $b^\delta|w|$. Apparently, these two arcs of two different circle are connected by two linear segments $T_2$ and $T_3$. Precisely, $T_2$ connects $b^{R-id}w$ to $b^{\delta-id}w$ and $T_3$ connects $b^{\delta+id}w$ to $b^{R+id}w$. It is evident from the definition of $T_2$ ($T_3$) that it is the image of $[R-id, \delta -id]$ ($[\delta+id, R+id]$) under the above map. Due to the bottom-to-top orientation of the contour $D_{R,d,\delta}$ in Theorem~\ref{FredRep}, $D_{R,d,\delta,w}$ will be anti clockwise oriented. See Figure~\ref{fig:4} for further details.
\ed

\begin{figure}
\begin{center}
\begin{tikzpicture}
  \draw[line width = 0.25mm, <->] (-3,0) -- (2,0) node {};
  \draw[line width = 0.25mm, <->] (-0.5,-3) -- (-0.5,3) node {};
  \draw[line width = 0.4mm, -] (0.649-0.5,0.375) -- (1.732-0.5, 1) node {};
  \draw[line width = 0.4mm, -] (0.375-0.5, 0.649) -- (1-0.5,1.732) node {};
  \draw[line width = 0.4mm, -] (5, 0.5) -- (8, 0.5) node{};
  \draw[line width = 0.4mm, -] (5,-0.5) -- (8, -0.5) node {};
  \draw[line width = 0.4mm, -] (8,-0.5) -- (8, -1.5) node {};
  \draw[line width = 0.4mm, -] (5,0.5) -- (5, -0.5) node {};
  \draw[line width = 0.2mm, ->] (5,0.75) .. controls (4,1.25) and (3,1.15) .. (2.2,0.95);
  \node[] at (-1.2,-0.7) (a) {$T_1$};
  \node[] at (1.3,1.4) (b) {$T_4$};
  \node[] at (1,0.45) (c) {$T_2$};
  \node[] at (-0.05,1.25) (d) {$T_3$};
  \node[] at (3.5,1.35) (e) {$wb^{h}$};
  \node[] at (8.5,0.5) (f) {$R+id$};
  \node[] at (8.5,-0.5) (g) {$R-id$};
  \node[] at (8.5,-1.95) (h) {$R-i(d+2\pi\log(b)^{-1})$};
   \draw[domain=1:1.732, variable=\x, line width = 0.2mm] plot ({\x-0.5},{sqrt(4-\x^2)}) ;
   \draw[domain=0:0.375, variable=\x, line width = 0.2mm] plot ({\x -0.5},{sqrt(0.5625-\x^2)}) ;
   \draw[domain=0:0.75, variable=\x, line width = 0.2mm] plot ({\x -0.5},{-sqrt(0.5625-\x^2)}) ;
   \draw[domain=0:0.75, variable=\x, line width = 0.2mm] plot ({-\x-0.5},{-sqrt(0.5625-\x^2)}) ;
   \draw[domain=0:0.75, variable=\x, line width = 0.2mm] plot ({-\x-0.5},{sqrt(0.5625-\x^2)}) ;
   \draw[domain=0.649:0.75, variable=\x, line width = 0.2mm] plot ({\x-0.5},{sqrt(0.5625-\x^2)}) ;
\end{tikzpicture}
\caption{Figure on the right side indicates the section of the contour $D_{R,d,\delta}$ which under the map $h\mapsto wb^{h}$ pertains to the contour $T_1\cup T_2\cup T_3\cup T_4$ shown in the left.}
\label{fig:4}
\end{center}
\end{figure}
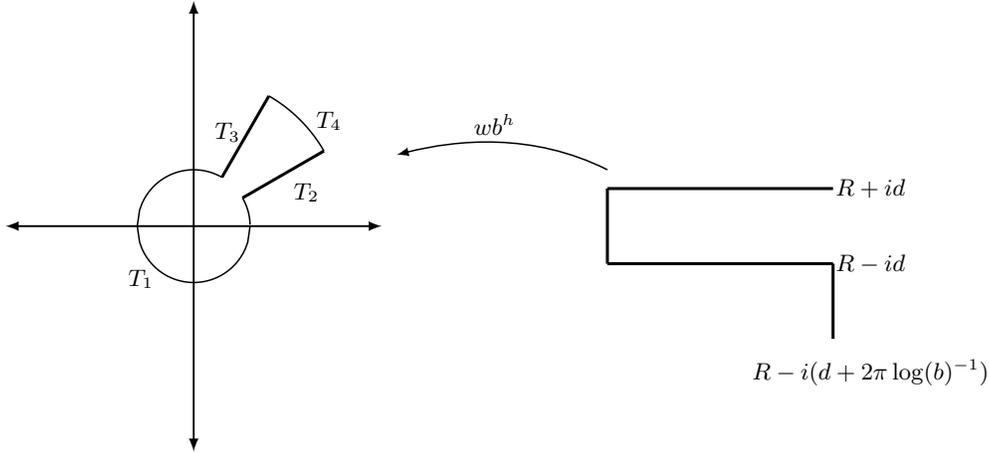

   For notational convenience, let us denote the circle of radius $s$ with center at origin in the complex plane $\CC$ by $\mathcal{C}_{s}$ (note the difference from $C_r$). In fact, $\mathcal{C}_{b^R|w|}$ is the image of each linear segments $\mathcal{I}_l$ for $l\in \ZZ\backslash \{-1\}$. Thus, one can write
   \begin{align}\label{eq:KernelIntMedForm2}
  K^{b}_{\zeta_t}(w,w^\prime)&=\sum_{\substack{k=-\infty\\ k\neq -1}}^\infty\frac{1}{2 i}\int_{\mathcal{C}_{b^R|w|}}\frac{(v/w)^{\beta}}{\sin\left(\frac{\pi}{\log (b)}(\log(v/w))+2\pi ik\right)}\cdot\frac{g(w;1,b,\nu t,t)}{g(v;1,b,\nu t,t)}\cdot\frac{dv}{v\log(b)(v-w^\prime)}\nonumber\\
  &+ \frac{1}{2 i}\int_{D_{R,d,\delta,w}}\frac{(v/w)^{\beta}}{\sin\left(\frac{\pi}{\log (b)}(\log(v/w))-2\pi i\right)}\cdot\frac{g(w;1,b,\nu t,t)}{g(v;1,b,\nu t,t)}\cdot\frac{dv}{v\log(b)(v-w^\prime)}.
  \end{align}
  At this point, we would make a surgery over the contour $D_{R,d,w,k}$. Let us connect a minor arc from $ b^{R-id}w$ to $ b^{R-i(d+2\pi\log(b)^{-1})}w$ in the circle $\mathcal{C}_{b^{R}|w|}$ (major arc being $T_1$) and call it $T^\prime_1$. We use $-T^\prime_1$ to indicate clockwise orientation. Therefore, $\tilde{D}_{R,d,w,k}:=T_2\cup T_3\cup T_4 \cup -T^\prime_1$ denotes a closed contour oriented anti clockwise. Using the properties of the contour integration, it can be said that integral over $D_{R,d,w,k}$ in \eqref{eq:KernelIntMedForm2} is a sum of the integrals over $\tilde{D}_{R,d,w,k}$ and $\mathcal{C}_{b^R|w|}$. Further, notice that there is no pole of the integrand sitting inside the contour $\tilde{D}_{R,d,w,k}$ for any $w\in C_r$ (see, Theorem~\ref{FredRep} for $w$) thanks to the conditions $(A)$ and $(B)$ in \eqref{eq:UniformBoundCond}. Thus, integral over $\tilde{D}_{R,d,w,k}$ vanishes. To this effect, one can write
  \begin{align}\label{eq:KernelIntMedForm3}
  K^{b}_{\zeta_t}(w,w^\prime)&=\sum_{k=-\infty}^\infty\frac{1}{2 i}\int_{\mathcal{C}_{b^R|w|}}\frac{(v/w)^{\beta}}{\sin\left(\frac{\pi}{\log (b)}(\log(v/w))+2\pi ik\right)}\cdot\frac{g(w;1,b,\nu t,t)}{g(v;1,b,\nu t,t)}\cdot\frac{dv}{v\log(b)(v-w^\prime)}.
  \end{align}

 \subsection{Case $\rho >1-\left(\nu(1-b)\right)^{-1/2}$, Tracy-Widom Fluctuations (GUE).}
 Fix $\nu >(1-b)^{-1}$. Set $\beta= -m_\nu t+s\sigma_\nu t^{1/3}$. Define a function $\Lambda$ as
  \begin{align}\label{eq:LambdaFuncDefine}
  \Lambda(z)= \frac{z}{b(1-b)}-\nu\log(1+z/b)+m_\nu \log(z).
  \end{align}
  Thus, one can further simplify \eqref{eq:KernelIntMedForm3} into
  \begin{align}\label{eq:KernelIntMedForm4}
  K^b_{\zeta_t}(w,w^\prime)=\frac{1}{2 i\log(b)}\sum_{k=-\infty}^{\infty}\int_{\mathcal{C}_{b^R|w|}}&\frac{\exp(t(\Lambda(w)-\Lambda(v))+t^{1/3}s\sigma_\nu(\log (v)-\log(w)))dv}{(v-w^\prime)\sin\left(\frac{\pi}{\log (b)}(\log(v/w))+2\pi ik\right)}\nonumber\\
  &\times \frac{(\alpha b^{-1}v;b)_\infty }{(\alpha b^{-1}w;b)_\infty }\times \frac{dv}{v}.
  \end{align}
  In a nutshell, the asymptotic behavior of the kernel is governed by the variations of the real part of $\Lambda$. In what follows, we exhibit the steepest descent contours, which helps in localizing the main contribution in a neighborhood of a critical point of $\mbox{Re}(\Lambda)$. Thereafter, using Taylor expansion of the arguments of the exponential inside the kernel in \eqref{eq:KernelIntMedForm4} and appropriate estimates for the kernel in the region away from the critical points, we prove the limit.
   To this end, notice that
  \[\Lambda^\prime(z)=\frac{1}{b(1-b)}-\frac{\nu b^{-1}}{1+z/b}+\frac{m_\nu}{z}.\]
  Plugging $m_\nu=\frac{(\sqrt{\nu(1-b)}-1)^2}{1-b}$, we get
  \[\Lambda^\prime(z)=\frac{(z-\varrho)^2}{b^2(1-b)(1+\varrho/b)\varrho},\]
  where $\varrho=b(\sqrt{\nu(1-b)}-1)$.  This shows
  that $\Lambda^\prime$ has an unique double critical point at $\varrho$. In particular, in a neighborhood of $\varrho$, we have
  \begin{align}\label{eq:TaylorExpan}
  \Lambda(z) = \Lambda(\varrho) + \left(\frac{\sigma_\nu}{\varrho}\right)^3(z-\varrho)^3 + \mathfrak{R}((z-\varrho))
  \end{align}
  where $(z-\varrho)^{-3}\mathfrak{R}((z-\varrho))\to 0$ when $z\to \varrho$.
  thanks to fact that
  \[\left(\frac{\sigma_\nu}{\varrho}\right)^3 = \left(b^2\left(\sqrt{\nu(1-b)}-1\right)(1-b)^{3/2}\nu^{1/2}\right)^{-1}=\Lambda^{\prime\prime\prime}(\varrho).\]
  Now, we would like to deform the contours of $v$ and $w$ suitably so that $\mbox{Re}(\Lambda(w)-\lambda(v))<-c$ for some $c>0$ whenever $v$ and $w$ are away from the critical point $\varrho$. In the following lemma, we describe the level curves of $\mbox{Re}(\Lambda(z))=\mbox{Re}(\Lambda(\varrho))$ in details.  We call a closed contour (let's say $\mathfrak{C}$) \emph{star shaped} if for any $\phi\in [0,2\pi)$, there exists exactly one point $z\in \mathfrak{C}$ such that $\mbox{Arg}(z)=\phi$.
  \bl\label{LevelCurveCharacteristic}
   There are exactly three level curves $\mathcal{L}_1$, $\mathcal{L}_2$ and $\mathcal{L}_3$ of $\mbox{Re}(\Lambda(z))=\mbox{Re}(\Lambda(\varrho))$ in the complex plane $\CC$ satisfying the following properties:
   \ba
   \ii both $\mathcal{L}_1$ and $\mathcal{L}_2$ are simple closed contours meeting with $\mathcal{L}_3$ at the point $\varrho$. Further, $\mathcal{L}_1$ ($\mathcal{L}_2$) cuts the negative half of the real line at some point $d_1$ ($d_2$) in the interval $(-b,0)$ ($(-\infty, -b)$),
   \ii $\mathcal{L}_3$ resides in the right half of the complex plane, escapes to infinity as $|z|\to \infty$,
   \ii both $\mathcal{L}_1$ and $\mathcal{L}_2$ are star shaped,
   \ii $\mathcal{L}_1\backslash \{\varrho\}$ is completely contained inside the interior of the region enclosed by $\mathcal{L}_2$. Furthermore, $\mathcal{L}_3$ meet the contours $\mathcal{L}_1$ and $\mathcal{L}_2$ at no other point point except $\varrho$,
   \ii $\mathcal{L}_1$ ($\mathcal{L}_2$) leaves the $x$-axis at an angle $5\pi/6$ ($\pi/2$) from the point $\varrho$ and meets again at $-5\pi/6$ ($-\pi/2$). Moreover, upper half (lower half) of the contour $\mathcal{L}_3$ leaves the real line from $\varrho$ at an angle $\pi/6$ ($-5\pi/6$).
   \ee
  \el
  \begin{proof}

 We closely adapt here the proofs of the related results in the context of stochastic six vertex model (see, \cite[Section~5.1]{BCG2016}).
 \ba
  \ii  It can be observed in \eqref{eq:LambdaFuncDefine} that $\Lambda(z)$ has singularities at $0$ and $-b$ on $\RR$. Thus, there are at least two level curves $\mathcal{L}_1$ and $\mathcal{L}_2$ of $\mbox{Re}(\Lambda(z))=\mbox{Re}(\Lambda(\varrho))$ which originate from the point $\varrho$ and in the way back to their origin, they cut $\RR^-$ at some points $d_1\in (-b,0)$ and $d_2\in (-\infty,-b)$ respectively.

  \ii Tracing the changes in sign of $\mbox{Re}(\Lambda(z))-\mbox{Re}(\Lambda(\varrho))$ in Figure~\ref{fig:5}, one can infer that outside $\mathcal{L}_2$, there should be a region where $\mbox{Re}(\Lambda(z))-\mbox{Re}(\Lambda(\varrho))<0$. But, if we take a large circle $\mathcal{C}_M$ of radius $M$ in $\CC$, then for $z\in \mathcal{C}_M$ and $Re(z)<0$, we have
  \begin{align}\label{eq:LevelCurveEquation}
  \mbox{Re}(\Lambda(z))&=\frac{Re(z)}{b(1-b)} - \nu\log|zb^{-1}(1+z^{-1}b)|+m_\nu\log|z/b|\\
  &=\frac{Re(z)}{b(1-b)} -(\nu -m_\nu)\log|z| +\nu\log(b)-m_\nu\log(b)+o(1)
  <0
\end{align}
     because $m_\nu<\nu$. Consequently, there doesn't exists any other level curve which contains  $0$ or $-b$ inside. But, there are another level curve $\mathcal{L}_3$ escaping to the infinity in the right half of the complex plane after taking off from $\varrho$. This behavior can also be predicted from \eqref{eq:LevelCurveEquation}. Let us sketch the necessary arguments behind such predictions. Note that if we increase the value $\mbox{Re}(z)$ to the infinity with a fixed value of $\mbox{Im}(z)$, then the right side of \eqref{eq:LevelCurveEquation} goes to infinity as well. This shows that at least in the right half of $\CC$, some more solutions (other than those on $\mathcal{L}_1$ and $\mathcal{L}_2$) of the equation $\mbox{Re}(\Lambda(z))=\mbox{Re}(\Lambda(\varrho))$ exist. This verifies the second claim.

   \ii  To show that both $\mathcal{L}_1$ and $\mathcal{L}_2$ are star shaped, we need to prove that any line from origin can cut $\mathcal{L}_1$ or $\mathcal{L}_2$ exactly at two points. Let us fix a line $rz_0$ in $\CC$ where $z_0$ is fixed with $\mbox{Re}(z_0)=1$ and $r$ varies over the real line. To this end, we have
     \begin{align}\label{eq:DerivLevelCurve}
     \frac{\partial }{\partial r}\mbox{Re}(\Lambda(z)) = \mbox{Re}\left[\frac{(z_0)}{b(1-b)}- \frac{\nu z_0b^{-1}}{1+rz_ob^{-1}}+m_\nu\frac{1}{r}\right].
\end{align}
One can note that solutions of the equation that the right side of \eqref{eq:DerivLevelCurve} equals to some fixed value is given by a cubic polynomial. But any cubic polynomial over $\RR$ can have  either one or three real roots. Next, we show $\frac{\partial }{\partial r}\mbox{Re}(\Lambda(rz_0)))=0$ has exactly three real roots unless $\mbox{Im}(z_0)=0$. It can be readily verified that $(i)$ when $r\downarrow 0$, then right side of \eqref{eq:DerivLevelCurve} increases to  infinity, $(ii)$ when $r\to \infty$, then $\frac{\partial }{\partial r}\mbox{Re}(\Lambda(rz_0)))\to b^{-1}(1-b)^{-1}$ and when $r=\varrho$, then $\frac{\partial }{\partial r}\mbox{Re}(\Lambda(rz_0)))<0$ unless $\mbox{Im}(z_0)=0$. This shows the equation $\frac{\partial }{\partial r}\mbox{Re}(\Lambda(rz_0)))=0$ has at least two distinct solutions on $\RR^{+}$ unless $\mbox{Im}(z_0)=0$. Further, co-efficient of $r^3$ and the constant term in the denominator of \eqref{eq:DerivLevelCurve} have same sign. This proves the existence of three real roots out of which two must be positive and last one should be negative. Even those two positive roots are distinct unless $\mbox{Im}(z_0)=0$. Thus, except for the case when $\mbox{Im}(z_0)=0$, none of the roots is an inflection point. We have seen $\mbox{Re}(\Lambda(z))=\mbox{Re}(\Lambda(\varrho))$ has at least five roots along the line $rz_0$ (four on $\mathcal{L}_1$ \& $\mathcal{L}_2$ and fifth one on $\mathcal{L}_3$). As we know $\frac{\partial }{\partial r}\mbox{Re}(\Lambda(rz_0))$ vanishes in between any two roots of the same sign, therefore, number of roots of $\mbox{Re}(\Lambda(rz_0))=\mbox{Re}(\Lambda(\varrho))$ cannot be greater than five. Otherwise, pairs of consecutive roots with the same sign must be greater than three contradicting the fact that the numerator of the right hand side in \eqref{eq:DerivLevelCurve} is a polynomial of degree $3$ in $r$. Thus, the third claim is proved.

 \ii Fourth claim again follows from the fact that the equation $\frac{\partial }{\partial r}\mbox{Re}(\Lambda(rz_0)))=0$ has exactly three solutions in $r$ except when $\mbox{Im}(z_0)=0$. If $\mathcal{L}_1$ meets $\mathcal{L}_2$ at some point $\varrho^\prime \neq \varrho$, then $\frac{\partial }{\partial r}\mbox{Re}(\Lambda(r\varrho^\prime)))=0$ is supposed to have at most two real solutions. Thus, $\varrho^\prime$ cannot be different from $\varrho$. But in the latter case, $\mbox{Re}(\Lambda(r\varrho))=\mbox(\Lambda(\varrho))$ is satisfied for three distinct values of $r$, namely, $1$, $d_1\varrho^{-1}$ and $d_2\varrho^{-1}$. As a matter of fact, we know $d_1\neq d_2$. Consequently, $\mathcal{L}_1$ has no other point of intersection with $\mathcal{L}_2$ except at $\varrho$. Similarly, one can prove $\mathcal{L}_3$ doesn't meet $\mathcal{L}_1$ or $\mathcal{L}_2$ at any other point besides $\varrho$. Moreover, positions of $d_1$ and $d_2$ over the real line implies $\mathcal{L}_1\backslash\{\varrho\}$ is completely contained inside the interior of the region bounded by $\mathcal{L}_2$.

  \ii To begin with, note that $\Lambda(z)$ in a neighbouhood of the point $\varrho$ behaves as $\Lambda(\varrho)+\sigma^3_{\nu}\varrho^{-3}(z-\varrho)^3$. To that effect, possible choices of the tangents for any of the level curves of the equation $\mbox{Re}(\Lambda(z))=\mbox{Re}(\Lambda(\varrho))$ are given by $\{\pi/6\pm 2\pi/3\}\cup \{\pi/6\pm 4\pi/3\}=\{\pm \pi/6,\pm \pi/2, \pm 5\pi/6\}$. Moreover, we know $\mathcal{L}_1\backslash \{0\}$ is completely contained inside the interior of the region bounded by the contour $\mathcal{L}_2$. Also, $\mathcal{L}_3$ stays outside of both the contours $\mathcal{L}_1$ and $\mathcal{L}_2$. These together justify the claims in part $(e)$.
\ee
  \end{proof}
  \bd\label{wvContourDef}
  We define two new contours $\mathcal{L}_w$ and $\mathcal{L}_v$. Both $\mathcal{L}_v$ and  $\mathcal{L}_w$ consist of a piecewise linear segment ($\mathcal{L}^{(1)}_v$ and $\mathcal{L}^{(1)}_w$) and a curved segment ($\mathcal{L}^{(2)}_v$ and $\mathcal{L}^{(2)}_w$). To begin with, $\mathcal{L}^{(1)}_w$ extends linearly from $\varrho+\epsilon e^{-i\pi/3}$ to $\varrho$  and from there to $\varrho+\epsilon e^{i\pi/3}$. Similarly, $\mathcal{L}^{(1)}_v$ goes from $\varrho-\varrho\sigma^{-1}_\nu t^{-1/3}+\epsilon e^{-2i\pi/3}$ to $\varrho-\varrho\sigma^{-1}_\nu t^{-1/3}$ to $\varrho-\varrho\sigma^{-1}_\nu t^{-1/3}+\epsilon e^{2i\pi/3}$. Here, $\epsilon$ will be chosen in such a way that $(i)$ $\mathcal{L}^{(1)}_v$ is completely contained inside $\mathcal{L}_2$ and $\mathcal{L}^{(1)}_w\backslash \{\varrho\}$ does not intersect $\mathcal{L}_2$ anywhere else, $(ii)$ $|\mathfrak{R}((z-\varrho))|$ is bounded above by $c|z-\varrho|^3$ for some small constant $c$. Next, $\mathcal{L}^{(2)}_w$ starts from the point $\varrho+\epsilon e^{i\pi/3}$ and encircles around $\mathcal{L}_2$ to join $\varrho+\epsilon e^{-i\pi/3}$. On the contrary, $\mathcal{L}^{(2)}_v$ starts from the point $\varrho-\varrho\sigma^{-1}_\nu t^{-1/3}+\epsilon e^{2i\pi/3}$ and curves around $0$ and $-b$ to meet $\varrho-\varrho\sigma^{-1}_\nu t^{-1/3}+\epsilon e^{-2i\pi/3}$. But, it never  goes out of the contour $\mathcal{L}_2$. It can be ascertained that both $\mathcal{L}_v$ and $\mathcal{L}_w$ are bounded away from $\mathcal{L}_2$ in a $\epsilon$-neighborhood of $\varrho$ by choosing $\epsilon$ appropriately.
  \ed

  \begin{figure}[H]
\begin{center}
\begin{tikzpicture}
  \draw[line width = 0.25mm, <->] (-4,0) -- (6,0) node {};
  \draw[line width = 0.25mm, <->] (0,-4) -- (0,4) node {};
   \draw[domain=0:2, variable=\x, dashed, line width = 0.2mm] plot ({\x},{sqrt(4-\x^2)}) ;
   \draw[domain=0:2, variable=\x, dashed, line width = 0.2mm] plot ({\x},{-sqrt(4-\x^2)}) ;
   \draw[domain=0:2, variable=\x, dashed, line width = 0.2mm] plot ({-\x},{sqrt(4-\x^2)}) ;
   \draw[domain=0:2, variable=\x, dashed, line width = 0.2mm] plot ({-\x},{-sqrt(4-\x^2)}) ;
   \draw[domain=0:2, variable=\x, dashed, line width = 0.2mm] plot ({\x},{sqrt(6.25-(\x-0.5)^2)-2}) ;
   \draw[domain=0:1, variable=\x, dashed, line width = 0.2mm] plot ({-\x},{sqrt(6.25-(\x+0.5)^2)-2}) ;
   \draw[domain=0:2, variable=\x, dashed, line width = 0.2mm] plot ({\x},{-sqrt(6.25-(\x-0.5)^2)+2}) ;
   \draw[domain=0:1, variable=\x, dashed, line width = 0.2mm] plot ({-\x},{-sqrt(6.25-(\x+0.5)^2)+2}) ;
   \draw[line width = 0.2mm, dashed, ->] (2,0) .. controls (2.5,0.35) and (3,1.75) .. (5,2.9);
   \draw[line width = 0.2mm, dashed, ->] (2,0) .. controls (2.5,-0.35) and (3,-1.75) .. (5,-2.9);
   \draw[line width =  0.4mm, -] (1.95,0) -- (1.65,0.45) node {};
   \draw[line width =  0.4mm, -] (1.95,0) -- (1.65,-0.45) node {};
   \draw[line width =  0.4mm, -] (2,0) -- (2.2,0.45) node {};
   \draw[line width =  0.4mm, -] (2,0) -- (2.2,-0.45) node {};
   \draw[domain=0:1.65, variable=\x, line width = 0.4mm] plot ({\x},{sqrt(2.925-\x^2)}) ;
   \draw[domain=0:1.65, variable=\x, line width = 0.4mm] plot ({\x},{-sqrt(2.925-\x^2)}) ;
   \draw[domain=0:1.71, variable=\x,  line width = 0.4mm] plot ({-\x},{sqrt(2.925-\x^2)}) ;
   \draw[domain=0:1.71, variable=\x, line width = 0.4mm] plot ({-\x},{-sqrt(2.925-\x^2)}) ;
   \draw[domain=0:2.2, variable=\x, line width = 0.4mm] plot ({\x},{sqrt(5.0425-\x^2)}) ;
   \draw[domain=0:2.2, variable=\x, line width = 0.4mm] plot ({\x},{-sqrt(5.0425-\x^2)}) ;
   \draw[domain=0:2.2455, variable=\x,  line width = 0.4mm] plot ({-\x},{sqrt(5.0425-\x^2)}) ;
   \draw[domain=0:2.2455, variable=\x, line width = 0.4mm] plot ({-\x},{-sqrt(5.0425-\x^2)}) ;
  %
  \node[] at (-0.5,-0.5) (a) {$\mathcal{L}_1$};
  \node[] at (1.0,1.0) (b) {$\mathcal{L}_v$};
  \node[] at (1.3,2.3) (c) {$\mathcal{L}_w$};
  \node[] at (-1.41,1.41) (d) {$\mathcal{L}_2$};
  \node[] at (1,0.15) (e) {$-$};
  \node[] at (-1,1) (f) {$+$};
  \node[] at (-2.4,-1.5) (g) {$-$};
  \node[] at (4,1.25) (h) {$+$};
  \node[] at (4,2.5) (h) {$\mathcal{L}_3$};
  \node[] at (4,-2.5) (i) {$\mathcal{L}_3$};
  \draw[-] (-1.5,-0.05) -- (-1.5,0.05) node[below] {$-b$};
   \draw[-] (3.7,-0.05) -- (3.7,0.05) node[below] {$b\alpha^{-1}$};
   \node[] at (2.25,0.15) (j) {$\varrho$};
\end{tikzpicture}
\caption{In the above plot, dashed lines denote the level curves of the equation $\mbox{Re}(\Lambda(z))=\mbox{Re}(\Lambda(\varrho))$ whereas solid lines indicate the contour $\mathcal{L}_v$ and $\mathcal{L}_w$ respectively. Signature of the function $\mbox{Re}(\Lambda(z)-\Lambda(\varrho))$ has been shown in the respective regions with $\pm$ signs.}
\label{fig:5}
\end{center}
\end{figure}
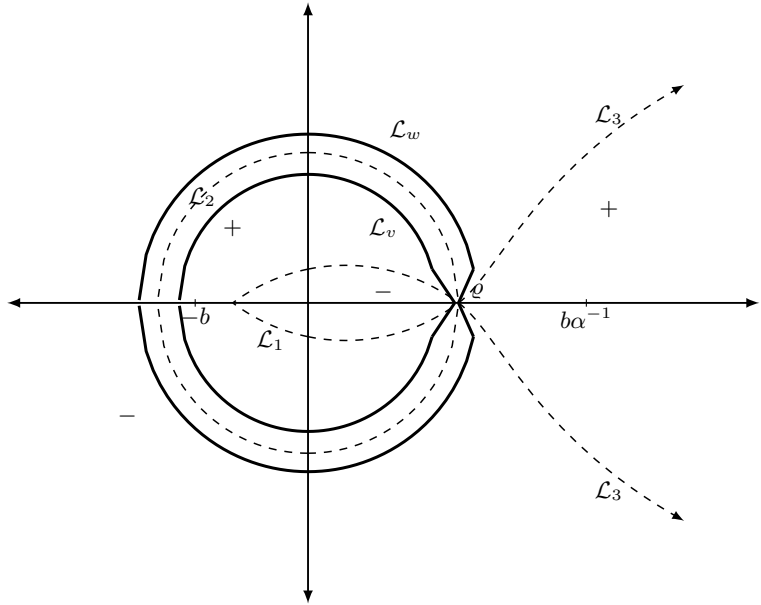

  We deform the contours of $v$ and $w$ to $\mathcal{L}_v$ and $\mathcal{L}_w$ respectively. Such deformations is possible because $\mathcal{L}_v$ sits completely inside $\mathcal{L}_w$ whereas previously, contour for $v$ was $\mathcal{C}_{b^R|w|}$. It can be also noted that we didn't allow any poles to cross in between in doing the deformations of the contours. Consequently, we have $\mbox{Re}(\Lambda(w))<\mbox{Re}(\Lambda(v))$ for all $w\in \mathcal{L}_w$ and for all $v\in \mathcal{L}_v$. In what follows, we show how the contribution of the integrals of $K_b(\zeta_t)$ over the contour $\mathcal{L}_w$ dominates the rest in the Fredholm determinant in Proposition~\ref{FinalAsympProp}~\eqref{eq:TW1limit}.
  \bl\label{TailBoundLemma}
 There exists some $c, C>0$ and $t_0>0$, such that for all $t\geq t_0$
\begin{align}\label{eq:KernelBound}
|K^{b}_{\zeta_t}(w,w^\prime)|\leq C \frac{\max\{\epsilon^{-1}, \varrho^{-1} \sigma_\nu t^{1/2}\}}{\log(1+\epsilon)}\exp\left(-c(t+|w|)+s\tilde{\sigma}_\nu t^{1/3}|w|\vee |w|^{-1}\right)
\end{align}
where $w\in \mathcal{L}^{(2)}_w$, $w^\prime\in \mathcal{L}_{w}$ and $\epsilon$ is defined in Definition~\ref{wvContourDef}. Moreover, the following convergence holds
\begin{align}\label{eq:TailKernelDecay}
\lim_{t\to \infty}\left|\det(1+K^b_{\zeta_t})_{L^2(\mathcal{L}_w)}-\det(1+K^b_{\zeta_t})_{L^2(\mathcal{L}^{(1)}_w)}\right| = 0
\end{align}
\el
\begin{proof}
To show \eqref{eq:KernelBound}, we need to bound the integrand in \eqref{eq:KernelEstimate4} with suitable integrable functions. In the following discussion, we enlist the upper bounds on the absolute value of each of the components of the integrand.
\bei
\ii For any two points $v$ and $w$ on the contours $\mathcal{L}_v$ and $\mathcal{L}^{(2)}_w$ respectively, we know $\mbox{Re}(\Lambda(w))<\mbox{Re}(\Lambda(v))$. In fact, using the compactness of the respective contours, we can further say that there exists $c>0$ such that $\mbox{Re}(\Lambda(w))-\mbox{Re}(\Lambda(v))<-c$ whenever $v\in \mathcal{L}_v$  and $w\in\mathcal{L}^{(1)}_w$. Furthermore, for any complex number $z$, we have $\log|z|\leq |z|\vee |z|^{-1}$. To sum up, we get
\begin{align}
\left|\exp(t(\Lambda(w)-\Lambda(v))+t^{1/3}s\sigma_\nu(\log (v)-\log(w)))\right|\leq \exp(-tc+s\sigma_\nu t^{1/3}|w|\vee |w|^{-1}).\\\label{eq:BoundsOnACExpTerm}
\end{align}
\ii For any $k\in \ZZ\backslash \{0\}$, a crude upper bound on the sine function in the integrand is given by
\begin{align}\label{eq:BoundOnSinFuncAC}
|\sin\left(\frac{\pi}{\log (b)}(\log(v/w))+2\pi ik\right)|\geq C^\prime\exp(2\pi |k|)
\end{align}
for some constant $C^\prime >0$. For the case when $k=0$,  $|\sin((\pi/\log(b))\log(v/w))|$ is bounded below by $C^{\prime\prime}\min_{v\in \mathcal{L}_v, w\in \mathcal{L}^{(2)}_w}|\log |(v/w)||$. Recall that $\mathcal{L}^{(2)}_w$ stays at least $\epsilon$ away from the contour $\mathcal{L}_2$. Thus, we have
\begin{align}\label{eq:BoundOnSinACWhenK0}
\min_{v\in \mathcal{L}_v, w\in \mathcal{L}^{(2)}_w}|\log |(v/w)||\geq \log (1+\epsilon).
\end{align}
Henceforth, \eqref{eq:BoundOnSinFuncAC} and \eqref{eq:BoundOnSinACWhenK0} implies that
\begin{align}\label{eq:BoundOnSeriesAC}
\sum_{k\in \ZZ}\left|\frac{1}{sin\left(\frac{\pi}{\log (b)}(\log(v/w))+2\pi ik\right)}\right|\leq \frac{2}{\log(1+\epsilon)}\sum_{k=0}^\infty \exp(-2\pi k).
\end{align}

\ii Among the other components, the following bounds can be readily proved.
\begin{align}\label{eq:InfiniteProdBoundAC}
\left|\frac{(\alpha b^{-1}v;b)_\infty }{(\alpha b^{-1}w;b)_\infty }\right|&\leq \exp\left( \frac{\alpha b^{-1}}{1-b}(|v|+|w|)\right),\\
\left|\frac{1}{v-w}\right|&\leq \max_{v\in \mathcal{L}_v, w^\prime \in \mathcal{L}_w}\left|\frac{1}{v-w^\prime}\right|\leq \max\{\epsilon^{-1}, \varrho^{-1} \sigma_\nu t^{1/2}\}.
\end{align}
\ee
Thus, these completes the proof of \eqref{eq:KernelBound}. Next, we show the limit in \eqref{eq:TailKernelDecay}. For that, decompose $K^b_{\zeta}$ as
\begin{align}\label{eq:PreSubDecom}
K^{b}_{\zeta_t}(w,w^\prime)
& := K^{b,1}_{\zeta_t}(w, w^\prime) +K^{b}_{\zeta_t}(w,w^\prime),\\
&= K^{b}_{\zeta_t}(w, w^\prime) \mathbbm{1}(w\in \mathcal{L}^{(1)}_{w})+K^{b}_{\zeta_t}(w,w^\prime)\mathbbm{1}(w\in \mathcal{L}^{(2)}_{w}).
\end{align}
Using the definition of the Fredholm determinant in \eqref{eq:FredDef}, one can write
\begin{align}\label{eq:SubtractionSurgery}
 \left|\det(1+K^{b}_{\zeta_t})_{L^2(\mathcal{L}_w)}\right. &\left.-\det(1+K^b_{\zeta_t})_{L^2(\mathcal{L}^{(1)}_w)}\right|\\
 & \leq \sum_{n\geq 0}\frac{1}{n!(2\pi i)^n}\int_{\mathcal{L}_{w}}\ldots \int_{\mathcal{L}_{w}} |\det(K^{b}_{\zeta_t}(w_i, w_j))-\det(K^{b,1}_{\zeta_t}(w_i,w_j))|dw_1\ldots dw_n.
\end{align}
Furthermore, \eqref{eq:PreSubDecom} implies
\begin{align}\label{eq:matrixdetdiff}
|\det(K^{b}_{\zeta_t}(w_i, w_j))&-\det(K^{b,1}_{\zeta_t}(w_i,w_j))|\leq \sum_{\substack{\tau_1,\ldots ,\tau_n\in \{1,2\}\\ (\tau_1,\ldots ,\tau_n)\\\neq (1,\ldots,1)}}|\det(K^{b, \tau_j}_{\zeta_t}(w_i,w_j))|
\end{align}
where for any fixed $(\tau_1,\ldots ,\tau_n)\in \{1,2\}^n$, $(K^{b,\tau_j}_{\zeta_t}(w_i,w_j))_{i,j=1}^n$ denotes a matrix whose $j$-th column is given by the $j$-th column of $K^{b,\tau_j}_{\zeta_t}$.  Using the bound on the kernel in \eqref{eq:KernelBound}, we have
\begin{align}\label{eq:BoundOnKernelDet}
|\det(K^{b, \tau_j}_{\zeta_t}(w_i,w_j))|\leq C \frac{\epsilon^{-n}t^{n/2}}{(\log(1+\epsilon))^n}\exp\left(-tc+c\sum_{i=1}^n |w_i|+s\sigma_\nu t^{1/2}\sum_{i=1}^n |w_i|\vee |w_i|^{-1}\right).
\end{align}
Moreover, there are $2^n-1$ many terms in the sum in the right side of \eqref{eq:matrixdetdiff}. These implies one can further bound the right side of \eqref{eq:SubtractionSurgery} by
\begin{align}\label{eq:FinalboundOndIffFredDet}
\left|\det(1+K^{b}_{\zeta_t})_{L^2(\mathcal{L}_w)}\right. &\left.-\det(1+K^b_{\zeta_t})_{L^2(\mathcal{L}^{(1)}_w)}\right|\leq C\exp(-ct +c^\prime t^{1/2})
\end{align}
for some constant $C,c,c^\prime >0$. Notice that the right side of \eqref{eq:FinalboundOndIffFredDet} converges to $0$ as $t$ tends to $\infty$. This implies that the difference of two Fredholm determinant in \eqref{eq:TailKernelDecay} is exponentially decaying to $0$.
\end{proof}
   We have decomposed the kernel $K^b_{\zeta_t}$ in \eqref{eq:PreSubDecom}. Apart from that, there is another kind of canonical decomposition of $K^b_{\zeta_t}$ based on the division of the integral in \eqref{eq:KernelIntMedForm4} into two parts. For instance, let us consider
   \begin{align}\label{eq:CrucialDecomp}
   K^{b}_{\zeta_t}(w,w^\prime)= K^{(1)}_{\zeta_1}(w,w^\prime)+K^{(2)}_{\zeta_t}(w,w^\prime)
   \end{align}
   where for $l=1,2$
   \begin{align}\label{eq:KernelTwoPartDescription}
   K^{(1)}_{\zeta_t}:= \frac{1}{2i}\sum_{k=-\infty}^{\infty}\int_{\mathcal{L}^{(l)}_{v}}&
  \frac{\exp\left(t(\Lambda(v)-\Lambda(w))+s\sigma_\nu t^{1/3}(\log (v)-\log (w))\right)}{\sin
  \left(\frac{\pi}{\log (b)}(\log(v)-\log(w))+2\pi ik\right)}\\
  & \times \frac{\exp\left(s\left(\log(v)-\log(w)\right)\right)}{\log(b)(v-w^\prime)} \times \frac{\left(\alpha b^{-1}v;b\right)_\infty}{\left(\alpha b^{-1}w;b\right)_\infty}\times \frac{dv}{v}.
   \end{align}
  In what follows, we show that the contribution of the kernel $K^{(2)}_{\zeta_2}$ in $\det(1+K^b_{\zeta_t})_{L^2(\mathcal{L}_w)}$ fades way as $t$ grows up.

  \bl\label{KernelRemainingPart}
   There exist $c,C>0$ such that
\begin{align}\label{eq:KernelResuidualBoundAC}
|K^{(2)}_{\zeta_t}(w,w^\prime)|\leq \frac{\epsilon^{-1}}{\log(1+\epsilon)}\exp\left(-c(t+|w|)+s\sigma_\nu t^{1/3}|w|\vee |w|^{-1}\right)
\end{align}
   holds for all $w, w^\prime \in \mathcal{L}^{(1)}_w$. Thus, one would get
   \begin{align}\label{eq:KernelRsPartDecay}
  \lim_{t\to \infty} \left|\det(1+K^b_{\zeta_t})_{L^2(\mathcal{L}^{(1)}_w)}-\det(1+K^{(1)}_{\zeta_t})_{L^2(\mathcal{L}^{(1)}_w)}\right|= 0.
   \end{align}
  \el
    \begin{proof}
   Like in the proof of \eqref{eq:KernelBound}, \eqref{eq:KernelResuidualBoundAC} can be proved using the bound on each of the components of the integrand of $K^{(2)}_{\zeta_t}$. Bound for the exponential term in the integrand would be exactly same as what had been derived in \eqref{eq:BoundsOnACExpTerm}. For the series over inverse sin function, we use the bounds in \eqref{eq:BoundOnSinFuncAC} and \eqref{eq:BoundOnSinACWhenK0} together to get \eqref{eq:BoundOnSeriesAC}. We must point out here that minimum value of $|\log|v/w||$ for any $v\in \mathcal{L}^{(2)}_v$ and $w^\prime\in\mathcal{L}^{(1)}_{w}$ is again bounded below by $\log(1+\epsilon)$, thus, verifying \eqref{eq:BoundOnSinACWhenK0} in the present scenario. Bound on the infinite product term in \eqref{eq:InfiniteProdBoundAC} will again be the same. Minimum distance between any two points on the contour $\mathcal{L}^{(2)}_v$ and $\mathcal{L}^{(1)}_w$ is $\epsilon$. Thus, bound on $|1/(v-w^\prime)|$ for the present scenario is only $\epsilon^{-1}$. These all together establish the inequality in \eqref{eq:KernelResuidualBoundAC}. Using the expansion of the Fredholm determinant as in \eqref{eq:SubtractionSurgery} and bound on difference between two determinants like in \eqref{eq:matrixdetdiff}, one can complete the proof of \eqref{eq:KernelRsPartDecay} given the bound on the kernel $K^{(2)}_{\zeta_t}$. As the bounds in \eqref{eq:KernelBound} and \eqref{eq:KernelResuidualBoundAC} are almost same, thus, it implies that difference of two Fredholm determiant in \eqref{eq:KernelRsPartDecay} also decays exponentially fast  as $t$ escapes to $\infty$.
    \end{proof}
  Now, we make a substitutions of the variables $v,w$ and $w^\prime$ in the $\epsilon$-neighborhood of $\varrho$. In particular, we transform
  \begin{align}\label{eq:Localization}
  v\mapsto \varrho+ t^{-1/3}\varrho \sigma^{-1}_\nu \tilde{v}, \quad w\mapsto \varrho+ t^{-1/3}\varrho \sigma^{-1}_\nu \tilde{w},\quad w^\prime\mapsto \varrho+ t^{-1/3}\varrho \sigma^{-1}_\nu \tilde{w^\prime} .
\end{align}
To the effect of these substitutions, transformed kernel $K^{(1)}_{\zeta_t}$ looks like
  \begin{align}
  K^{(1)}_{\zeta_t}= \frac{t^{-1/3}\varrho \sigma^{-1}_\nu}{2i}\sum_{k=-\infty}^{\infty}\int_{\mathcal{L}^{(1)}_{\tilde{v}}}&
  \frac{\exp\left(t(\Lambda(1+t^{-1/3}\sigma^{-1}_\nu \tilde{v})-\Lambda(1+t^{-1/3}\sigma^{-1}_\nu \tilde{w}))\right)}{\sin
  \left(\frac{\pi}{\log (b)}(\log(1+t^{-1/3}\sigma^{-1}_\nu \tilde{v})-\log(1+t^{-1/3}\sigma^{-1}_\nu \tilde{w}))+2\pi ik\right)}\nonumber\\
  & \times \frac{\exp\left(s\sigma_\nu t^{1/3}\left(\log(1+t^{-1/3}\sigma^{-1}_\nu \tilde{v})-\log(1+t^{-1/3}\sigma^{-1}_\nu \tilde{w})\right)\right)}{\log(b)(\tilde{v}-\tilde{w}^\prime)}\nonumber\\
  & \times \frac{\left(\alpha b^{-1}(\varrho+t^{-1/3}\varrho\sigma^{-1}_\nu \tilde{v});b\right)_\infty}{\left(\alpha b^{-1}(\varrho+t^{-1/3}\varrho\sigma^{-1}_\nu \tilde{w});b\right)_\infty}\times \frac{d\tilde{v}}{\varrho+t^{-1/3}\varrho\sigma^{-1}_\nu \tilde{v}}.\label{eq:KernelIntMedForm5LocVer}
\end{align}
  We must point out here that $\mathcal{L}^{(1)}_{\tilde{v}}, \mathcal{L}^{(2)}_{\tilde{v}}$ have been used to denote the contours for the transformed variable $\tilde{v}$. Similar change in notations for the contours of $\tilde{w}$ and $\tilde{w}^\prime$ will also be made in the subsequent discussion.
Next, we state two lemmas which essentially prove part $(a)$ of Proposition~\ref{FinalAsympProp}.
\bl\label{KernelLimitLemma}
 For all $\tilde{w},\tilde{w}^{\prime}\in \mathcal{L}^{(1)}_{\tilde{w}}$, we have
 \begin{align}\label{eq:LocalKernelTailEstimate}
 |K^{(1)}_{\zeta_t}(\tilde{w},\tilde{w}^\prime)|\leq C\exp\left(-c_1|\mbox{Re}(\tilde{w}^3)|-s\sigma_\nu|w|\right)
\end{align}
 for some constants $ C, c_1,c_2>0$. Moreover, the following convergence
 \begin{align}\label{eq:KernelConv}
 \lim_{t\to \infty}K^{(1)}_{\zeta_t}(\tilde{w},\tilde{w}^\prime)=\frac{1}{2\pi i}\int_{\infty e^{-i\pi/3}}^{\infty e^{i\pi/3}}\frac{\exp\left(\frac{\tilde{w}^3}{3}-\frac{\tilde{v}^3}{3}+s\tilde{v}-s\tilde{w}\right)d\tilde{v}}{(\tilde{v}-\tilde{w}^\prime)(\tilde{v}-\tilde{w})}
 \end{align}
 holds.
\el
\begin{proof}
For all $\tilde{w}\in \mathcal{L}^{(1)}_{\tilde{w}}$ and $\tilde{v}\in \mathcal{L}^{(1)}_{\tilde{v}}$, we have
\begin{align}\label{eq:LambdaLimit}
t\left(\Lambda(\varrho+t^{-13}\varrho\sigma^{-1}_\nu\tilde{w})-\Lambda(\varrho+t^{-13}\varrho\sigma^{-1}_\nu\tilde{w})\right)&=\frac{\tilde{w}^3}{3}-\frac{\tilde{v}^3}{3}\\&+t\left(\mathfrak{R}(t^{-1/3}\varrho\sigma^{-1}_\nu\tilde{w})-\mathfrak{R}(t^{-1/3}\varrho\sigma^{-1}_\nu\tilde{v})\right)
\end{align}
where $t(|\mathfrak{R}(t^{-1/3}\varrho\sigma^{-1}_\nu\tilde{w})|+|\mathfrak{R}(t^{-1/3}\varrho\sigma^{-1}_\nu\tilde{v})|)$ converges to $0$ as $t\to \infty$. One can note that $\epsilon$ in Definition~\ref{wvContourDef} has been chosen in such a way that for $\tilde{v}\in \mathcal{L}^{(1)}_{\tilde{v}}$ and $ \tilde{w}\in \mathcal{L}^{(1)}_{\tilde{w}}$, we have
\begin{align}\label{eq:RemainderBound}
t(|\mathfrak{R}(t^{-1/3}\varrho\sigma^{-1}_\nu\tilde{w})|+|\mathfrak{R}(t^{-1/3}\varrho\sigma^{-1}_\nu\tilde{v})|)\leq c\left(|\tilde{v}|^3+|\tilde{w}|^3\right).
\end{align}
for some small constant $c>0$.
Further, we know $\mbox{Re}(w^3-v^3)<0$ for $w\in \mathcal{L}^{(1)}_w$ and $v\in \mathcal{L}^{(1)}_v$. That implies real part of the right side in \eqref{eq:LambdaLimit} is bounded above by $-c_1\left(|\mbox{Re}(\tilde{v})^3|+|\mbox{Re}(\tilde{w})^3|\right)$ for some constant $c_1>0$. Among the other components of the integrand in \eqref{eq:KernelIntMedForm5LocVer}, we have
\bei
\ii
\begin{align}
\left|\exp\left(s\sigma_\nu t^{1/3}\left(\log(1+t^{-1/3}\sigma^{-1}_\nu \tilde{v}) -\log(1+t^{-1/3}\sigma^{-1}_\nu \tilde{w})\right)\right)\right|
\leq \exp\left(s\sigma_\nu (|v|-|w|)\right)\\
\label{eq:LinearcompBoundAC}
\end{align}
and, as $t\to \infty$,
\begin{align}\label{eq:LinearComp}
\exp\left(s\sigma_\nu t^{1/3}\left(\log(1+t^{-1/3}\sigma^{-1}_\nu \tilde{v})-\log(1+t^{-1/3}\sigma^{-1}_\nu \tilde{w})\right)\right)\to \exp(s(\tilde{v}-\tilde{w})),
\end{align}
\ii
\begin{align}\label{eq:SinCompBoundAC}
\left|\frac{t^{1/3}}{\sin\left(\frac{\pi}{\log(b)}\left(\log(1+t^{-1/3}\sigma^{-1}_\nu \tilde{v})-\log(1+t^{-1/3}\sigma^{-1}_\nu \tilde{w})\right)\right)}\right|\leq \frac{\sigma_\nu|\log(b)|}{2|\tilde{v}-\tilde{w}|}
\end{align}
and,
\begin{align}
\frac{\log(b)}{\pi t^{-1/3}\varrho\sigma^{-1}_{\nu}}\sin\left(\frac{\pi}{\log(b)}\left(\log(1+t^{-1/3}\sigma^{-1}_\nu \tilde{v})-\log(1+t^{-1/3}\sigma^{-1}_\nu \tilde{w})\right)\right)\stackrel{t\to\infty}{\longrightarrow} \tilde{v}-\tilde{w},\\
\label{eq:SinComp}
\end{align}
\ii
\begin{align}
t^{-1/3} \sum_{\substack{k=-\infty\\ k\neq 0}}^\infty\frac{1}{\sin\left(\frac{\pi}{\log(b)}\left(\log(1+t^{-1/3}\sigma^{-1}_\nu \tilde{v})-\log(1+t^{-1/3}\sigma^{-1}_\nu \tilde{w})\right)+2\pi ik\right)}\stackrel{t\to\infty}{\longrightarrow} 0,\\
\label{eq:SinSumComp}
\end{align}
\ii
\begin{align}\label{eq:InfiniteProductComp}
\frac{(\alpha b^{-1}(\varrho+ \varrho \sigma^{-1}_\nu t^{-1/3}\tilde{v});b)_\infty}{(\alpha b^{-1}(\varrho+ \varrho \sigma^{-1}_\nu t^{-1/3}\tilde{w});b)_\infty} \to 1,
\end{align}
whenever $t\to \infty$. Convergence in \eqref{eq:LinearComp} and \eqref{eq:SinComp} are easy to see. To prove the inequality in \eqref{eq:LinearcompBoundAC}, note that we have
\[\mbox{Re}\left(t^{1/3}(\log(1+t^{-1/3}\sigma^{-1}_\nu \tilde{v}) -\log(1+t^{-1/3}\sigma^{-1}_\nu \tilde{w}))\right)\leq \sigma^{-1}_\nu(|v|-|w|)\]
 for all large $t$. For showing \eqref{eq:SinCompBoundAC}, we used the fact that $|\sin(x)|\geq 2x/\pi$ for all $x\in [-\pi/2,\pi/2]$. For large values of $t$, $\log(1+t^{-/3\sigma^{-1}_\nu\tilde{v}})$ and $\log(1+t^{-/3\sigma^{-1}_\nu\tilde{v}})$ will be close enough to $0$. This validates the bound in \eqref{eq:SinCompBoundAC}. In the third case, when $k\neq 0$, $\sin(.+2\pi ik)$ grows exponentially fast as $|k|$ increases. This makes the series in \eqref{eq:SinSumComp} convergent and hence, when multiplied with $t^{-1/3}$, it goes to $0$. To see \eqref{eq:InfiniteProductComp}, notice that $\alpha b^{-1} \varrho <1$ thanks to the condition $\rho >1-(\nu (1-b))^{-1/2}$. Therefore,
\begin{align}\label{eq:InfiniteProductLimitProof}
 \frac{(\alpha b^{-1}(\varrho+ \varrho \sigma^{-1}_\nu t^{-1/3}\tilde{v});b)_m}{(\alpha b^{-1}(\varrho+ \varrho \sigma^{-1}_\nu t^{-1/3}\tilde{v});b)_m} \stackrel{t\to \infty}{\longrightarrow} 1.
\end{align}
 uniformly over $m$. Moreover, using some simple inequalities like $|(1+z)|\leq \exp(|z|)$ for any $z\in \CC$ and $(1+z)^{-1}\leq \exp(|z|)$ for $z\in \CC$ with $|\mbox{Re}(z)|< 1$, it is also easy to see following crude bound
\begin{align}\label{eq:BoundOnInfiniteProdAC}
\frac{(\alpha b^{-1}(\varrho+ \varrho \sigma^{-1}_\nu t^{-1/3}\tilde{v});b)_\infty}{(\alpha b^{-1}(\varrho+ \varrho \sigma^{-1}_\nu t^{-1/3}\tilde{v});b)_\infty} \leq \exp(\sigma^{-1}_\nu t^{-1/3}(|v|+|w|)).
\end{align}
 Thus, the integrand in \eqref{eq:KernelIntMedForm5LocVer} is uniformly bounded by an integrable function \[C\frac{1}{1+|\tilde{v}|}\exp(-c_1(|\mbox{Re}(\tilde{v}^3)|+|\mbox{Re}(\tilde{w})^3|)+c_2(|v|-|w|)).\] This integrable envelop proves the bound on the kernel in \eqref{eq:LocalKernelTailEstimate}. Further, using dominated convergence theorem, we get the convergence in \eqref{eq:KernelConv} whereas \eqref{eq:LinearComp}-\eqref{eq:InfiniteProductComp} shows the limiting value.
\ee
\end{proof}

\textsc{Proof of Proposition~\ref{FinalAsympProp}$(a)$. }
 Using Lemma~\ref{TailBoundLemma} and \ref{KernelRemainingPart}, we have
 \begin{align}\label{eq:LastEquivalenceOfFredDetAC}
 \lim_{t\to \infty}\det (1+K^b_{\zeta_t})_{L^2(C_r)}= \lim_{t\to \infty}\det (1+K^b_{\zeta_t})_{L^2(L^{(1)}_{\tilde{w}})}.
 \end{align}
Furthermore, the kernel $K^{(1)}_{\zeta_t}$ has an integrable upper bound as shown in \eqref{eq:LocalKernelTailEstimate} for all $\tilde{w},\tilde{w}^\prime\in \mathcal{L}^{(1)}_{\tilde{w}}$. Thus, using Corollary~\ref{DominatedConvCor}, one can establish that the limit of the Fredholm determinant $\det (1+K^b_{\zeta_t})_{L^2(L^{(1)}_{\tilde{w}})}$ is same as the  Fredholm determinant of the limiting kernel in \eqref{eq:KernelConv}. This completes the proof.

\subsection{Case $\rho = 1-\left(\nu(1-b)\right)^{-1/2}$, $F_{GOE}^2$ Fluctuations.}
In this scenario, the double critical point $\varrho$ of the function $\Lambda$ becomes equal to $b\alpha^{-1}$. Thus, the deformation of the contour of $w$ to $\mathcal{L}_w$ is no more viable because in Definition~\ref{wvContourDef}, $\mathcal{L}_w$ is described to pass through $\varrho$ whereas the contour of $w$ must avoid the pole at $b\alpha^{-1}$ as stated in Theorem~\ref{FredRep}. Furthermore, the result in \eqref{eq:InfiniteProductComp} is also not going hold. In what follows, we make a mild change in the definition of contours $\mathcal{L}_w$ and $\mathcal{L}_v$ to remedy those issues.

\bd\label{Criticalwvcontourdef}
Recall all the notations from Definition~\ref{wvContourDef}. To avoid the pole at $b\alpha^{-1}$, we redefine $\mathcal{L}^{(1)}_{w}$ to go from $\varrho-t^{-1/3}\delta\varrho\sigma^{-1}_\nu+\epsilon e^{-i\pi/3}$ to $\varrho-t^{-1/3}\delta\varrho\sigma^{-1}_\nu$ to $\varrho-t^{-1/3}\delta\varrho\sigma^{-1}_\nu+\epsilon e^{i\pi/3}$. Similarly, we shift the contour $\mathcal{L}^{(1)}_{v}$ to the left in the complex plane by $t^{-1/3}\delta\varrho\sigma^{-1}_\nu$. For instance, in this new settings, $\mathcal{L}^{(1)}_{v}$ extends from $\varrho-t^{-1/3}(1+\delta)\varrho\sigma^{-1}_\nu+\epsilon e^{-i2\pi/3}$ to $\varrho-t^{-1/3}\delta\varrho\sigma^{-1}_\nu$ to $\varrho-t^{-1/3}(1+\delta)\varrho\sigma^{-1}_\nu+\epsilon e^{i2\pi/3}$. Here, $\delta$ can be any positive real number. Furthermore, $\epsilon$ must be chosen is such a way that following conditions, given as
\be
\ii both end points $\varrho-t^{-1/3}\delta\varrho\sigma^{-1}_\nu+\epsilon e^{-i\pi/3}$ and $\varrho-t^{-1/3}\delta\varrho\sigma^{-1}_\nu+\epsilon e^{i\pi/3}$ of $\mathcal{L}^{(1)}_w$ lie outside the level curve $\mathcal{L}_2$ (see, Theorem~\ref{LevelCurveCharacteristic} for definition),
\ii $\mbox{Re}(\Lambda(w))<\mbox{Re}(\Lambda(w))$ for all $v\in\mathcal{L}^{(1)}_v$ and $\mathcal{L}^{(1)}_w$,
 \ii $|\mathfrak{R}((z-\varrho))|<c|z-\varrho|^3$ for any $z$ in an $\epsilon$-neighbourhood of  $\varrho$,
\ee
 are satisfied.  The residual parts, namely $\mathcal{L}^{(2)}_w$ and $\mathcal{L}^{(2)}_v$ are defined exactly the same way as in Definition~\ref{wvContourDef} except the necessary adjustments of their origins and end points. For instance, $\mathcal{L}^{(2)}_w$ now starts from $\varrho-t^{-1/3}\delta\varrho\sigma^{-1}_\nu+\epsilon e^{-i\pi/3}$ and ends at $\varrho-t^{-1/3}\delta\varrho\sigma^{-1}_\nu+\epsilon e^{i\pi/3}$.
\ed
 Thanks to the new definitions above, we can deform of the contour of $v$ and $w$ to $\mathcal{L}_v$ and $\mathcal{L}_w$ without crossing any singularities. Moreover, one can note that $b\alpha^{-1}$ lies outside contour $\mathcal{L}_w$. To this end, one can prove the results like  Lemma~\ref{TailBoundLemma} and Lemma~\ref{KernelRemainingPart} almost in the same way in the present context except possibly with a different treatment for the quantity $\left(\alpha b^{-1}v;b\right)_\infty\left(\alpha b^{-1}w;b\right)^{-1}_\infty$ in \eqref{eq:InfiniteProdBoundAC}. In what follows, we carry out such changes.
\bl\label{InfiniteProductAtCriticality}
Assume $\rho=1-\left(\nu(1-b)\right)^{-1/2}$.
 There exist $c^\prime, C^\prime >0$ such that whenever $\max\{|w-\varrho|,|v-\varrho|\}<c^\prime$, then
\begin{align}\label{eq:InfiniteProdCritBound}
\left|\frac{\left(\alpha b^{-1}v;b\right)_\infty}{\left(\alpha b^{-1}w;b\right)_\infty}\right|<C^\prime \left|\frac{v-\varrho}{w-\varrho}\right|.
\end{align}
 Furthermore, consider the substitutions in \eqref{eq:Localization}. Denote the contours $\mathcal{L}^{(1)}_v$ and $\mathcal{L}^{(1)}_w$ after the substitution by $\mathcal{L}^{(1)}_{\tilde{v}}$ and $\mathcal{L}^{(1)}_{\tilde{w}}$ respectively. Then, for any $\tilde{v}\in \mathcal{L}^{(1)}_{\tilde{v}}$ and $\tilde{w}\in \mathcal{L}^{(1)}_{\tilde{w}}$, we have
\begin{align}\label{eq:InfiniteProdCritLimit}
\lim_{t\to \infty} \frac{\left(\alpha b^{-1}(\varrho+t^{-1/3}\varrho\sigma^{-1}_\nu \tilde{v});b\right)_\infty}{\left(\alpha b^{-1}(\varrho+t^{-1/3}\varrho\sigma^{-1}_\nu \tilde{w});b\right)_\infty}=\frac{\tilde{v}}{\tilde{w}}.
\end{align}
\el

\begin{proof}

 To prove \eqref{eq:InfiniteProdCritBound}, let us write
 \begin{align}\label{eq:InfiniteprodPostAnalysis}
 \frac{\left(\alpha b^{-1}v;b\right)_\infty}{\left(\alpha b^{-1}w;b\right)_\infty}= A\cdot B\cdot C
\end{align}
where
\[A:=\frac{v-\varrho}{w-\varrho},\quad B:= \frac{\prod_{j=1}^m (1-\alpha b^{j-1}v)}{\prod_{j=1}^m (1-\alpha b^{j-1}w)},\quad C:= \frac{(\alpha b^{m}v;b)_\infty}{(\alpha b^{m}w;b)_\infty}.\]
Further, $m$ is chosen in such  a way that the inequality $b^m\alpha\left(\varrho+c^\prime\right)<1/2$ is satisfied. To that effect, we have
\begin{align}\label{eq:ConTrolTailOfInfiniteProduct}
|C|&=\left|\exp\left(\sum_{j=1}^\infty \log\left(\frac{1+\alpha b^{m+j-1}v}{1+\alpha b^{m+j-1}v}\right)\right)\right|= \exp\left(\alpha b^m\mbox{Re}\left(\sum_{j=1}^\infty \frac{b^{j-1}(v-w)}{1+\psi(v,w,j)}\right)\right)\nonumber\\
&\leq \exp\left(\alpha b^m\left(\sum_{j=1}^\infty \frac{b^{j-1}|v-w|}{1-\frac{1}{2}}\right)\right)\leq \exp\left(\frac{2\alpha b^m|v-w|}{1-b} \right).
\end{align}
In the second equality above, we have used the mean value theorem for the complex variables. Thus, $\psi(v,w,j)$ denotes a complex number in the line joining $\alpha b^{m+j} v$ and $\alpha b^{m+j} w$. Consequently, one have
\[\mbox{Re}\left(\frac{v-w}{1+\psi(v,w,j)}\right)\leq \frac{|v-w|}{1-|\psi(v,w,j)|}\leq \frac{|v-w|}{1-b^m\alpha\left(\varrho+c^\prime\right)}.\]
This shows the inequalities in \eqref{eq:ConTrolTailOfInfiniteProduct}. Lastly, $B$ involves product of finitely many terms which are bounded above. Hence, this completes the proof of \eqref{eq:InfiniteProdCritBound}.

Further, one can write
\begin{align}\label{eq:InfiniteProd}
\frac{\left(\alpha b^{-1}(\varrho+t^{-1/3}\varrho\sigma^{-1}_\nu \tilde{v});b\right)_\infty}{\left(\alpha b^{-1}(\varrho+t^{-1/3}\varrho\sigma^{-1}_\nu \tilde{w});b\right)_\infty}=\frac{(-t^{-1/3}\varrho\sigma^{-1}_\nu\tilde{v})}{(-t^{-1/3}\varrho\sigma^{-1}_\nu\tilde{w})}\times\frac{\left(\alpha (\varrho+t^{-1/3}b\varrho\sigma^{-1}_\nu \tilde{v});b\right)_\infty}{\left(\alpha (\varrho+t^{-1/3}b\varrho\sigma^{-1}_\nu \tilde{w});b\right)_\infty}.
\end{align}
Note that $b^n\alpha\varrho$ is bounded away from $1$ for all non-negative powers of $b$. Thus, the second term in the product of the right side in \eqref{eq:InfiniteProd} converges to $1$ as $t\to \infty$. This shows the limit in \eqref{eq:InfiniteProdCritLimit}.
\end{proof}

 Now, we note down two results in the same spirit of Lemma~\ref{TailBoundLemma} and Lemma~\ref{KernelRemainingPart}. Proofs are exactly same except one has to use the bound in \eqref{eq:InfiniteProdCritBound} instead of \eqref{eq:InfiniteProdBoundAC}. For brevity, we omit the proofs.
\bl\label{TailBoundLemmaAtCriticality}
 There exists some $c, C>0$ and $t_0>0$, such that for all $t\geq t_0$
\begin{align}\label{eq:KernelBoundAtCriticality}
|K^{b}_{\zeta_t}(w,w^\prime)|\leq C \frac{1}{|w-\varrho|}\frac{\max\{\epsilon^{-1}, \varrho^{-1} \sigma_\nu t^{1/3}\}}{\log(1+\epsilon)}\exp\left(-c(t+|w|)+s\sigma_\nu t^{1/3}|w|\vee |w|^{-1}\right)
\end{align}
where $w\in \mathcal{L}^{(2)}_w$, $w^\prime\in \mathcal{L}_{w}$ and $\epsilon$ is defined in Definition~\ref{wvContourDef}. Moreover, the following convergence
\begin{align}\label{eq:TailKernelDecayAtCriticality}
\lim_{t\to \infty}\left|\det(1+K^b_{\zeta_t})_{L^2(\mathcal{L}_w)}-\det(1+K^b_{\zeta_t})_{\mathcal{L}^{(1)}_w}\right| = 0
\end{align}
also holds in the present scenario.
\el
\bl\label{KernelRemainingPartAtCriticality}
   There exist $c,C>0$ such that
\begin{align}\label{eq:KernelResuidualBoundAtCriticality}
|K^{(2)}_{\zeta_t}(w,w^\prime)|\leq \frac{1}{|w-\varrho|}\frac{\epsilon^{-1}}{\log(1+\epsilon)}\exp\left(-c(t+|w|)+s\sigma_\nu t^{1/3}|w|\vee |w|^{-1}\right)
\end{align}
   holds for all $w, w^\prime \in \mathcal{L}^{(1)}_w$. Thus, one would get
   \begin{align}\label{eq:KernelRsPartDecayAtCriticality}
  \lim_{t\to \infty} \left|\det(1+K^b_{\zeta_t})_{L^2(\mathcal{L}^{(1)}_w)}-\det(1+K^{(1)}_{\zeta_t})_{L^2(\mathcal{L}^{(1)}_w)}\right|= 0.
   \end{align}
  \el

\bl\label{KernelLimitLemmaAtCriticality}
 Recall all the notations from the previous section and Definition~\ref{Criticalwvcontourdef}. For all $\tilde{w},\tilde{w}^\prime\in \mathcal{L}^{(1)}_{\tilde{w}}$, we have
 \begin{align}\label{eq:BoundByIntegrableFunc}
 |\tilde{K}^{(1)}_{\zeta_t}(\tilde{w},\tilde{w}^\prime)|\leq C\frac{1}{\delta}\exp\left(-c_1|\mbox{Re}(\tilde{w}^3)|-s\sigma_\nu|w|\right)
\end{align}
for some positive constants $c_1,c_2,C$ and
\begin{align}\label{eq:CriticalKernelLimit}
\lim_{t\to \infty}\tilde{K}^{(1)}_{\zeta_t}(\tilde{w},\tilde{w}^\prime) = \frac{1}{2\pi i}\int_{-(1+\delta)+e^{-i2\pi/3}}^{-(1+\delta)+e^{i2\pi/3}}\frac{\exp\left(\frac{\tilde{w}^3}{3}-\frac{\tilde{v}^3}{3}+s(\tilde{v}-\tilde{w})\right)}{(\tilde{v}-\tilde{w})(\tilde{v}-\tilde{w}^\prime)}\times \frac{\tilde{v}d\tilde{v}}{\tilde{w}}.
\end{align}
\el
\begin{proof}
 To get the bound in \eqref{eq:BoundByIntegrableFunc}, note that integrand in \eqref{eq:KernelIntMedForm5LocVer} is  bounded above by
\[C\frac{1}{(1+|\tilde{v}|)}\exp(-c_1(|\mbox{Re}(\tilde{v}^3)|+|\mbox{Re}(\tilde{w})^3|)+c_2(|v|-|w|))\left|\frac{\left(\alpha b^{-1}(\varrho+t^{-1/3}\varrho\sigma^{-1}_\nu \tilde{v});b\right)_\infty}{\left(\alpha b^{-1}(\varrho+t^{-1/3}\varrho\sigma^{-1}_\nu \tilde{w});b\right)_\infty}\right|.\]
We have illustrated this fact in Lemma~\ref{KernelLimitLemma}. Thus, we omit further details on it. Let us  recall that $\mathcal{L}_w$ is shifted to the left by $t^{-1/3}\delta\varrho\sigma^{-1}_\nu$ from $\varrho$. To this effect, $|\tilde{w}|\gtrsim\delta$ and henceforth, the ratio of the infinite products in \eqref{eq:InfiniteProdCritLimit} gets bounded above by $C^\prime|\tilde{v}|/\delta$ thanks to \eqref{eq:InfiniteProdCritBound}. This accounts for the term $1/\delta$ in the bound on the right side in \eqref{eq:BoundByIntegrableFunc}. Proof of \eqref{eq:CriticalKernelLimit} follows from the similar arguments like in Lemma~\ref{KernelLimitLemma} except the limit in \eqref{eq:InfiniteProductComp} should be replaced by the result \eqref{eq:InfiniteProdCritLimit} of Lemma~\ref{InfiniteProductAtCriticality}.
\end{proof}
 To summarize, all these results above prove
\begin{align}\label{eq:CrtitDetLimitNotExact}
\lim_{t \to \infty}\det(1+K^b_{\zeta_t})_{L^2(C_r)}=\det(1+K_{crit})_{L^2(\bar{\mathcal{L}}_{\tilde{w}})}
\end{align}
where
\begin{align}\label{eq:KcritLimit}
K_{crit}(\tilde{w},\tilde{w}^\prime):=\frac{1}{2\pi i}\int_{-(1+\delta)+e^{-i2\pi/3}}^{-(1+\delta)+e^{i2\pi/3}}\frac{\exp\left(\frac{\tilde{w}^3}{3}-\frac{\tilde{v}^3}{3}+s(\tilde{v}-\tilde{w})\right)}{(\tilde{v}-\tilde{w})(\tilde{v}-\tilde{w}^\prime)}\times \frac{\tilde{v}d\tilde{v}}{\tilde{w}}
\end{align}
and $\bar{\mathcal{L}}_{\tilde{w}}$ denotes a piecewise linear contour which extends linearly from $-\delta+\infty e^{-i\pi/3}$ to $-\delta$ and from there again linearly to $-\delta+\infty e^{i\pi/3}$. Hence, the proof of Proposition~\ref{FinalAsympProp} $(b)$ follows.

\subsection{Case $\rho<1-(\nu(1-b))^{-1/2}$, Gaussian Fluctuation}
 In this section, we prove part $(c)$ of Proposition~\ref{FinalAsympProp}. Main reason behind a separate analysis is that whenever $\rho<1-(\nu(1-b))^{-1/2}$, contours like $\mathcal{L}_v$ and $\mathcal{L}_w$ in Definition~\ref{wvContourDef} or even their modifications in Definition~\ref{Criticalwvcontourdef} contain $b\alpha^{-1}$ inside. To circumvent this illegal inclusion, we adjust the choice of $m_\nu$ in the function $\Lambda$ (see \eqref{eq:LambdaFuncDefine}) in such a way that $b\alpha^{-1}$ turns out to be a new critical point of $\Lambda$. To begin with, the form of the kernel $K^b_{\zeta_t}$ is given by
\begin{align}\label{eq:KernelIntMedForm6}
K^b_{\zeta_t}(w,w^\prime)=\frac{1}{2 i\log(b)}\sum_{k=-\infty}^{\infty}\int_{\mathcal{C}_{b^R|w|}}&\frac{\exp(t(\Lambda(w)-\Lambda(v))+t^{1/2}s\tilde{\sigma}_\nu(\log (v)-\log(w)))dv}{(v-w^\prime)\sin\left(\frac{\pi}{\log (b)}(\log(v/w))+2\pi ik\right)}\nonumber\\
  &\times \frac{(\alpha b^{-1}v;b)_\infty }{(\alpha b^{-1}w;b)_\infty }\times \frac{dv}{v}.
\end{align}
 where  $\zeta_t$ and $b^h w$ are substituted with $b^{-\tilde{m}_\nu t+s\tilde{\sigma}_\nu t^{1/2}}$ and $v$ respectively. Under the these recalibration, the explicit form of the function $\Lambda$ is given by
\begin{align}\label{eq:LambdaReWrite}
\Lambda(z)= \frac{z}{b(1-b)}-\nu \log(1+zb^{-1})+\tilde{m}_\nu\log(z).
\end{align}
  In what follows, we perform some calculations \emph{de rigueur} the asymptotic analysis.
 Plugging $\tilde{m}_\nu = \frac{1}{1+\alpha}\nu - \frac{1}{\alpha(1-b)}$ in \eqref{eq:LambdaReWrite}, one can write
 \begin{align}\label{eq:LambdaDeriveReWrite}
 \Lambda^\prime(z)=\frac{(z-b\alpha^{-1})(z-\varrho^\prime)}{zb(z+b)(1-b)}.
 \end{align}
   where $\varrho^\prime := b\left(\nu(1-\rho)(1-b)-1\right)$. Under the condition $\rho<1-(\nu(1-b))^{-1/2}$, it can be noted that $\varrho^\prime>\varrho$. To this end, we have
   \begin{align}\label{eq:LambdaDDeriveReWrite}
   \Lambda^{\prime\prime}(b \alpha^{-1}) = -\frac{(1-\rho)^2\nu(1-b)-1}{b^2\rho(1-b)}.
\end{align}
Further, using Taylor's expansion, one would get
\begin{align}\label{eq:Taylor'sExpan}
\Lambda(z) = \Lambda(b\alpha^{-1})-\frac{1}{2}\left(\frac{\tilde{\sigma}_{\nu}}{b\alpha^{-1}}\right)^2 (z-b\alpha^{-1})^2+\tilde{\mathfrak{R}}(z-b\alpha^{-1})
\end{align}
where $|z-b\alpha^{-1}|^2\tilde{\mathfrak{R}}(z-b\alpha^{-1})$ converges to $0$ whenever $z\to b\alpha^{-1}$. In the following lemma, we illustrate on the properties of the level curves $\mbox{Re}(\Lambda(z))=\mbox{Re}(\Lambda(b\alpha^{-1}))$.

\bl\label{LevelCurveBelowCriticality}
 There exists two simple closed contours $\mathcal{L}_1$ and $\mathcal{L}_2$ such that for all $z\in \mathcal{L}_1 \cup \mathcal{L}_2$, we have $\mbox{Re}(\Lambda(z))=\mbox{Re}(\Lambda(b\alpha^{-1}))$.  They further satisfy the following properties.
 \ba
 \ii Both the contours originate at the point $b\alpha^{-1}$ and encircle around $0$ to meet again at their origin. For instance, $\mathcal{L}_1$ ($\mathcal{L}_2$) in its way back to $b\alpha^{-1}$ crosses the real line for the second time at a point $d_1\in (-b,0)$ ($d_2\in (-\infty, -b)$).
 \ii Contour $\mathcal{L}_1\backslash \{b\alpha^{-1}\}$ is completely contained in the interior of  the region bounded by $\mathcal{L}_2$.
 \ii Contour $\mathcal{L}_1$ ($\mathcal{L}_2$) leaves $x$ axis from the point $b\alpha^{-1}$ at an angle $3\pi/4$ ($\pi/4$) and meet again at angle $-3\pi/4$ ($-\pi/4$).
 \ee
 Apart from $\mathcal{L}_1$ and $\mathcal{L}_2$, there exists another level curve $\mathcal{L}_3$  of the equation $\mbox{Re}(\Lambda(z))=\mbox{Re}(\Lambda(b\alpha^{-1}))$ lying completely in the right half of the complex plane. Furthermore, $\mathcal{L}_3$ doesn't intersect $\mathcal{L}_1$ or $\mathcal{L}_2$ and escapes to infinity whenever $\mbox{Re}(z)\to \infty$.
\el
\begin{proof}
Proof of the existence of $\mathcal{L}_1,\mathcal{L}_2$ and $\mathcal{L}_3$ are similar to that in  Lemma~\ref{LevelCurveCharacteristic}. Moreover, property $(a)$ and $(b)$ can also be proved in the same way. To see part $(c)$, note that $\Lambda(z)-\Lambda(b\alpha^{-1})$ behaves as $-2^{-1}(\tilde{\sigma}_\nu \alpha b^{-1})^{2}\nu(z-b\alpha^{-1})^2$ in a neighborhood around $b\alpha^{-1}$. Thus, $\mbox{Re}\left(\Lambda(z)-\Lambda(b\alpha^{-1})\right)=0$  implies $\mbox{Arg}(z-b\alpha^{-1})=\pm\frac{\pi}{4}\pm \frac{2\pi }{2}$. This indicates part $(c)$. We are left to show that $\mathcal{L}_3$ doesn't intersect $\mathcal{L}_1$ or $\mathcal{L}_2$. To start with, note that $\mathcal{L}_3$ doesn't pass through $b\alpha^{-1}$. This is because there are only four possible choices of the angles ($\pm\frac{\pi}{4}\pm \frac{\pi}{2}$) that any level curve of the equation $\mbox{Re}\left(\Lambda(z)-\Lambda(b\alpha^{-1})\right)=0$ can make at $b\alpha^{-1}$. Thus, there can be at most two level curves which pass through the point $b\alpha^{-1}$.  In fact, $\mathcal{L}_3$ crosses the real line at a point further right of $b\alpha^{-1}$. To see this, notice that $\mbox{Re}\left(\Lambda(z)-\Lambda(b\alpha^{-1})\right)<0$ if $\mbox{Im}(z)=0$ and $b\alpha^{-1}<\mbox{Re}(z)<b\alpha^{-1}+\delta$ for some $\delta>0$. Additionally, $\mbox{Re}(\Lambda(z))$ increases to infinity with increasing $\mbox{Re}(z)$. Thus, there must exists another point on the real line to the right of $b\alpha^{-1}$ where once again $\mbox{Re}\left(\Lambda(z)-\Lambda(b\alpha^{-1})\right)=0$ holds. This shows the claim. Moreover, one can mimick almost the same proof as in Lemma~\ref{LevelCurveCharacteristic} to show that $\mathcal{L}_3$ doesn't intersect $\mathcal{L}_1$ or $\mathcal{L}_2$ at any other point. This completes the proof.
\end{proof}
We must point out here that $\mathcal{L}_2$ is no more star shaped (see Figure~\ref{fig:6}). This urges for the necessary changes in the definitions of the contours $\mathcal{L}_v$ and $\mathcal{L}_w$ in the present scenario. In what follows, we enunciate those changes formally.
\bd\label{wvCounyoursBelowCriticality}
 We define the linear segment $\mathcal{L}^{(1)}_{v}$ of $\mathcal{L}_v$ to extend linearly from $-\epsilon i+ \mbox{Int}_{b\alpha^{-1},2\delta}$ to $\epsilon i+\mbox{Int}_{b\alpha^{-1},2\delta}$ crossing the real line at $\mbox{Int}_{b\alpha^{-1}, 2\delta}$ where
 \begin{align}\label{eq:DefLinSegvIntsect}
 \mbox{Int}_{x,\tau}=x-t^{-1/3}\tau x\delta (\tilde{\sigma}_\nu)^{-1} \quad \text{ for any } x\in \CC, \tau \in \RR^+.
 \end{align}
On the flip side, we define the piecewise linear segment $\mathcal{L}^{(1)}_{w}$ to extend linearly from $\epsilon e^{-i\pi/6}+\mbox{Int}_{b\alpha^{-1},\delta}$  to $\mbox{Int}_{b\alpha^{-1},\delta}$ and from there to $\epsilon e^{i\pi/6}+\mbox{Int}_{b\alpha^{-1},\delta}$. At this junction, it is worthwhile to point out that any particular choice of $\epsilon$ must satisfies the following properties: $(a)$ $\epsilon e^{-i\pi/6}+\mbox{Int}_{b\alpha^{-1},\delta}$ and $\epsilon e^{i\pi/6}+\mbox{Int}_{b\alpha^{-1},\delta}$ should lie outside the contour $\mathcal{L}_2$, $(b)$ in the $\epsilon$ - neighborhood of the point $b\alpha^{-1}$, $|\mathfrak{R}(z-b\alpha^{-1})|$ in \eqref{eq:Taylor'sExpan} must be less than $c|z-b\alpha^{-1}|^2$ for small constant $c$. In particular, $c=1/4$ suitably fits into our analysis. Curved segments $\mathcal{L}^{(2)}_v$ and $\mathcal{L}^{(2)}_w$ are defined similarly as in Definition~\ref{wvContourDef}. For instance, $\mathcal{L}^{(2)}_v$ starts from the top end point of $\mathcal{L}^{(1)}_v$ and encircles around $0$ to meet at the bottom end point of $\mathcal{L}^{(1)}_v$. To be precise, it remains always inside of the contour $\mathcal{L}_2$. Unlike $\mathcal{L}^{(1)}_v$, starting from the top end point $\mathcal{L}^{(1)}_w$, $\mathcal{L}^{(2)}_w$ circles outside the contour $\mathcal{L}_2$ to meet the bottom end point of $\mathcal{L}^{(1)}_w$.
\ed

 \begin{figure}
\begin{center}
\begin{tikzpicture}[scale=0.8]
  \draw[line width = 0.25mm, <->] (-4,0) -- (8,0) node {};
  \draw[line width = 0.25mm, <->] (0,-4) -- (0,4) node {};
   \draw[domain=0:3.444, variable=\y, dashed, line width = 0.2mm] plot ({3.414-0.707*(\y-1.414)^2},{\y}) ;
   \draw[domain=0:3.444, variable=\y, dashed, line width = 0.2mm] plot ({3.414-0.707*(\y-1.414)^2},{-\y}) ;
   \draw[domain=0.5:2.747, variable=\x, dashed, line width = 0.2mm] plot ({-\x},{sqrt(14.041-(\x+1)^2)}) ;
   \draw[domain=0.5:2.747, variable=\x, dashed, line width = 0.2mm] plot ({-\x},{-sqrt(14.041-(\x+1)^2)}) ;
   \draw[line width = 0.2mm, dashed, -] (0.5,3.444) .. controls (0.25,3.529) and (-0.25,3.532) .. (-0.5,3.433);
   \draw[line width = 0.2mm, dashed, -] (0.5,-3.444) .. controls (0.25,-3.529) and (-0.25,-3.532) .. (-0.5,-3.433);
   \draw[domain=0:2, variable=\x, dashed, line width = 0.2mm] plot ({\x},{sqrt(6.25-(\x-0.5)^2)-2}) ;
   \draw[domain=0:1, variable=\x, dashed, line width = 0.2mm] plot ({-\x},{sqrt(6.25-(\x+0.5)^2)-2}) ;
   \draw[domain=0:2, variable=\x, dashed, line width = 0.2mm] plot ({\x},{-sqrt(6.25-(\x-0.5)^2)+2}) ;
   \draw[domain=0:1, variable=\x, dashed, line width = 0.2mm] plot ({-\x},{-sqrt(6.25-(\x+0.5)^2)+2}) ;
   \draw[line width = 0.2mm, dashed, ->] (5,0) .. controls (5.2,0.75) and (6,1.25) .. (7,4.75);
   \draw[line width = 0.2mm, dashed, ->] (5,0) .. controls (5.2,-0.75) and (6,-1.25) .. (7,-4.75);
   \draw[line width =  0.4mm, -] (1.65,1.523) -- (1.65,-1.523) node {};
   \draw[line width =  0.4mm, -] (1.85,0) -- (2.95,0.25) node {};
   \draw[line width =  0.4mm, -] (1.85,0) -- (2.95,-0.25) node {};
    \draw[domain=0.25:3.609, variable=\y, line width = 0.4mm] plot ({3.908-0.707*(\y-1.414)^2},{\y}) ;
   \draw[domain=0.25:3.609, variable=\y, line width = 0.4mm] plot ({3.908-0.707*(\y-1.414)^2},{-\y}) ;
   \draw[domain=0.5:2.8952, variable=\x, line width = 0.4mm] plot ({-\x},{sqrt(15.175-(\x+1)^2)}) ;
   \draw[domain=0.5:2.8952, variable=\x,  line width = 0.4mm] plot ({-\x},{-sqrt(15.175-(\x+1)^2)}) ;
   \draw[line width = 0.4mm, -] (0.5,3.609) .. controls (0.25,3.688) and (-0.25,3.822) .. (-0.5,3.613);
   \draw[line width = 0.4mm, -] (0.5,-3.609) .. controls (0.25,-3.688) and (-0.25,-3.822) .. (-0.5,-3.613);
   \draw[domain=0:1.65, variable=\x, line width = 0.4mm] plot ({\x},{sqrt(5.0425-\x^2)}) ;
   \draw[domain=0:1.65, variable=\x, line width = 0.4mm] plot ({\x},{-sqrt(5.0425-\x^2)}) ;
   \draw[domain=0:2.2455, variable=\x,  line width = 0.4mm] plot ({-\x},{sqrt(5.0425-\x^2)}) ;
   \draw[domain=0:2.2455, variable=\x, line width = 0.4mm] plot ({-\x},{-sqrt(5.0425-\x^2)}) ;
  %
  \node[] at (-0.5,-0.5) (a) {$\mathcal{L}_1$};
  \node[] at (1.2,1.2) (b) {$\mathcal{L}_v$};
  \node[] at (4,2.3) (c) {$\mathcal{L}_w$};
  \node[] at (-1.41,1.41) (d) {$\mathcal{L}_2$};
  \node[] at (1,0.15) (e) {$+$};
  \node[] at (2.2,2.2) (f) {$-$};
  \node[] at (7.25,1.5) (g) {$-$};
  \node[] at (7.25,-1.5) (k) {$-$};
  \node[] at (4.5,1.25) (h) {$+$};
  \node[] at (7,2.5) (h) {$\mathcal{L}_3$};
  \node[] at (7,-2.5) (i) {$\mathcal{L}_3$};
  \draw[-] (-1.5,-0.05) -- (-1.5,0.05) node[below] {$-b$};
   \draw[-] (4,-0.05) -- (4,0.05) node[below] {$\varrho$};
   \node[] at (2,-0.35) (j) {$\alpha^{-1}b$};
\end{tikzpicture}
\caption{Solid lines are used to denote the contour $\mathcal{L}_w$ and $\mathcal{L}_v$ whereas dotted lines indicate the level curves of the equation $\mbox{Re}(\Lambda(z)-\Lambda(b\alpha^{-1}))=0$.}
\label{fig:6}
\end{center}
\end{figure}
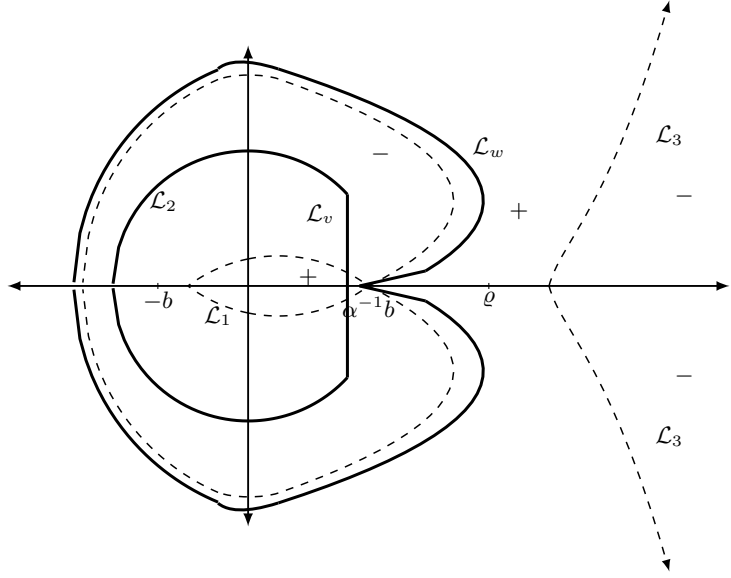

In the next lemma, we show that the contribution of the kernel $K^b_{\zeta_t}$ in the Fredholm determinant $\det (1+K^b_{\zeta_t})_{L^2(\mathcal{L}_w)}$ over the segment $\mathcal{L}^{(2)}_{w}$ decays down exponentially.
\bl\label{eq:CurvedSectionContribution}
For any $w\in \mathcal{L}^{(2)}_w$ and $w^\prime\in \mathcal{L}_w$, we have
\begin{align}\label{eq:KernelBoundOfCursection}
|K^b_{\zeta_t}(w,w^\prime)|\leq C \max\{\epsilon^{-1},\delta^{-1} t^{1/2}\}\exp\left(-tc+s\tilde{\sigma}_{\nu} t^{1/2}|w|\vee |w|^{-1}+\frac{\alpha}{1-b}|w|\right)
\end{align}
for some constants $c,C>0$. Here, $\epsilon$ is same as in Definition~\ref{wvCounyoursBelowCriticality}. To that effect, we have
\begin{align}\label{eq:LimitOfFredDet}
\lim_{t\to \infty}\left|\det(1+K^{b}_{\zeta_t})_{L^2(\mathcal{L}_w)}-\det(1+K^{b}_{\zeta_t})_{L^2(\mathcal{L}^{(1)}_w)}\right| = 0.
\end{align}
\el
\begin{proof}
To begin with, note that for all $w\in \mathcal{L}^{(2)}_w$ and all $v\in \mathcal{L}_v$ there exists a constant  $c>0$ such that $\mbox{Re}\left(\Lambda(w)-\Lambda(v)\right)< -c$. This holds due to the particular way the contours $\mathcal{L}^{(2)}_w,\mathcal{L}_v$ are chosen and the compactness of each of them. Other components of the integrand in \eqref{eq:KernelIntMedForm6} are bounded in the following ways -
\bei
\ii
\begin{align}\label{eq:BoundOnQuadratic}
\left|\exp\left(s\tilde{\sigma}_\nu t^{1/2}\left(\log(v)-\log(w)\right)\right)\right|\leq \exp\left(s\tilde{\sigma}_\nu t^{1/2}(|v|\vee |v|^{-1} +|w|\vee |w|^{-1})\right),
\end{align}

\ii
\begin{align}\label{eq:BoundOnSinSum}
\left|\sum_{k=-\infty}^\infty \frac{1}{\sin\left(\frac{\pi}{\log (b)}\log(v/w)+2\pi ik\right)}\right|\leq \sum_{k=0}^\infty 2\exp(-2\pi k)\exp\left(\frac{\pi}{\log(b)}(|v|\vee |v|^{-1} +|w|\vee |w|^{-1})\right),
\end{align}

\ii
\begin{align}\label{eq:BoundOnLinTermInDen}
\left|\frac{1}{v(v-w)}\right|\leq C\frac{1}{\min_{v\in \mathcal{L}_v,w^\prime\in\mathcal{L}^{(2)}_w}|v-w^\prime|}\leq C\max\{\epsilon^{-1}\vee\delta^{-1} t^{1/2}\}
\end{align}
for some large constant $C$,

\ii
\begin{align}\label{eq:BoundOnRatioOfInfProd}
\left|\frac{(\alpha b^{-1}v;b)_\infty}{\alpha b^{-1}w;b)_\infty}\right|&\leq \prod_{k=1}^\infty (1+b^{k-1}\alpha|v|)\prod_{k=0}^m |(1-b^{m-1}\alpha|w|)|^{-1}\prod_{k=m+1}^\infty (1+2b^{k-1}\alpha|w|)\nonumber\\
&\leq C\exp\left(\frac{\alpha}{1-b}(|v|+2|w|)\right)
\end{align}
where for all $k\geq m$, one has $(1+2b^{k-1}\alpha|w|)(1-b^{k-1}\alpha|w|)\leq 1$.
\ee
 Inequality in \eqref{eq:BoundOnQuadratic} is a simple consequence of the fact that  \begin{align}
|\exp(\log(v)-\log(w))|&\leq \exp\left(\mbox{Re}(v)\vee \mbox{Re}(v^{-1})+\mbox{Re}(w)\vee \mbox{Re}(w^{-1})\right)\\
&\leq \exp\left(|v|\vee |v|^{-1} +|w|\vee |w|^{-1}\right)
\end{align}
for any two complex numbers $v,w$. To see the bound in the right side of \eqref{eq:BoundOnSinSum}, note that $\sin(u+2\pi ik)$ is bounded above by $2\exp(|u|-2\pi k)$ for all $u\in \CC$ and $k\in \ZZ$. It is easy to comprehend the first inequality in \eqref{eq:BoundOnLinTermInDen}. For the second inequality, recall that both the end points of  $\mathcal{L}_w^{(1)}$ are at a distance $\epsilon$ away from the point $b\alpha^{-1}$ for large enough $t$. Thus, it is possible to construct the contour $\mathcal{L}^{(2)}_w$ at  a distance of $\epsilon$ far away from $\mathcal{L}_2$. As we know $\mathcal{L}_v$ is contained in the region enclosed by $\mathcal{L}_2$, thus the minimum value of $|v-w^\prime|$ must be greater than $\epsilon$ for all $v\in \mathcal{L}_v$ and $w^\prime\in \mathcal{L}^{(2)}_w$. On the flip side, if $w^\prime\in \mathcal{L}^{(1)}_w$, then minimum value of $|v-w|$ should be greater than $\delta t^{-1/2}$, thus shows the claim in \eqref{eq:BoundOnLinTermInDen}. First inequality in \eqref{eq:BoundOnRatioOfInfProd} is quite straightforward. To prove the second inequality, note that the finite product $\prod_{k=0}^m |(1-b^{m-1}\alpha|w|)|^{-1}$ is bounded by a large constant. To bound the other components apart from the finite product term, we use a simple inequality $(1+x)\leq e^{x}$ which holds for any $x\in \RR$. These all together imply the inequality in \eqref{eq:KernelBoundOfCursection}.
 Next, we prove the limit in \eqref{eq:LimitOfFredDet}. Let us write down
\begin{align}\label{eq:KernelDecomposition}
K^b_{\zeta_t}(w,w^\prime) = K^b_{\zeta_t}(w,w^\prime)\mathbbm{1}(w\in \mathcal{L}^{(1)}_w)+ K^b_{\zeta_t}(w,w^\prime)\mathbbm{1}(w\in \mathcal{L}^{(2)}_w).
\end{align}
Denote $K^b_{\zeta_t}(w,w^\prime)\mathbbm{1}(w\in \mathcal{L}^{(1)}_w)$ by $K^{b,res}_{\zeta_t}$. Then, using the inequality in \eqref{eq:matrixdetdiff} and the bound on $K^b_{\zeta_t}(w,w^\prime)\mathbbm{1}(w\in \mathcal{L}^{(2)}_w)$ from \eqref{eq:KernelBoundOfCursection}, we have
\begin{align}
\left|\det\left(K^b_{\zeta_t}(w_i,w_j)\right)^n_{i,j=1}-\det\left(K^{b,res}_{\zeta_t}(w_i,w_j)\right)^n_{i,j=1}\right|\leq C 2^n \max\{\epsilon^{-n},(\delta^{-2}t)^{n/2}\}\exp(-ct).
\end{align}
Further, using series expansion of Fredholm determinant, we get the following bound
\begin{align}
\left|\det(1+K^{b}_{\zeta_t})_{L^2(\mathcal{L}_w)}-\det(1+K^{b}_{\zeta_t})_{L^2(\mathcal{L}^{(1)}_w)}\right|&= \left|\det(1+K^{b}_{\zeta_t})_{L^2(\mathcal{L}_w)}-\det(1+K^{b,res}_{\zeta_t})_{L^2(\mathcal{L}_w)}\right|\nonumber\\
& \leq C(1+\exp(-2\epsilon^{-1}))\exp(-ct+2\delta^{-1}t^{1/2})\label{eq:BoundsOnAbsoluteDtDifference}
\end{align}
where the right side of the inequality above goes to $0$ as $t\to \infty$. This completes the proof.
\end{proof}

 Likewise in \eqref{eq:KernelIntMedForm5LocVer}, we further represent $K^b_{\zeta_t}$ as the sum of two different kernels to identify the leading contribution in the computation of the Fredholm determinant $\det(1+K^{b}_{\zeta_t})_{L^2(\mathcal{L}^{(1)}_w)}$. In particular, we set
 \begin{align}\label{eq:KernelIntMedForm7}
 K^b_{\zeta_t} = K^{(1)}_{\zeta_t}+K^{(2)}_{\zeta_t}
 \end{align}
 where for $l=1,2$,
 \begin{align}\label{eq:KernelIntMedForm7LocVer}
  K^{(l)}_{\zeta_t}:= \frac{1}{2i}\sum_{k=-\infty}^{\infty}\int_{\mathcal{L}^{(l)}_{v}}&
  \frac{\exp\left(t(\Lambda(v)-\Lambda(w))\right)}{\sin
  \left(\frac{\pi}{\log (b)}(\log(v)-\log(w))+2\pi ik\right)}\\
  & \times \frac{\exp\left(s\left(\log(v)-\log(w)\right)\right)}{\log(b)(v-w^\prime)} \times \frac{\left(\alpha b^{-1}v;b\right)_\infty}{\left(\alpha b^{-1}w;b\right)_\infty}\times \frac{dv}{v}.
\end{align}
 Following lemma which is in the same spirit of the last Lemma~\ref{eq:BoundsOnAbsoluteDtDifference} proves that one can even ignore the effect of $K^{(2)}_{\zeta_t}$ as $t$ increases to infinity.
 \bl
  For any $w,w^\prime \in \mathcal{L}^{(1)}_w$, we have
  \begin{align}\label{eq:KernelBoundOfLinsection}
  |K^b_{\zeta_t}(w,w^\prime)|\leq C (\epsilon\delta)^{-1} t^{1/2}\exp\left(-tc+s\tilde{\sigma}_{\nu} t^{1/2}|w|\vee |w|^{-1}+\frac{\alpha}{1-b}|w|\right)
  \end{align}
  for some constants $c,C>0$. Thus, we claim
  \begin{align}\label{eq:LimitOfFredDetYetAnother}
\lim_{t\to \infty}\left|\det(1+K^{b}_{\zeta_t})_{L^2(\mathcal{L}^{(1)}_w)}-\det(1+K^{(1)}_{\zeta_t})_{L^2(\mathcal{L}^{(1)}_w)}\right| = 0.
\end{align}
 \el
 \begin{proof}
 Proof of \eqref{eq:KernelBoundOfLinsection} hinges upon the fact that $\mbox{Re}(\Lambda(w)-\Lambda(v))<-c$ for all $v\in \mathcal{L}^{(2)}_v$ and all $w\in \mathcal{L}^{(1)}_w$. This again follows from the compactness of both of the contours, continuity of the map $\mbox{Re}(\Lambda(.))$ and more importantly, strategic position of the contours around $\mathcal{L}_2$. Bounds on the other component of the integrand in \eqref{eq:KernelIntMedForm7LocVer} can be derived as in \eqref{eq:BoundOnQuadratic} - \eqref{eq:BoundOnRatioOfInfProd}. For instance, bounds in \eqref{eq:BoundOnQuadratic} and \eqref{eq:BoundOnSinSum} work perfectly for the current result. In what follows, we recollect bounds on other two terms and sharpen it for the present scenario.
 \bei
 \ii
\begin{align}\label{eq:BoundOnLinTermInDenYetAnother}
\left|\frac{1}{v(v-w)}\right|\leq C\frac{1}{\min_{v\in \mathcal{L}^{(2)}_v,w^\prime\in\mathcal{L}^{(1)}_w}|v-w^\prime|}\leq C\epsilon^{-1}
\end{align}
for some large constant $C$
\ii
\begin{align}\label{eq:BoundOnRatioOfInfProdYetAnother}
\left|\frac{(\alpha b^{-1}v;b)_\infty}{\alpha b^{-1}w;b)_\infty}\right|&\leq \prod_{k=1}^\infty (1+b^{k-1}\alpha|v|)\prod_{k=0}^m |(1-b^{m-1}\alpha|w|)|^{-1}\prod_{k=m+1}^\infty (1+2b^{k-1}\alpha|w|)\nonumber\\
&\leq C\delta^{-1}t^{1/2}\exp\left(\frac{\alpha}{1-b}(|v|+2|w|)\right)
\end{align}
where for all $k\geq m$, one has $(1+2b^{k-1}\alpha|w|)(1-b^{k-1}\alpha|w|)\leq 1$.
 \ee
 One can prove the inequality in \eqref{eq:BoundOnLinTermInDenYetAnother} by recalling that the minimum value of $|v-w^\prime|$ over all $v\in \mathcal{L}^{(2)}_v,w^\prime\in\mathcal{L}^{(1)}_w$ is bounded below by $\epsilon$ (see, Definition~\ref{wvCounyoursBelowCriticality}). In \eqref{eq:BoundOnRatioOfInfProdYetAnother}, we cannot now bound the finite product term by a large constant uniformly over all choices of $t$ and $w$. This is because \[\min_{w\in \mathcal{L}^{(1)}_w} |1-\alpha b^{-1} w|=\delta t^{-1/2}.\]
 This causes the factor $\delta^{-1}t^{1/2}$ to come in the right side of the second inequality. To show \eqref{eq:LimitOfFredDetYetAnother}, we use the inequality in \eqref{eq:matrixdetdiff} to conclude
 \begin{align}
\left|\det\left(K^b_{\zeta_t}(w_i,w_j)\right)^n_{i,j=1}-\det\left(K^{(1)}_{\zeta_t}(w_i,w_j)\right)^n_{i,j=1}\right|\leq C 2^n (\delta\epsilon)^{-n}t^{n/2}\exp(-ct).
\end{align}
whenever $w_1,\ldots ,w_n\in \mathcal{L}^{(1)}_w$. To that effect, Fredholm determinant expansion provides the similar bound as in \eqref{eq:BoundsOnAbsoluteDtDifference} on the difference between $\det(1+K^{b}_{\zeta_t})_{L^2(\mathcal{L}^{(1)}_w)}$ and $\det(1+K^{(1)}_{\zeta_t})_{L^2(\mathcal{L}^{(1)}_w)}$. This ends the proof.
 \end{proof}

 Above two lemmas show that in order to find out the limit of the Fredholm determinant in \eqref{eq:Gaussianlimit}, it suffices to study $K^{(1)}_{\zeta_t}$ over the segment $\mathcal{L}^{(1)}_w$. Thus, like in \eqref{eq:Localization}, we make the following substitutions -
\begin{align}\label{eq:LocalizationBelowCriticality}
v\mapsto b\alpha^{-1}+t^{-1/2}b(\alpha\tilde{\sigma}_\nu)^{-1} \tilde{v}, \quad w\mapsto b\alpha^{-1}+t^{-1/2}b(\alpha\tilde{\sigma}_\nu)^{-1} \tilde{w}, \quad w^\prime\mapsto b\alpha^{-1}+t^{-1/2}b(\alpha\tilde{\sigma}_\nu)^{-1} \tilde{w}^\prime.
\end{align}
  To this end, one can further write $K^{(1)}_{\zeta_t}$ in the following form -

  \begin{align}\label{eq:KernelIntMedForm8LocVer}
  K^{(1)}_{\zeta_t}= \frac{t^{-1/2}b(\alpha \tilde{\sigma}_\nu)^{-1}}{2i}\int_{\mathcal{L}^{(1)}_{\tilde{v}}}&
 \sum_{k=-\infty}^{\infty} \frac{\exp\left(\frac{\tilde{v}^2}{2}-\frac{\tilde{w}^2}{2}+t\left(\tilde{\mathfrak{R}}(t^{-1/2}b(\alpha\tilde{\sigma}_\nu)^{-1}\tilde{w})-\tilde{\mathfrak{R}}(t^{-1/2}b(\alpha\tilde{\sigma}_\nu)^{-1}\tilde{w})\right)\right)}{\sin
  \left(\frac{\pi}{\log (b)}(\log(1+t^{-1/2}(\tilde{\sigma}_\nu)^{-1} \tilde{v})-\log(1+t^{-1/2}(\tilde{\sigma}_\nu)^\prime \tilde{w}))+2\pi ik\right)}\\
  & \times \frac{\exp\left(s\tilde{\sigma}_\nu t^{1/2}\left(\log(1+t^{-1/2}(\tilde{\sigma}_\nu)^{-1} \tilde{v})-\log(1+t^{-1/2}(\tilde{\sigma}_\nu)^{-1} \tilde{w})\right)\right)}{\log(b)(\tilde{v}-\tilde{w}^\prime)}\\
  & \times \frac{\left(\alpha b^{-1}(b\alpha^{-1}+t^{-1/2}b(\alpha\tilde{\sigma}_\nu)^{-1} \tilde{v});b\right)_\infty}{\left(\alpha b^{-1}(b\alpha^{-1}+t^{-1/2}b(\alpha\tilde{\sigma}_\nu)^{-1} \tilde{w});b\right)_\infty}\times \frac{d\tilde{v}}{b\alpha^{-1}+t^{-1/2}b(\alpha\tilde{\sigma}_\nu)^{-1} \tilde{v}}
\end{align}
  where $\mathcal{L}^{(1)}_{\tilde{v}}$ and $\mathcal{L}^{(1)}_{\tilde{w}}$ denote the segments $\mathcal{L}^{(1)}_v$ and $\mathcal{L}^{(1)}_w$ respectively after the substitutions made in \eqref{eq:LocalizationBelowCriticality}. In the next result, we show that in fact the Fredholm determinant of the kernel $K^{(1)}_{\zeta_1}$ converges to the normal distribution.

 \bt\label{GaussianLimitLemma}
 For all $s\in \RR$, we have
 \begin{align}\label{eq:GaussianLimitLemmaStatement}
\lim_{t\to \infty} \det \left(1+K^{(1)}_{\zeta_t}\right)_{L^2(\mathcal{L}^{(1)}_{w})}= G(s).
\end{align}
 \et
  \begin{proof}
  To begin with, we enlist below how each of the components of the integrand in \eqref{eq:KernelIntMedForm8LocVer} is bounded above by some integrable function and what are their limits as $t\to \infty$.
  \bei
  \ii We first look at the inequalities shown below -
  \begin{align}\label{eq:ExpoQuadratic}
  \Big|\exp&\left(\frac{\tilde{v}^2}{2}-\frac{\tilde{w}^2}{2}+t\left(\tilde{\mathfrak{R}}(t^{-1/2}b(\alpha\tilde{\sigma}_\nu)^{-1}\tilde{w})-\tilde{\mathfrak{R}}(t^{-1/2}b(\alpha\tilde{\sigma}_\nu)^{-1}\tilde{w})\right)\right)\Big|\\
  &\leq \exp\left(\mbox{Re}\left(\frac{\tilde{v}^2}{2}-\frac{\tilde{w}^2}{2}\right)+c|\tilde{v}|^2+c|\tilde{w}|^2\right)\leq C^\prime\exp\left(-c^\prime_1|\tilde{v}|^2+c^\prime_2\mbox{Re}(\tilde{w}^2)\right)
\end{align}
for some constants $C^\prime,c^\prime_1>0$ and $c^\prime_2$. To see the first inequality, recall that $|\tilde{\mathcal{R}}(z-b\alpha^{-1})|$ is less than $c|z-b\alpha^{-1}|^2$ in the $\epsilon$-neighborhood of $b\alpha^{-1}$ (see, Definition~\ref{wvCounyoursBelowCriticality}). For showing the second inequality in \eqref{eq:ExpoQuadratic}, one can argue that $\mbox{Re}(\tilde{v}^2)+c|\tilde{v}|^2\leq -c^\prime_1|\tilde{v}|^2$ for some $c^\prime_1>0$. This is because the complex variable $\tilde{v}$ varies over $-2\delta+i\RR$. Thus, $\mbox{Re}(\tilde{v}^2)$ converges to $-\infty$ as $|\tilde{v}^2|\to \infty$. In fact, we get $\mbox{Re}(\tilde{v}^2)/|\tilde{v}|^2\to -1$ whenever $|\tilde{v}|$ goes to $\infty$. This completes the proof of the inequalities above.  Furthermore, as $t\to \infty$, both $t\tilde{\mathfrak{R}}(t^{-1/2}b(\alpha\tilde{\sigma}_\nu)^{-1}\tilde{v})$ and $\tilde{\mathfrak{R}}(t^{-1/2}b(\alpha\tilde{\sigma}_\nu)^{-1}\tilde{w})$ goes to $0$. Thus, we get
\begin{align}\label{eq:ExpoQuadraticLimit}
\exp\left(\frac{\tilde{v}^2}{2}-\frac{\tilde{w}^2}{2}+t\left(\tilde{\mathfrak{R}}(t^{-1/2}b(\alpha\tilde{\sigma}_\nu)^{-1}\tilde{w})-\tilde{\mathfrak{R}}(t^{-1/2}b(\alpha\tilde{\sigma}_\nu)^{-1}\tilde{w})\right)\right)\to \exp\left(\frac{\tilde{v}^2}{2}-\frac{\tilde{w}^2}{2}\right).
\end{align}
\ii Next, we illustrate on the limit of the infinite sum over $k\in \ZZ$. Note that for any $k\in \ZZ$, we have
\begin{align}\label{eq:SinTermBound}
&\left|\sin
  \left(\frac{\pi}{\log (b)}(\log(1+t^{-1/2}(\tilde{\sigma}_\nu)^{-1} \tilde{v})-\log(1+t^{-1/2}(\tilde{\sigma}_\nu)^\prime \tilde{w}))+2\pi ik\right)\right|\\
  &\phantom{======} \geq C(1+t^{-1/2}|\tilde{v}|)\exp(2\pi |k|).
\end{align}
Thus, it shows the series in \eqref{eq:KernelIntMedForm8LocVer} converges absolutely and uniformly over $\tilde{v}$. Hence, as $t$ goes to infinity,
\[\lim_{t\to \infty}\sum_{k=-\infty}^\infty\frac{t^{-1/2}\pi(\log(b))^{-1}}{\sin
  \left(\frac{\pi}{\log (b)}(\log(1+t^{-1/2}(\tilde{\sigma}_\nu)^{-1} \tilde{v})-\log(1+t^{-1/2}(\tilde{\sigma}_\nu)^\prime \tilde{w}))+2\pi ik\right)}=\frac{1}{\tilde{v}-\tilde{w}}.\]
  \ii Likewise in Lemma~\ref{InfiniteProductAtCriticality}, one can prove that for all large $t$, there exists constants $C^\prime$ such that
  \begin{align}\label{eq:InfiniteProductBelowCriticality}
  \left|\frac{\left(\alpha b^{-1}(b\alpha^{-1}+t^{-1/2}b(\alpha\tilde{\sigma}_\nu)^{-1} \tilde{v});b\right)_\infty}{\left(\alpha b^{-1}(b\alpha^{-1}+t^{-1/2}b(\alpha\tilde{\sigma}_\nu)^{-1} \tilde{w});b\right)_\infty}\right|\leq C^\prime \left|\frac{\tilde{v}}{\tilde{w}}\right|.
  \end{align}
  Moreover, we have
  \begin{align}\label{eq:InfiniteProductLimitBelowCriticality}
  \frac{\left(\alpha b^{-1}(b\alpha^{-1}+t^{-1/2}b(\alpha\tilde{\sigma}_\nu)^{-1} \tilde{v});b\right)_\infty}{\left(\alpha b^{-1}(b\alpha^{-1}+t^{-1/2}b(\alpha\tilde{\sigma}_\nu)^{-1} \tilde{w});b\right)_\infty}\stackrel{t\to \infty}{\longrightarrow}\frac{\tilde{v}}{\tilde{w}}.
  \end{align}.
  \ii Lastly, we can put a very crude bound like
  \begin{align}\label{eq:RemainCompBound}
  &\left|\exp\left(s\tilde{\sigma}_\nu t^{1/2}\left(\log(1+t^{-1/2}(\tilde{\sigma}_\nu)^{-1} \tilde{v})-\log(1+t^{-1/2}(\tilde{\sigma}_\nu)^{-1} \tilde{w})\right)\right)\frac{b\alpha^{-1}}{b\alpha^{-1}+t^{-1/2}b(\alpha\tilde{\sigma}_\nu)^{-1} \tilde{v}}\right|\\
  &\phantom{========}\leq C^{\prime\prime}\exp(c(|\tilde{v}|+|\tilde{w}|))\frac{1}{1+|\tilde{v}|}.
  \end{align}
  for some positive constant $C^{\prime\prime},c^{\prime\prime}$.
  Notice that even this crude bound doesn't affect integrability of the product of all bounds together because of the dominating term $\exp(-c^\prime|\tilde{v}|^2)$ in \eqref{eq:ExpoQuadratic}. Further, we have the following limit
\begin{align}\label{eq:LimitOfRemains}
&\exp\left(s\tilde{\sigma}_\nu t^{1/2}\left(\log(1+t^{-1/2}(\tilde{\sigma}_\nu)^{-1} \tilde{v})-\log(1+t^{-1/2}(\tilde{\sigma}_\nu)^{-1} \tilde{w})\right)\right)\frac{b\alpha^{-1}}{b\alpha^{-1}+t^{-1/2}b(\alpha\tilde{\sigma}_\nu)^{-1} \tilde{v}}\\&\phantom{========}\stackrel{t\to \infty}{\longrightarrow}\exp\left(s(\tilde{v}-\tilde{w})\right)
\end{align}
  for the remaining component of the integrand in \eqref{eq:KernelIntMedForm8LocVer}.
  \ee
  To this end, using dominated convergence theorem, one can write
  \begin{align}\label{eq:FinalLimitOfKernelBelowCriticality}
  \lim_{t\to \infty}K^{(1)}_{\zeta_t}(\tilde{w},\tilde{w}^\prime)=\frac{1}{2\pi i}\int_{-2\delta+i\RR}\frac{\exp\left(\frac{\tilde{w}^2}{2}-\frac{\tilde{v}^2}{2}+s(\tilde{w}-\tilde{v})\right)}{(\tilde{v}-\tilde{w}^\prime)(\tilde{v}-\tilde{w})}\cdot \frac{\tilde{v}}{\tilde{w}}\cdot d\tilde{v}
  \end{align}
  for all $\tilde{w},\tilde{w}^\prime$ in the limiting form of the contour $\mathcal{L}^{(1)}_{\tilde{w}}$. In Appendix~\ref{App:ApendixA}, we prove that Fredholm determinant of the limiting kernel is indeed Gaussian distribution function.
  \end{proof}

  \textsc{Proof of Theorem~\ref{Asymptheo2Mainstate}.}
   To show \eqref{eq:HighDenLimit}, note that
   \begin{align}\label{eq:EventsEquivalence}
   \{x_{m_\nu t}\leq \nu t + \frac{ \partial m_\nu}{\partial \nu} sL t^{1/3}\}=\{N_{\nu_s}(t)\geq m_\nu t\}.
   \end{align}
  where $\nu_s:= \nu(1+sL\sigma_\nu \frac{\partial m_\nu}{\partial \nu} t^{-2/3})$ and $L,s\in \RR$. It is easy to see that  $\nu_s$ equals to $\nu$ when $s=0$. Moreover, using Taylor's theorem on $m_{\nu_s}$, one can obtain \[m_{\nu_s}= m_\nu + sL\sigma_\nu\frac{\partial m_\nu}{\partial \nu} L t^{-2/3}+o(t^{-2/3}).\]
 Let us choose $L= \left(\frac{\partial m_\nu}{\partial \nu}\right)^{-1}$. To that effect, we have
  \begin{align}\label{eq:ProbabilityEquivalence}
\mathbb{P}\left( N_{\nu_s}(t)\geq m_\nu t \right)& =  \mathbb{P}\left(m_{\nu_s} t -N_{\nu_s}t\leq s\sigma_\nu  t^{1/3}+o(t^{1/3})\right)
 \end{align}
  Theorem~\ref{eq:HighDenLimit} implies that the right side in \eqref{eq:ProbabilityEquivalence} converges to $F_{GUE}(s)$. Thus, it completes the proof of \eqref{eq:HighDenLimit2}. Similarly, the proofs of \eqref{eq:CriticalDenLimit2} and \eqref{eq:LowDenLimit2} follow from \eqref{eq:CriticalDenLimit} and \eqref{eq:LowDenLimit} respectively.
  \qed

\textsc{Proof of Corollary~\ref{LimitShapesCorollary}}

   Theorem~\ref{Asymptheo1Mainstate} implies $(t^{-1}N_{\nu t}(t)-m_\nu)=O_p(t^{-2/3})$ when $\nu >\nu_0$. Thus, it proves the first case. Further, note that for any $\epsilon>0$, there exists $\delta_\epsilon>0$ such that whenever $\nu\leq \nu_0+\delta_\epsilon$, limit shape $m_\nu <\epsilon$. To that effect, for any $\nu\leq \nu_0$, we have
  \begin{align}\label{eq:TrivialUpperBound}
  \limsup_{t\to \infty}\frac{N_{\nu t}(t)}{t}\leq \lim_{t\to \infty}\frac{N_{(\nu_0+\delta_\epsilon)t}(t)}{t}\leq \epsilon
  \end{align}
  with high probability. Henceforth, this shows the second claim.

\section{Weak Scaling Limit}\label{WeakScaleLim}

Our goal here is to prove the Theorem \ref{KPZLimitTheo}. Following two propositions provide the necessary
arguments needed to complete the proof.
\bp\label{Tightness}
The collection of processes $\{Z_\epsilon\}_\epsilon$ is tight in $C(\RR_+\times \RR)$.
\ep

\bp\label{SubseqArguements}
Any limiting point $\mathcal{Z}$ of $\{Z_\epsilon\}_\epsilon$ solves the SHE.
\ep

 To prove Proposition~\ref{Tightness}-\ref{SubseqArguements}, we first show that the G\"artner transform (recall Definition~\ref{GartnerTransform}) of HL-PushTASEP satisfies a discrete SHE. This step is in harmony with
 the other works like \cite[Section~3.2]{Bertini1995}, \cite[Section~2]{Dembo2016}, \cite[Section~3]{CL15}. Let us note that here onwards, constants will change from line to line. Unless specified, they do not depend on anything else.

   \subsection{Discrete SHE}\label{DiscreteSHEconstruction}

   First, we derive a dynamical equation for $Z(t,x)$. Recall that any particle which is at left side of any site $x$ when gets excited can move across $x$ with certain probability. For instance, once the clock of the particle at site $y(<x)$ rings at time $t$, it jumps across $x$ with probability $b^{x-y-N_x(t)+N_y(t)}$. For notational convenience, let us introduce the following notation 
\begin{align}\label{eq:CharacterisLineOfKPZ}
\kappa(x,t) = x+\frac{t}{1-b^{\nu_\epsilon}}.
\end{align} From the definition of $Z(t,x)$ in \eqref{eq:HopfColeTransform}, it is clear that contributions of the all the right hops across $x$ in $dZ(t,x)$ during the interval $[t,t+dt]$ is given by
$(b^{-1}-1)Z(t,x)d\mathcal{P}^x_t$
where $\{\mathcal{P}^x_t\}_{-\infty <x<\infty}$ denotes a sequence of Poisson point processes with rates
\[\sum_{l=1}^{ \kappa(x,t)}\eta_t(l)b^{\kappa(x,t)-l-N_{\kappa(x,t)}(t)+N_l(t)}\]
and $\eta_t(l)$ is the occupation variable for the position $l$ (see, Definition~\ref{eq:ParticleCountDefinition}).  
There is also a contribution coming from the exponential growth term $\exp(-\mu_\epsilon t)$.
Collecting all the terms together, we can write down the following
\begin{align}\label{eq:SHEDecomposition}
dZ(t,x)&=\left((b^{-1}-1)\sum_{l=1}^{\kappa(x,t)}\eta_t(l)b^{\kappa(x,t)-l-N_{\kappa(x,t)}(t)+N_l(t)}-\frac{(1-\nu_\epsilon)\log b}{(1-b^{\nu_\epsilon})}-\mu_\epsilon\right)Z(t,x)dt \nonumber\\
&+(b^{-1}-1)Z(t,x)\left(d\mathcal{P}^x_t-\sum^{\kappa(x,t)}_{l=1}\eta_t(l)b^{\kappa(x,t)-l-N_{\kappa(x,t)}(t)+N_l(t)}dt\right).
\end{align}
 For the sake of simplicity, denote
\begin{align}\label{eq:DefMartingale}
dM(t,x)= (b^{-1}-1)\left(d\mathcal{P}^x_t-\sum_{l=1}\eta_t(l)b^{\kappa(x,t)-l-N_{\kappa(x,t)}(t)+N_l(t)}dt\right).
\end{align}
Further, let us define
\begin{align}\label{eq:OmegaDef}
\Omega_x:= (b^{-1}-1)\sum_{l=1}\eta_t(l)b^{\kappa(x,t)-l-N_{\kappa(x,t)}(t)+N_l(t)}.
\end{align}
\bl\label{RepresentationLemma}
We have
\begin{align}\label{eq:HeatKernelLaplacian}
\Omega_x Z(t,x)= (b^{-1}-1)\sum_{l< x} b^{\nu_\epsilon\left(x-l\right)}Z(t,l)-Z(t,x).
\end{align}
\el
\begin{proof}
Using the form of $\Omega_x$ given in \eqref{eq:OmegaDef}, we get 
\[\Omega_xZ(t,x)=(b^{-1}-1)\sum_{l\leq \kappa(x,t)}\eta_t(l) b^{\nu_\epsilon\left(\kappa(x,t)-l\right)}b^{N_l(t)-(1-\nu_\epsilon)l}\exp(-\mu_\epsilon t).\]
We will show
\begin{align}\label{eq:InterludeEst}
\sum_{l\leq  \kappa(x,t)} b^{\nu_\epsilon\left(\kappa(x,t)-l\right)}b^{N_l(t)-(1-\nu_\epsilon)l}&+\exp(\mu_\epsilon t)\Omega_xZ(x,t)=b^{-1}\sum_{l< \kappa(x,t)} b^{\nu_\epsilon\left(\kappa(x,t)-l\right)}b^{N_l(t)-l/2}.
\end{align}
One can easily see that above result implies \eqref{eq:HeatKernelLaplacian}.
For proving \eqref{eq:InterludeEst}, we have
$\Omega_x Z(t,x)=\lambda_\epsilon\sqrt{\epsilon}(1+C\sqrt{\epsilon})\sum_{l\leq \kappa(x,t)}\eta_t(l) b^{\nu_\epsilon\left(\kappa(x,t)-l\right)}b^{N_l(t)-(1-\nu_\epsilon)l}$ by expanding $b^{-1}-1$ in terms of $\epsilon$. Note that the value of $C$ is $\sum_{n=1}^\infty \frac{\lambda^{n+1}_\epsilon\epsilon^{n}}{(n+1)!}$. To this effect, we have
\begin{align}
\sum_{l\leq  \kappa(x,t)}& b^{\nu_\epsilon\left(\kappa(x,t)-l\right)}b^{N_l(t)-(1-\nu_\epsilon)l}+\Omega_xZ(x,t)\exp(\mu_bt)\\
&=\sum_{l\leq  \kappa(x,t)} (1+\sqrt{\epsilon}\lambda_\epsilon\eta_t(l))b^{\nu_\epsilon\left(\kappa(x,t)-l\right)}b^{N_l(t)-(1-\nu_\epsilon)l}+\frac{C\lambda_\epsilon\epsilon}{1+C\sqrt{\epsilon}}\Omega_xZ(t,x)\\
&= \sum_{l\leq \kappa(x,t)} [\exp(\lambda_\epsilon\sqrt{\epsilon}\eta_t(l))-C\lambda_\epsilon\epsilon\eta_t(l)]b^{\nu_\epsilon\left(\kappa(x,t)-l\right)}b^{N_l(t)-(1-\nu_\epsilon)l}\\&\phantom{=====}+\frac{C\lambda_\epsilon\epsilon}{1+C\sqrt{\epsilon}}\Omega_xZ(t,x)\\
&= b^{-1}\sum_{l< \kappa(x,t)} b^{\nu_\epsilon\left(\kappa(x,t)-l\right)}b^{N_l(t)-(1-\nu_\epsilon)l}.
\end{align}
We wrote down the equality in the third line by identifying $(1+\eta_t(l)\lambda_\epsilon\sqrt{\epsilon})$ as $\exp(\eta_t(l)\lambda_\epsilon\sqrt{\epsilon})-C\lambda_\epsilon\epsilon\eta_t(l)$. In the last line, we used the fact $N_x(t)-\eta_t(x)=N_{x-1}(t)$ for any $x\in\ZZ_{\geq 0}$. This completes the proof.
\end{proof}

Using the definition of  $\mu_\epsilon$ in \eqref{eq:MuEpsilon}, one can further write down
\begin{align}\label{eq:SimpleHopfCole}
dZ(t,x) = \frac{1-b}{b(1-b^{\nu_\epsilon})}\left(\sum_{k=0}^{\infty}(1-b^{\nu_\epsilon}) b^{\nu_\epsilon
k}Z(t,x-k)-Z(t,x)\right)dt+Z(t,x)dM(t,x).
\end{align}

 Let $\{R^{}_i\}_{i\geq 1}$ be  iid random variables such that

\begin{align}\label{eq:RSupportMassFunction}
\mbox{supp}(R_i)=\ZZ_{\geq 0} - \frac{1}{1-b^{\nu_\epsilon}},\quad \text{and }\mathbb{P}\left(R_i=k-\frac{1}{1-b^{\nu_\epsilon}}\right)=(1-b^{\nu_\epsilon})b^{\nu_\epsilon k} \quad \text{for }k\in \ZZ_{\geq 0}.
\end{align}
Note that $R_i$ satisfies $\mathbb{E}(R_i)=0$.  Under weak noise scaling, one can note that for small $\epsilon $, $R_i+(\nu_\epsilon\lambda_\epsilon\sqrt{\epsilon})^{-1}$ follows $\mbox{Geo}(\lambda_\epsilon\nu_\epsilon\sqrt{\epsilon})$. Let us denote a Poisson point process of rate $\frac{b^{-1}-1}{1-b^{\nu_\epsilon}}$ associated with $R_i$'s by $\Gamma_t$. Its increments over the interval $[t_1,t_2]$ will be denoted as $\Gamma_{t_1,t_2}$. Further, denote 
 $\Xi(t_1,t_2):=\ZZ_{\geq 0}- \phi_\epsilon $
for $t_1\leq t_2$ where 
\begin{align}\label{eq:SateSpaceForRandomWalkMean}
\phi_\epsilon:=\frac{(b^{-1}-1)(t_2-t_1)}{(1-b^{\nu_\epsilon})^2}.
\end{align}
 In this context, define the transition probability
\begin{align}\label{eq:TransitionProbabilityForPoissonProcess}
p_\epsilon(t_1,t_2,\xi):=\mathbb{P}\left(\sum_{i=1}^{\Gamma_{t_1,t_2}}R_i+\frac{\Gamma_{t_1,t_2}}{1-b^{\nu_\epsilon}}- \phi_\epsilon= \xi\right)
\end{align}
for $\xi\in \Xi(t_1,t_2)$. Note that support space of the transition probability $p_\epsilon(t_1,t_2,.)$ is $\Xi(t_1,t_2)$. Moreover, one can see that $\{p_\epsilon(t_1,t_2,\xi)\}_{t_1\leq t_2,\xi\in \mathfrak{t_1,t_2}}$ form a semi-group. We use the notation $[p_\epsilon(t_1,t_2)*f(t_1)]$ to denote the convolution $\sum_{\zeta\in \Xi(t_1,t_2)}p_\epsilon(t_1,t_2,\xi-\zeta)f(t_1,\zeta)$ of any function $f$ with the semi-group. One can further note that $p_\epsilon(t,.)$ satisfy the following differential equation,
\begin{align}\label{eq:HeatKernelSDE}
\frac{d}{dt}p_\epsilon(t,x)&=\frac{(1-b)}{b(1-b^{\nu_\epsilon})}\left(\sum_{k=0}^{\infty}(1-b^{\nu_\epsilon})b^{\nu_\epsilon k}p_\epsilon(t,x-k)-p_\epsilon(t,x)\right)\\
p_\epsilon(0,x)&=\delta_0(x).
\end{align}
To put this into words, $p_\epsilon$ behaves as a discrete heat kernel. To that effect, we get the following simple form of a discrete SHE
\begin{align}\label{eq:DiscreteSHE}
Z(t,x)&=[p_\epsilon(t_1,t)*Z(t_1)](x)+\int_{t_1}^t [p_\epsilon(s,t)*Z(s)dM(s)](x)\quad \text{where }0\leq t_1 <t.
\end{align}
Further, utilizing the particle dynamics of HL-PushTASEP, we derive the time derivative of the quadratic variation of the martingale $M(t,x)$ defined in \eqref{eq:DefMartingale}.

\bl
We have
\begin{align}\label{eq:TotVarDiff}
(b^{-1}-&1)^{-2}Z(t,x_1)  Z(t,x_2)d\langle M(t,x_1),M(t,x_2)\rangle \\&= \frac{b^{\nu_\epsilon|x_1-x_2|}}{1-b^{\nu_\epsilon}}Z(t,x_1\wedge x_2)\left([p_\epsilon(t,t+dt)*Z(t)](x_1\wedge x_2)-\frac{(1-b^{\nu_\epsilon})b^{-1}}{b^{-1}-1}Z(t,x_1\wedge x_2)dt\right)
\end{align}
\el
\begin{proof}
We know that $d\langle M(t,x_1), M(t,x_2)\rangle$ equals to the rate at which both the Poisson process associated with $M(t,x_1)$ and $M(t,x_2)$ acknowledges a jump in the interval $[t,t+dt]$. Thus, it must be equals to the rate at which any particle to the left of $\kappa(x_1\wedge x_2, t)$ hops across $\kappa(x_1,t)$ and $\kappa(x_2, t)$ at time $t$. Here, we denote $\min\{x_1,x_2\}$ and $\max\{x_1,x_2\}$ by $x_1\wedge x_2$ and $x_1\vee x_2$ respectively. To this effect, one can write
\[d\langle M(t,x_1),M(t,x_2)\rangle = (b^{-1}-1)^2\sum_{l\leq \kappa(x_1\wedge x_2,t)}\eta_t(l)b^{\kappa(x_1\vee x_2,t)-l-N_{\kappa(x_1\vee x_2,t)}(t)+N_{l}(t)}dt.\]
Notice that for $l\leq \kappa(x_1\wedge x_2,t)$, we have \[b^{\kappa(x_1\vee x_2,t)-l-N_{\kappa(x_1\vee x_2,t)}(t)+N_{l}(t)}=Z(t,x_1\wedge x_2)Z^{-1}(x_1\vee x_2)b^{\nu_\epsilon|x_1-x_2|}b^{\kappa(x_1\wedge x_2,t)-l-N_{\kappa(x_1\wedge x_2,t)}(t)+N_{l}(t)}.\]
Therefore, we can say
\begin{align}\label{eq:TotVarRep}
(b^{-1}-1)^{-2}&Z(t,x_1)Z(t,x_2) d\langle M(t,x_1),M(t,x_2)\rangle \\& =Z^2(t,x_1\wedge x_2)b^{\nu_\epsilon|x_1-x_2|} \sum_{l\leq \kappa(x_1\wedge x_2,t)} \eta_t(l) b^{\kappa(x_1,x_2,l)l-N_{\kappa(x_1\wedge x_2,t)}(t)+N_l(t)}dt.
\end{align}
One can recognize the following
\[\sum_{l\leq \kappa(x_1\wedge x_2,t)} \eta_t(l) b^{\kappa(x_1,x_2,l)l-N_{\kappa(x_1\wedge x_2,t)}(t)+N_l(t)}=(b^{-1}-1)^{-1}\Omega_{x_1\wedge x_2}Z(t,x).\]
Thus, using Lemma~\ref{RepresentationLemma}, we can further say that $(b^{-1}-1)^{-2}Z(t,x_1)Z(t,x_2) d\langle M(t,x_1),M(t,x_2)\rangle$ equals to
\[ (1-b^{\nu_\epsilon})^{-1}b^{\nu_\epsilon|x_1-x_2|}Z(t,x_1\wedge x_2)\left(\sum_{k=1}^\infty (1-b^{\nu_\epsilon})b^{\nu_\epsilon k}Z(t,x_1\wedge x_2-k)-\frac{1-b^{\nu_\epsilon}}{b^{-1}-1}Z(t,x_1\wedge x_2)\right) dt.\]
Note that one can simplifies the term inside the bracket in the above expression to
\begin{align}
 \left([p_\epsilon(t,t+dt)*Z(t)](x_1\wedge x_2)-\frac{(1-b^{\nu_\epsilon})b^{-1}}{b^{-1}-1}Z(t,x_1\wedge x_2)dt\right).
\end{align}
This shows the claim.
\end{proof}

\subsection{Heat Kernel Estimates \& Tightness}
Now, we aim at proving the tightness of the sequence $\{Z_\epsilon\}_\epsilon$. For that, we need to specialize on certain properties the discrete heat kernel $p_\epsilon$ (see, \eqref{eq:HeatKernelSDE} for definition). In the following proposition, we provides the necessary estimates for the heat kernel.
\bp\label{KernelEstimate}
Given any $T>0$, $u\in \RR$ and $v\in (0,1]$, there exists $C$ depending on $T$ and $u$ such that
\bei
\ii \begin{align}\label{eq:KernelEstimate1} \sum_{\zeta\in \Xi(t_1,t_2)}p_\epsilon(t_1,t_2,\zeta)\exp(u\epsilon |\zeta|)\leq C,
\end{align}
\ii \begin{align}\label{eq:KernelEstimate2}
\sum_{\zeta\in \Xi(t_1,t_2)}p(t_1,t_2,\zeta)|\zeta|^v\exp(u\epsilon|\zeta|)\leq C\epsilon^{-v/2}(|t_2-t_1|)^{v/2},
\end{align}
\ii \begin{align}\label{eq:KernelEstimate3}
p_\epsilon(t_1,t_2,v)\leq C\sqrt{\epsilon}\min\{1,(|t_2-t_1|)^{-1/2}\},
\end{align}
\ii \begin{align}\label{eq:KernelEstimate4}
|p_\epsilon(t_1,t_2,\zeta)-p_\epsilon(t_1,t_2,\zeta^\prime)|\leq C\epsilon^{(1+v)/2}|\zeta-\zeta^\prime|^{v}\min\{1,(t_2-t_1)^{-(1+v)/2}\},
\end{align}
for all $t_1\leq t_2\in (0,\epsilon^{-1}T]$ and $\zeta,\zeta^\prime \in \Xi(t_1,t_2)$.
\ee
\ep

\begin{proof}
\bei
\ii
 Recall that $\phi_\epsilon$ (see, \eqref{eq:SateSpaceForRandomWalkMean}) is the mean of a random walk at time $t$ of iid Geometric random variable with parameter $(1-b^{\nu_\epsilon})$ driven by a Poisson point process of rate $(1-b^{\nu_\epsilon})^{-1}(b^{-1}-1)$ over the interval $[t_1,t_2]$. We have noted that the support of the transition probability $p_\epsilon(t_1,t_2,.)$ is given by $\ZZ_{\geq 0}-\phi_\epsilon$. Also define $u^\prime:=\epsilon^{1/2}(\nu_\epsilon\lambda_\epsilon)^{-1} u$. For any $\zeta\in \Xi(t_1,t_2)$ and $n\in \NN$, let us denote 
\begin{align}\label{eq:SpecialBinomialNotation}
\mathfrak{B}(n,\zeta):= \binom{n+\zeta+\phi_\epsilon-1}{\zeta+\phi_\epsilon-1}.
\end{align} 
For small $\epsilon$, $(1-b^{\nu_\epsilon})^{-1}(b^{-1}-1)\leq \nu^{-1}_\epsilon$. Hence, Poisson process associated with the semi-group $\{p_\epsilon(t_1,t_2,.)\}$ will be dominated by a Poisson process of rate $\nu^{-1}_\epsilon$. For notational simplicity, we introduce the following 
\begin{align}\label{eq:PoissonProbMass}
m(\nu_\epsilon, t_2-t_1,n)=\exp\left(-\nu^{-1}_\epsilon(t_2-t_1)\right)\frac{\nu^{-1}_n(t_2-t_1)^n}{n}.
\end{align}
 To that effect, we can bound the sum $\sum_{\zeta\in \Xi(t_1,t_2)} p_\epsilon(t_1,t_2,\zeta)\exp(u\epsilon|\zeta|)$ by
\begin{align}
\sum_{n=0}^\infty\sum_{\zeta\in \Xi(t_1,t_2)}m(\nu_\epsilon, t_2-t_1,n) &\mathfrak{B}(n,\zeta)(1-b^{\nu_\epsilon})^n b^{\nu_\epsilon(\zeta+\phi_\epsilon)}\exp(u\epsilon|\zeta|)
\end{align}
Using a simple fact $e^{u\epsilon|\zeta|}\leq e^{u\epsilon \zeta}+ e^{-u\epsilon\zeta}$ and taking the sum over all $\zeta \in \Xi(t_1,t_2)$, one can further bound the sum above by $(A)+(B)$ where 
\begin{align}
(A):&= \sum_{n=0}^\infty m(\nu_\epsilon, t_2-t_1,n)\left(\frac{(1-b^{\nu_\epsilon})}{(1-b^{\nu^{-1}_\epsilon(1-u^\prime)})}\right)^n\exp(-u\epsilon\phi_\epsilon)\\
&= \exp\left(\nu^{-1}_\epsilon\frac{1-b^{\nu_\epsilon}}{1-b^{\nu_\epsilon(1-u^\prime)}}(t_2-t_1)-\nu^{-1}_\epsilon(t_2-t_1)-u\epsilon\phi_\epsilon\right),\label{eq:TypeI}\\
(B):&= \sum_{n=0}^\infty m(\nu_\epsilon, t_2-t_1,n) \left(\frac{(1-b^{\nu_\epsilon})}{(1-b^{\nu_\epsilon(1+u^\prime)})}\right)^n\exp(u\epsilon\phi_\epsilon)\\
&=\exp\left(\nu^{-1}_\epsilon\frac{1-b^{\nu_\epsilon}}{1-b^{\nu_\epsilon(1+u^\prime)}}(t_2-t_1)-\nu^{-1}_\epsilon(t_2-t_1)+u\epsilon\phi_\epsilon\right).\label{eq:TypeII}
\end{align}

 Now recall $b=e^{-\epsilon^{1/2}}$. Henceforth, one can say
\[\frac{1-b^{\nu_\epsilon}}{1-b^{\nu_\epsilon(1-u^\prime)}}=1+u\nu^{-1}_\epsilon \epsilon^{1/2}+O(\epsilon)\quad\text{and }\frac{1-b^{\nu_\epsilon}}{1-b^{\nu_\epsilon(1+u^\prime)}}=1-u\nu^{-1}_\epsilon\epsilon^{1/2} +O(\epsilon).\]
Also, we have \[u\epsilon\phi_\epsilon = u\epsilon(1-b^{\nu_\epsilon})^{-2}(b^{-1}-1)(t_2-t_1)= u(\nu^{-2}_\epsilon\epsilon^{1/2}+O(\epsilon))(t_2-t_1).\]
This enforces
$\nu^{-1}_\epsilon\frac{1-b^{1/2}}{1-b^{(1-u^\prime)/2}}(t_2-t_1)-\nu^{-1}_\epsilon(t_2-t_1)-u\epsilon\phi_\epsilon$ and $\nu^{-1}_\epsilon\frac{1-b^{1/2}}{1-b^{(1+u^\prime)/2}}(t_2-t_1)-\nu^{-1}_\epsilon(t_2-t_1)+u\epsilon\phi_\epsilon$ to be bounded by some constant whenever $t_1,t_2\in (0,\epsilon^{-1}T]$.
Thus, one can note right side in the last line of \eqref{eq:TypeI} and \eqref{eq:TypeII} are bounded by some constant $C$ which depends only on $u$ and $T$.
\ii
 We can bound $\sum_{\zeta\in \mathcal{E}(t_1,t_2)}p_\epsilon(t_1,t_2,\zeta)|\zeta|^{v}\exp(u\epsilon|\zeta|) $ by 
\begin{align}
 \sum_{n=0}^\infty m(\nu_\epsilon, t_2-t_1,n)&\sum_{\mathclap{\zeta\in \Xi(t_1,t_2)}}\mathfrak{B}(n,\zeta)(1-b^{\nu_\epsilon})^n b^{\nu_\epsilon(\zeta+\phi_\epsilon)}|\zeta|^v\exp(u\epsilon|\zeta|)\leq (A^\prime)+(B^\prime)
\end{align}
where 
\begin{flalign}
(A^\prime)&:= \sum_{n=0}^\infty m(\nu_\epsilon, t_2-t_1,n)\sum_{\mathclap{\zeta\in \Xi(t_1,t_2)}}\mathfrak{B}(n,\zeta)(1-b^{\nu_\epsilon})^n b^{\nu_\epsilon(\zeta+\phi_\epsilon)}|\zeta|^v\exp(u\epsilon\zeta)&\nonumber\\
&\leq \sum_{n=0}^\infty m(\nu_\epsilon, t_2-t_1,n)\left|\frac{n}{1-b^{\nu_\epsilon(1-u^\prime)}}-\phi_\epsilon\right|^v\left(\frac{1-b^{\nu_\epsilon}}{1-b^{\nu_\epsilon(1-u^\prime)}}\right)^n\exp(-u\epsilon\phi_\epsilon),&\label{eq:TypeISec}\\
(B^\prime)&:=\sum_{n=0}^\infty m(\nu_\epsilon, t_2-t_1,n)\sum_{\mathclap{\zeta\in \Xi(t_1,t_2)}}\mathfrak{B}(n,\zeta)(1-b^{1/2})^n b^{\nu_\epsilon(\zeta+\phi_\epsilon)}|\zeta|^v\exp(-u\epsilon\zeta)&\\
&\leq \sum_{n=0}^\infty  m(\nu_\epsilon, t_2-t_1,n)\left|\frac{n}{1-b^{\nu_\epsilon(1+u^\prime)}}-\phi_\epsilon\right|^v\left(\frac{1-b^{\nu_\epsilon}}{1-b^{\nu_\epsilon(1+u^\prime)}}\right)^n\exp(u\epsilon \phi_\epsilon).&\label{eq:TypeIISec}
\end{flalign}
Notice $g(x)=|x|^v$ is a concave function for $v\in (0,1)$. We used Jensen's inequality for concave function in \eqref{eq:TypeISec} - \eqref{eq:TypeIISec} of the above computation. To this end, it can be noted that
\[\frac{n}{1-b^{\nu_\epsilon(1\pm u^\prime)}}=\frac{(1-b^{\nu_\epsilon})n}{(1-b^{\nu_\epsilon(1\pm u^\prime)})^2} +(\pm C+O(\epsilon^{1/2}))n\]
for some constant $C$. Also, $\phi_\epsilon=\frac{\nu^{-1}_\epsilon(t_2-t_1)(1-b^{\nu_\epsilon})}{(1-b^{\nu_\epsilon(1\pm u^\prime)})^2}\pm C^{\prime}(t_1-t_1)+O(\epsilon)$ for some constant $C^\prime$. For any $v<1$, using H\"older's inequality, we can write
\[\mathbb{E}\left(|aX+b|^v\right)\leq (\mathbb{E}((ax+b)^2)^{v/2}\]
for any random variable $X$.
 Thus, we have
\begin{align}
\sum_{n=0}^\infty  m(\nu_\epsilon, t_2-t_1,n)&\left|\frac{n}{1-b^{\nu_\epsilon(1-u^\prime)}}-\phi_\epsilon\right|^v\left(\frac{1-b^{\nu_\epsilon}}{1-b^{\nu_\epsilon(1\pm u^\prime)}}\right)^n\exp(\pm u\epsilon\phi_\epsilon)\\&\leq \exp\left(\nu^{-1}_\epsilon\frac{1-b^{\nu_\epsilon}}{1-b^{\nu_\epsilon(1\pm u^\prime)}}(t_2-t_1)-\nu^{-1}_\epsilon(t_2-t_1)\pm u\epsilon\phi_\epsilon\right)\\&\times \left(\frac{\nu^{-1}_\epsilon(t_2-t_1)(1-b^{\nu_\epsilon})^2}{(1-b^{\nu_\epsilon(1\pm u^\prime)})^4}+C_1(t_2-t_1)^2\right)^{v/2}.\label{eq:JensenTrick}
\end{align}
To this end, we have
\[\frac{\nu^{-1}_\epsilon(t_2-t_1)(1-b^{\nu_\epsilon})^2}{(1-b^{\nu_\epsilon(1\pm u^\prime)})^4}+C_1(t_2-t_1)^2\lesssim \epsilon^{-1}(t_2-t_2).\]
Also, one can note that that other part of the product on the right side of \eqref{eq:JensenTrick} appeared in \eqref{eq:TypeI} and \eqref{eq:TypeII}. Thus, it is bounded by some constant which only depends on $T$. This completes the proof.
\ii Using the inversion formulae of the characteristic function, we can write down the following
\[p_\epsilon(t_1,t_2,\zeta)\leq \sum_{n=0}^\infty \exp(-\nu^{-1}_\epsilon(t_2-t_1))\frac{\nu^{-n}_\epsilon(t_2-t_1)^n}{n!}\frac{1}{2\pi i}\int_{-\pi}^\pi  |\phi(s,\epsilon)|^n ds\]
where
\begin{align}\label{eq:CharFuncDef}
\phi(s,\epsilon)=\sum_{\xi=0}^\infty \exp\left(is\left(\xi-(1-b^{\nu_\epsilon})^{-1}\right)\right)(1-b^{\nu_\epsilon})b^{\nu_\epsilon\xi}.
\end{align}
Now, using the estimate $|\phi(s,\epsilon)|\leq \frac{\sqrt{\epsilon}}{\sqrt{\epsilon+s^2}}$ (obtained by Taylor series expansion), we get
\begin{align}\label{eq:CharFunEst}
\int_{-\pi}^{\pi}|\phi(s,\epsilon)|^n ds\leq  \int_{[-\pi,\pi]} \left(\frac{\sqrt{2}}{1+|s|/\sqrt{\epsilon}}\right)^nds\leq 2\sqrt{\epsilon} \int_0^\infty \left(\frac{\sqrt{2}}{1+s}\right)^nds.
\end{align}
Notice that when $n=1$, transition probability is just $(1-b)b^{\zeta+\phi_\epsilon}$ which is anyway less that $\epsilon$. But for $n>1$, using \eqref{eq:CharFunEst} one can further say
\[\int_{-\pi}^{\pi}|\phi(s,\epsilon)|^n ds\leq 2^{1+n/2}\sqrt{\epsilon} \frac{1}{n-1}\]
 which implies
\begin{align}
p_\epsilon(t_1,t_2,\zeta)\leq C \sqrt{\epsilon} \exp(-\nu^{-1}_\epsilon(t_2-t_1))\left[\nu^{-1}_\epsilon(t_2-t_1)+2\sum_{n=2}^\infty\frac{(\sqrt{2}\nu^{-1}_\epsilon)^{n}(t_2-t_1)^n}{n!(n-1)}\right].\\\label{eq:PossonBound}
\end{align}
 If $X\sim \mbox{Poisson}(\tau)$, then using Bernstein's inequality, we have
 \[\mathbb{P}\left(|X-\tau|\geq t\sqrt{\tau}\right)\leq 2\exp(-Ct^2)\]
 for some constant. This being said, one can deduce that $\mathbb{E}\left(\frac{1}{X}\right)\leq C\min \{1,\tau^{-1/2}\}$. One can now see second term in \eqref{eq:PossonBound} has an upper bound $\min \{1,(t_2-t_1)^{-1/2}\}$. Also, $(t_2-t_1)\exp(-2(t_2-t_1))$ can be bounded by $(t_2-t_1)^{-1/2}$ upto some multiplicative constant when $\epsilon$ is small enough. Thus, we get
 \[p_\epsilon(t_1,t_2,\zeta)\leq C \sqrt{\epsilon}\min \{1,(t_2-t_1)^{-1/2}\}.\]
\ii Using uniform $v$-Holder continuity of the map $x\mapsto e^{ix}$ for $x\in \RR$ and inversion formulae of the characteristic functions, one can write down
\begin{align}\label{eq:TPDiffCrudeBd}
|p_\epsilon(t_1,t_2,\zeta)-p_\epsilon(t_1,t_2,\zeta^\prime)|\leq \sum_{n=0}^\infty \exp(-(t_2-t_1))\frac{(t_2-t_1)^n}{n!}\frac{1}{2\pi i}\int_{-\pi}^\pi |s(\zeta-\zeta^\prime)|^v |\phi(s,\epsilon)|^n ds.
\end{align}
where $\phi(.,.)$ is defined in \eqref{eq:CharFuncDef}. Using the estimate of the characteristic function in the last part, we can have the following bound
\begin{align}\label{eq:TPDiffCrudeBd1}
\int_{-\pi}^{\pi}s^v|\phi(s,\epsilon)|^n ds\leq 2\epsilon^{(1+v)/2}\int^\infty_0 s^v\left(\frac{\sqrt{2}}{1+s}\right)^nds
  \end{align}
  when $n>2$. Note that $(n-1)(1+s)^{-n}$ is a density on $\RR_{+}$. It is easy to see that expectation under that density is $1/(n-2)$ whenever $n>2$. Thus, using Jensen's inequality, one can conclude that right side of \eqref{eq:TPDiffCrudeBd1} is bounded above by $2^{1+n/2}\epsilon^{(1+v)/2}((n-1)(n-2)^v)^{-1}$ when $n$ exceeds $2$. For $n=1,2$, \[\nu^{-1}_\epsilon\exp(-\nu^{-1}_\epsilon(t_2-t_1))(t_2-t_1)^{n}(1-b)|b^{n\zeta}-b^{n\zeta^\prime}|\leq C\min\{1,(t_2-t_1)^{-(1+v)/2}\}\epsilon^{(1+v)/2}|\zeta-\zeta^\prime|^v\] for small enough $\epsilon$. Once again, using the property of the confidence band of Poisson distribution,
  \[\sum_{n=2}^\infty\exp(-\nu^{-1}_\epsilon(t_2-t_1))\frac{2^{n/2}\nu^{-n}_\epsilon(t_2-t_1)^n}{n!(n-1)(n-2)^v}\leq C^\prime\min \{1,(t_2-t_1)^{-(1+v)/2}\},\]
 one can bound right side in \eqref{eq:TPDiffCrudeBd} by $C|\zeta-\zeta^\prime|^v\epsilon^{(1+v)/2}\min \{1, (t_2-t_1)^{-(1+v)/2}\}$. This completes the proof.
\ee
\end{proof}

For the notational convenience, we use $Z(t,x)$ instead of $Z_\epsilon(t,x)$ for the next two lemmas. Let us note that one can modify the scaling appropriately to get necessary results for $Z_\epsilon(t,x)$. In a nutshell, the following results portray the key role in showing the tightness. One can also find similar result in the context of ASEP \cite[Lemma~4.2,~4.3]{Bertini1995}, higher spin exclusion process \cite[Lemma~4.3]{CL15}. Recall the decomposition of $Z(t,x)$ in \eqref{eq:DiscreteSHE}. We define $Z_{mg}(t_1,t_2,x)$  and $Z_{\nabla, mg}(t_1,t_2,x)$ as follows 
\begin{align}\label{eq:NewTwoDef}
Z_{mg}(t_1,t_2,x):&= \int_{t_1}^{t_2}[p_\epsilon(s,t_2)*Z(s)dM(s)](x)\\
Z_{\nabla, mg}(t_1,t_2,x,x^\prime):&= \int_{t_1}^{t_2}[p_\epsilon(s,t_2)*Z(s)dM(s)](x)-\int_{t_1}^{t_2}[p_\epsilon(s,t_2)*Z(s)dM(s)](x^\prime)
\end{align}
In the next lemma, we probe into the higher order norms of the processes $Z_{mg}(t_1,t_2,x)$ and $Z_{\nabla,mg}(t_1,t_2,x)$.

\bl\label{NormEstimateLemma}
For any $k\geq 1$ and $v\in (0,1]$, there exists a large constant $C$ which depends only on $k$ such that for all $t_1\leq t_2\in \NN$ and $\zeta,\zeta^\prime\in \Xi(t_1,t_2)$,
\begin{align}
\|Z_{mg}(t_1,& t_2,\zeta)\|^2_{2k}\\&\leq C\epsilon^{1/2}\int_{t_1}^{t_2}\min \{1,(t_2-s)^{-1/2}\}\left[(\delta_0+P_\epsilon)*[p_\epsilon(s,t)*\|Z(s)\|^2_{2k}]\right](\zeta)ds\label{eq:NormEstimate1} \\
\|Z_{\nabla, mg}&(t_1,t_2,\zeta,\zeta^\prime)\|^2_{2k}\\ &\leq C \epsilon^{(1+v)/2} |\zeta-\zeta^\prime|^v \int_{t_1}^{t_2}\min \{1,(t_2-s)^{-(1+v)/2}\}\left(\left[(\delta_0+P_\epsilon)*[p_\epsilon(s,t)*\|Z(s)\|^2_{2k}]\right](\zeta)\right.\nonumber\\&\left.+\left[(\delta_0+P_\epsilon)*[p_\epsilon(s,t)*\|Z(s)\|^2_{2k}]\right](\zeta^\prime)\right)ds\label{eq:NormEstimate2}
\end{align}
where $\delta_0$ is the probability mass function of the distribution which is degenerate at $0$ and  $P_\epsilon$ denotes probability mass function of the random variables $R_i$ (see, \eqref{eq:RSupportMassFunction} for definition).
\el

\begin{proof}
  Fix any $t\in (t_1,t_2)$. Using Burkholder-Davis-Gundy's inequality, we have
  \begin{align}\label{eq:ZmgBDG}
  |Z_{mg}(t_1,t,\zeta)|^2_{2k}\leq \|[Z_{mg}(t_1,.,\zeta),Z_{mg}(t_1,.,\zeta)]_t\|_{k},
  \end{align}
  and
  \begin{align}\label{eq:ZNabMgBDG}
  |Z_{\nabla, mg}(t_1,t,\zeta)|^2_{2k}\leq \|[Z_{\nabla, mg}(t_1,.,\zeta),Z_{\nabla, mg}(t_1,.,\zeta)]_t\|_{k},
  \end{align}
  where $[.,.]$ denotes the quadratic variation. Let us denote set all of all time points at which particles hop across $\zeta,\zeta^\prime$ in the interval $(s_1,s_2]$ by $\mathcal{S}_{(s_1,s_2]}(\zeta,\zeta^\prime)$. Recall once again the form of the Martingale in \eqref{eq:DefMartingale} and the function $\kappa(.,.)$ in \eqref{eq:CharacterisLineOfKPZ}. Define \begin{align}\label{eq:HoppingEfectInTotVer}
  \mathcal{Q}(\xi_1,\xi_2,y,s) = (b^{-1}-1)^{2}\sum_{y\leq \kappa(\xi_1\wedge \xi_2,s)}\eta_s(y) Z^{-1}(\xi_1\vee\xi_2,s)b^{\nu_\epsilon\kappa(\xi_1\vee \xi_2,s)+N_y(s)-y}
  \end{align}
  which controls the contribution coming out of the Martingale $M(t,x)$ in the total variation terms on right side of \eqref{eq:ZmgBDG}
and \eqref{eq:ZNabMgBDG}. Using interdependence of the Poisson processes $\{\mathcal{P}_l\}_{-\infty <l<\infty}$, one can say $\mathcal{Q}(\xi_1,\xi_2,y,s)$ captures the hopping effect of particles at any position $y$ across $\kappa(\xi_1,s)$ and $\kappa(\xi_2,s)$ at time $s$. Thus, we have
  \begin{align}\label{eq:ZmgExpre}
  [Z_{mg}(t_1,.,\zeta),Z_{mg}(t_1,.,\zeta)]_t = \sum_{\xi_1,\xi_2} \sum _{s\in \mathcal{S}_{(t_1,t]}(\xi_1,\xi_2)}p^{\xi_1,\xi_2}_\epsilon(t_1,t,\zeta)Z^{\xi_1,\xi_2}(s)\mathcal{Q}(\xi_1,\xi_2,y,s)
\end{align}
and
\begin{align}\label{eq:ZNabExpre}
  [Z_{\nabla ,mg}(t_1,.,\zeta,\zeta^\prime),Z_{\nabla, mg}(t_1,.,\zeta,\zeta^\prime)]_t = \sum_{\xi_1,\xi_2} \sum _{s\in \mathcal{S}_{(t_1,t]}(\xi_1,\xi_2)}p^{\xi_1,\xi_2}_{\nabla,\epsilon}(t_1,t,\zeta,\zeta^\prime)Z^{\xi_1,\xi_2}(s)\mathcal{Q}(\xi_1,\xi_2,y,s)\nonumber\\
\end{align}
where $Z^{\xi_1,\xi_2}(s)=Z(\xi_1,s)Z(\xi_2,s)$, $p_{\epsilon}^{\xi_1,\xi_2}(s,t_2,\zeta)=p_{\epsilon}(s,t_2,\zeta-\xi_1)p_{\epsilon}(s,t_2,\zeta-\xi_2)$ and \[p_{\nabla,\epsilon}^{\xi_1,\xi_2}(s,t,\zeta)=(p_{\epsilon}(s,t_2,\zeta-\xi_1)-p_{\epsilon}(s,t_2,\zeta^\prime-\xi_1))(p_{\epsilon}(s,t_2,\zeta-\xi_2)-p_{\epsilon}(s,t_2,\zeta^\prime-\xi_2)).\]


For the notational simplicity, for any two point $m>n$ on the lattice $\ZZ_{\geq 0}$ and $s\in \RR^{+}$
\begin{align}\label{eq:GapDefine}
h(m,n,s):= m-n -N_n(s) +N_m(s).
\end{align}
 On simplifying the summands in \eqref{eq:ZmgExpre} and \eqref{eq:ZNabExpre}, one can get
 \[Z(\xi_1,s)Z(\xi,s)\mathcal{Q}(\xi_1,\xi_2,y,s)=\epsilon Z(\xi_1\wedge \xi_2,s)^2\sum_{y\leq \kappa(\xi_1\wedge \xi_2 ,s)}\eta_s(y)b^{h(y,\kappa(\xi_1\wedge\xi_2,s),s)}b^{\nu_\epsilon|\xi_1-\xi_2|}\]
where $h(y,\kappa(\xi_1\wedge \xi_2, s),s)$ counts the number of holes in between $y$ and $\kappa(\xi_1\wedge \xi_2,s)$ at time $s$.
Notice that $p_{\epsilon}^{\xi_1,\xi_2}(s,t_2,\zeta)\leq C\sqrt{\epsilon} \max\{1, (t_2-s)^{-1/2}\} p_{\epsilon}(s,t_2,\zeta-\xi_1\wedge \xi_2)$ thanks to the Proposition \ref{KernelEstimate}. 
  For any fixed $\xi_1$, one can note $\sqrt{\epsilon}\sum_{\xi_2\geq \xi_1} b^{\nu_\epsilon|\xi_1-\xi_2|}$ is constant.
Thus, we bound right side of \eqref{eq:ZmgExpre} by
\[C\epsilon\min\{1, (t_2-s)^{-1/2}\}\sum_{\xi}\sum_{s\in \mathcal{S}_{(t_1,t]}(\zeta)}p_\epsilon(s,t_2,\zeta-\xi)Z(\xi,s)^2\sum_{y\leq \kappa(\xi,s)}\eta_s(y)b^{h(y,\kappa(\xi,s),s)}.\]
Divide the interval $(t_1,t]$ into disjoint intervals $\cap_{i=1}^m (s_i,s_{i+1}]$ where each of them is of length less than or equals to $ 1$. Next, recall that
\begin{align}\label{eq:TotVarAtSamePoint}
d\langle M(s,\xi),M(s,\xi)\rangle = \sum_{y\leq \kappa(\xi,s)}\eta_s(y)b^{h(y,\kappa(\xi,s),s)} dt.
\end{align}
Using Lemma~\ref{eq:TotVarDiff}, one can bound $Z^2(s,\xi)$ times the right side of \eqref{eq:TotVarAtSamePoint} by
\[\frac{1}{1-b^{\nu_\epsilon}}Z(s,\xi)\left([P_\epsilon*Z(s)](\xi)-\frac{(1-b^{\nu_\epsilon})b^{-1}}{b^{-1}-1}Z(s,\xi)dt\right).\]
 Thus, recalling that $b=e^{-\lambda_\epsilon\sqrt{\epsilon}}$, we get the following bound
\[(\nu_\epsilon\lambda_\epsilon \sqrt{\epsilon})^{-1}\left(\frac{[P_\epsilon*\|Z(s)\|^2_{2k}](\xi)+\|Z(s,\xi)\|^2_{2k}}{2}+\nu_\epsilon\|Z(s,\xi)\|^2_{2k}\right)\]
on the $k$-th norm of $Z(\xi,s)^2\sum_{y\leq \kappa(\xi,s)}\eta_s(y)b^{h(y,\kappa(\xi,s), s)}$.
  Also note that, \[\max_{s\in (s_i,s_{i+1}]}Z(s,x)\leq \exp\left(2\sqrt{\epsilon}\lambda_\epsilon N_{(s_i,s_{i+1}]}(x)+4 (\nu_\epsilon)^{-1}\right)Z(s_i,x)\]
and
  \[\min_{s\in (s_i,s_{i+1}]}Z(s,x)\geq \exp\left(-2\sqrt{\epsilon}\lambda_\epsilon N_{(s_i,s_{i+1}]}(x)-4 (\nu_\epsilon)^{-1}\right)\max_{s\in (s_i,s_{i+1}]}Z(s_i,x)\]
   where $N_{(s_i,s_{i+1}]}(x)$ is independent of $Z(s_i,x)$ and dominated by a Poisson point process of rate $\sum_{y\leq \kappa(x,s_{i+1})} \eta_{s_{i+1}}(y)b^{h(y, \kappa(x,s_{i+1}),{s_{i+1}})}$. Once again, using the same analysis as above one can show that the rate is bounded above by $C\epsilon^{-1/2}$. Thus, we have \[\mathbb{E}\left(\exp\left(2\sqrt{\epsilon}k\lambda_\epsilon N_{(s_i,s_{i+1}]}(x)+4k (\nu_\epsilon)^{-1}\right)\right)\leq C^\prime\]
where $C^\prime$ is a large constant which only depends on $k$.
To this end, we can bound right side of \eqref{eq:ZmgBDG} by
\[C\sqrt{\epsilon}\sum_{\xi}p_\epsilon(s,t_1,\zeta-\xi)\sum_{k=1}^m\min\{1,(t_2-s)^{-1/2}\}\inf_{s\in(s_i,s_{i+1}]}\left[[P_\epsilon*\|Z(s)\|^2_{2k}](\xi)+\|Z(s,\xi)\|^2_{2k}\right].\]
 Recall that the length of each of the intervals $(s_i,s_{i+1}]$ is less than or equals to $1$. Applying the bound $\inf_{s\in (s_i,s_{i+1}]}\|Z(\xi,s)\|_{2k}\leq (s_i-s_{k+1})^{-1}\int \|Z(\xi,s)\|_{2k} ds=\int \|Z(\xi,s)\|_{2k} ds$ and using the semi-group property of $p_\epsilon$, we get
\[ \Vert Z_{mg}(t_1,t_2,\zeta)\Vert_{2k}^2\leq C\epsilon^{1/2}\int_{t_1}^{t_2} \min\{1,(t_2-s)^{-1/2}\}\left\{\left[(\delta_0+P_\epsilon)*[p_\epsilon(s,t_2)*\|Z(s)\|^2_{2k}]\right](\zeta)\right\}ds\]
where $\delta_0$ is the probability mass function which puts all its mass in $0$. Note that we also had made of the fact that under convolution operation, $p_\epsilon(t_1,t,.)$ and $P_\epsilon$   commute with each other. This holds because $p_\epsilon$ is essentially the semi-group of the random walk of the iid random variables coming from the distribution $P_\epsilon$.

In the case of $Z_{\nabla,mg}$, using Proposition~\ref{eq:KernelEstimate4}, we can bound $p^{\xi_1,\xi_2}_{\nabla,\epsilon}(t_1,s,\zeta,\zeta^\prime)$ in \eqref{eq:ZNabExpre} by \[C\epsilon^{(1+v)/2}\min\{1,(t_2-s)^{-(1+v)/2}\}|\zeta-\zeta^{\prime}|^v \left(p_{\epsilon}(t_1,s,\zeta-\xi_1\wedge \xi_2)+p_{\epsilon}(t_1,s,\zeta^\prime-\xi_1\wedge \xi_2)\right).\] Rest of the terms in \eqref{eq:ZmgExpre} can be controlled exactly in the same fashion as we have done for bounding $\|Z_{mg}(t_1,t_2,\zeta)^2\|_{2k}$. Thus, one can write down
\begin{align*}
\|Z_{\nabla, mg}&(t_1,t_2,\zeta)\|^2_{2k}\\&\leq C\epsilon^{(1+v)/2}|\zeta-\zeta^\prime|^{v}\int_{t_1}^{t_2}\min\{1,(t_2-s)^{-(1+v)/2}\}\left\{\left[(\delta_0+P_\epsilon)*[p_\epsilon(s,t_2)*\|Z(s)\|^2_{2k}]\right](\zeta)\right.\\&\left.+ \left[(\delta_0+P_\epsilon)*[p_\epsilon(s,t_2)*\|Z(s)\|^2_{2k}]\right](\zeta^\prime)\right\}ds.
\end{align*}
\end{proof}

In the next result, we present a chaos-series type bound for $Z(t,\zeta)$.
\bl\label{ChaosSeriesLemma}
 For any $k\geq 1$, $t> 0$ and $\zeta \in \Xi(t)$,
\begin{align}\label{eq:ChaoSeriesBD}
\|Z(t,\zeta)\|^2&\leq 2([p_\epsilon(0,t)*\|Z(0)\|^2_{2k}](\zeta))\\&+ 2\sum_{n=1}^\infty \int_{\vec{s}\in \Delta_{n+1}(t)}h_\epsilon(s_1,s_2)\ldots h_\epsilon(s_n,s_{n+1})\Psi(s_{n+1},\zeta) ds_1ds_2\ldots ds_{n+1}
\end{align}
where $h_{\epsilon}(s_i,s_{i+1})=C\epsilon^{1/2} \min\{1,(s_{i}-s)^{-1/2}\}1_{s_i\geq s_{i+1}}$, \[\Psi_n(s_{n+1},\zeta)=\left[(\delta_0+P_\epsilon)^n*[p_\epsilon(0,s_{n+1})*\|Z^2(0)\|_{2k}]\right](\zeta)\] and $\Delta_n(t)=\{\vec{s}\in (\RR_{\geq 0})^n| t\geq s_1\geq \ldots \geq s_{n+1}\geq 0\}$. Here, $P_\epsilon$ denotes probability mass function of $R_i$ (see, \eqref{eq:RSupportMassFunction} for definition). Note that here $C$ is a finite constant depended only on $k$.
\el
\begin{proof}
Recall the decomposition of $Z(t,x)$ in \eqref{eq:DiscreteSHE}. Using one simple inequality, $|x+y|^2\leq 2(x^2+y^2)$, we have
\begin{align}\label{eq:ChaosSeriesBD1}
\|Z_\epsilon(t,\zeta)\|^2_{2k} &\leq 2(\|[p_\epsilon(0,t)*Z_\epsilon(0)]\|_{2k})^2 +2\left\Vert\int_0^t [p_\epsilon(s,t)*Z_\epsilon(s)dM(s)](\zeta)ds\right\Vert^2_{2k}.
\end{align}
We use Jensen's inequality to bound the first term on the right side of \eqref{eq:ChaosSeriesBD1} by $2([p_\epsilon(0,t)*\|Z_\epsilon(0)\|^2_{2k}])$. Moreover, one can also bound second term in \eqref{eq:ChaosSeriesBD1} using \eqref{eq:NormEstimate1}. Thus, we get


\begin{align}\label{eq:ChaosSeriesBD2}
\|Z_\epsilon(t,\zeta)\|^2_{2k} &\leq 2([p_\epsilon(0,t)*\|Z_\epsilon(0)\|^2_{2k}]) \\&+C \epsilon^{1/2}\int_0^t \min\{1,(t-s)^{-1/2}\}\left[(\delta_0+P_\epsilon)*[p_\epsilon(s,t)*\|Z_\epsilon(s)\|^2_{2k}]\right](\zeta)ds. 	
\end{align}
 To this end, using successive recursion on $\|Z(s,.)\|^2_{2k}$ in RHS of \eqref{eq:ChaosSeriesBD2}, one can get an asymptotic expansion. Furthermore, semi-group property $[p_\epsilon(t_1,s)*p_\epsilon(s,t_2)](.)=p_\epsilon(t_1,t_2,.)$ and the fact that $P_\epsilon$ commutes with the semi-group $\{p_\epsilon\}$ helps in concluding the final step of the inequality.
\end{proof}

Our next goal is to prove Proposition~\ref{Tightness} with some refined estimates of the norms of  $Z_\epsilon(t,x)$. Below, we illustrate the idea in details.  

\textit{Proof of Proposition~\ref{Tightness}}:
 First, we prove the following moment estimates:
 \begin{align}
 \|Z_\epsilon(t,x)\|_{2k}&\leq Ce^{\epsilon\tau |x|}, \label{eq:MomentEstimate1}\\
 \|Z_\epsilon(t,x)-Z_\epsilon(t,x^\prime)\|_{2k}& \leq C  \left(\epsilon |x-x^\prime|\right)^{v}e^{\tau\epsilon(|x|+|x^\prime|)},\label{eq:MomentEstimate2}\\
 \|Z_\epsilon(t,x)-Z_\epsilon(t^\prime, x)\|_{2k}& \leq C \left(\epsilon |t-t^\prime|\right)^{\alpha}e^{\tau\epsilon|x|}.\label{eq:MomentEstimate3}
 \end{align}
 for some $\alpha:=v/2\wedge 1/4$, $t,t^\prime \in (t_1,\epsilon^{-1}T]$, $x,x^\prime\in \Xi(t)$ where $T\in (0,\infty)$ is fixed apriori and $v$ is exactly same as in \eqref{eq:EquiInitCond2}. Once those are proved, we can make use of Kolmogorov-Chentsov criteria to conclude the tightness of the sequence $\{Z_\epsilon\}_\epsilon$.

 \textit{Proof of \eqref{eq:MomentEstimate1}}:
  We have $\|Z(0,x)\|_{2k}\leq e^{\tau \epsilon|x|}$ from \eqref{eq:EquiInitCond1}. One can further bound $|\xi|$ by $|x-\xi|+|x|$. Thus, we have
  \[[p_\epsilon(0,t)*\|Z(0)\|^2_{2k}](x)\leq \sum_{\xi\in \Xi(t)}p_\epsilon(0,t,x-\xi)\exp(2\tau\epsilon|x-\xi|)\exp(2\tau\epsilon|x|)\leq C\exp(2\tau \epsilon|x|).\]
  thanks to \eqref{eq:KernelEstimate1}. Consequently, first term in RHS of \eqref{eq:ChaoSeriesBD} is bounded above by $C\exp(2\tau\epsilon|x|)$. Using the same arguments, one can also bound $[p_\epsilon(0,s_{n+1})*\|Z(0)\|^2_{2k}](x)$ of the second term in \eqref{eq:ChaoSeriesBD} by $C_1\exp(2\tau\epsilon|x|)$. Furthermore, we know
\[\sum_{\xi=-\infty}^{x}P_\epsilon(x-\xi)\exp(2\tau \epsilon|\xi|)\leq \exp(2\tau\epsilon|x|)\sqrt{\epsilon}\sum_{k=0}^\infty\exp(-k\sqrt{\epsilon}(\nu_\epsilon-2\sqrt{\epsilon}\tau))\leq C_2\exp(2\tau\epsilon|x|)\]
for some constant $C_2$ when $\epsilon$ is small enough.
Thus, we can conclude that $[(\delta_0+P_\epsilon)^n*[p_\epsilon(0,s_{n+1})*Z(s)]](x)$ is bounded by $C_1C_2^n \exp(2\tau\epsilon|x|)$ for all small $\epsilon$.
     Interestingly, we can now easily compute left over integral of the second term in the RHS of \eqref{eq:ChaoSeriesBD}. It turns out that
  \begin{align}\label{eq:IteratedIntEst}
  \int_{\vec{s}\in \Delta_{n+1}(t)} \prod_{i=1}^n h_\epsilon(s_i,s_{i+1})ds_{1}ds_{2}\ldots ds_{n+1} \leq \frac{\left((C\Gamma(1/2))^2\epsilon t\right)^{n/2}}{\Gamma(n/2)}.
\end{align}
This makes the sum in \eqref{eq:ChaoSeriesBD} finite and thus helps us to conclude \eqref{eq:MomentEstimate1}.

 \textit{Proof of \eqref{eq:MomentEstimate2}}:
    We can decompose $Z(t,x)-Z(t,x^\prime)$ as
   \begin{align}\label{eq:PosDiffDecom}
   Z(t,x)-Z(t,x^\prime)&=\sum_{\xi\in \mathfrak{t}} p_\epsilon(0,t,\xi)\left(Z(t_1,x-\xi)-Z(t_1,x^\prime-\xi)\right)\nonumber\\
   &+\int_0^t \sum_{\xi\in \Xi(s,t)}\left(p_\epsilon(s,t,x-\xi)-p_\epsilon(s,t,x^\prime-\xi)\right)Z(s,\xi)dM(s,\xi)
   \end{align}
   Notice that we have the bound on $p_\epsilon(0,t,\zeta)-p_\epsilon(0,t,\zeta^\prime)$ from Lemma~\ref{KernelEstimate}. One can use \eqref{eq:EquiInitCond2} to bound the first term on the RHS of \eqref{eq:PosDiffDecom} as noted down below.
   \begin{align}\label{eq:PosDiffEstAt0}
   \| [p(0,t)*Z(0)]&(x)-[p(0,t)*Z(0)](x^\prime)\|_{2k}\nonumber\\&\leq \sum_{\xi\in \Xi(0,t)}p_\epsilon(0,t,\xi)\|Z(0,x-\xi)-Z(0,x^\prime-\xi)\|_{2k}\nonumber\\
    &\leq \sum_{\xi \in \Xi(0,t)}p_\epsilon(0,t,\xi)\left(\epsilon|x-x^\prime|\right)^v \exp\left(\tau\epsilon(|x-\xi|+|x^\prime-\xi|)\right)\nonumber\\
   &\leq \left(\epsilon|x-x^\prime|\right)^v\exp\left(\tau\epsilon(|x|+|x^\prime|)\right)\sum_{2\xi \in \Xi(0,t)}p_\epsilon(0,t,\xi)\exp(2\tau\epsilon|\xi|)\nonumber\\
   &\leq C\left(\epsilon|x-x^\prime|\right)^v\exp\left(\tau\epsilon(|x|+|x^\prime|)\right).
   \end{align}
    One can identify the term in the last line of \eqref{eq:PosDiffDecom} with $ Z_{\nabla,mg}
   (t_1,t,x,x^\prime)$. It follows from \eqref{eq:MomentEstimate1} that $\|Z(s,\xi)\|_{2k}\leq C\exp(\tau
   \epsilon |\xi|)$ for all $\xi\in \Xi(t_1,s)$ and $0\leq s\leq \epsilon^{-1}T$ where $C$ only depends on
   $T$. Using the bound $|x-\xi|+|x|$ for $|\xi|$ in $\exp(\tau\epsilon|\xi|)$, one can bound $[p_
   \epsilon(s,t)*\|Z(s)\|^2_{2k}](x)$ above by $C\exp(2\tau \epsilon|x|)$. Similarly, $C\exp(2\tau
   \epsilon|x^\prime|)$ is also valid as an upper bound to $[p_\epsilon(s,t)*\|Z(s)\|^2_{2k}](x^
   \prime)$. So, the integrand in \eqref{eq:NormEstimate1} can be bounded above by $\min\{1,(t_2-s)^{-(1+v)/2}\}1_{s\leq t_2}C
   \exp\left(2\tau\epsilon(|x|+|x^\prime|)\right)$ where $t_2$ is equal to $t$ and $s\in(0,t)$.
   Consequently, we can bound the $\|Z_{\nabla,mg}(0,t,x,x^\prime)\|_{2k}$ by
   $C(\epsilon|x-x^\prime|)^v (\epsilon t)^{(1-v)/2}\exp\left(\tau \epsilon(|x|+|x^\prime|)\right)$. This
   proves \eqref{eq:MomentEstimate2} thanks to the fact $t\in (0,\epsilon^{-1}T]$.

   \textit{Proof of \eqref{eq:MomentEstimate3}}:
        Assume $ t^\prime<t$. We can write down the following:
        \begin{align}\label{eq:TimeDiffDecom}
        Z(t,x) - Z(t^\prime,x)&= \sum_{\xi\in \Xi(0,t)}p(0,t,\xi)Z(t_1,x-\xi)-\sum_{\xi^\prime\in \Xi(t^\prime)}p(0,t^\prime,\xi^\prime)Z(0,x-\xi^\prime)\nonumber\\
        &+\int_{t^\prime}^t [p_\epsilon(s,t)*Z(s)dM(s)](x).
        \end{align}
        Notice that the first term in the last line of \eqref{eq:TimeDiffDecom} is exactly $ Z_{ mg}(t^\prime,t,x)$. Furthermore, using semigroup property of $p_\epsilon$, one can say
        \begin{align}\label{eq:LittleExpansion}
        \sum_{\xi\in \Xi(0,t)}p_\epsilon(t_1,t,\xi)Z(0,x-\xi)&=\sum_{\xi\in \Xi(t)}\sum_{\xi^\prime\in \Xi(t^\prime)}p_\epsilon(0,t^\prime,\xi^\prime)p_\epsilon(t^\prime,t,\xi-\xi^\prime)Z(0,x-\xi)\nonumber\\
        &=\sum_{\xi^\prime\in \Xi(0,t^\prime)}p_\epsilon(0,t^\prime,\xi^\prime)\sum_{\zeta\in \Xi(t^\prime,t)}p_\epsilon(t^\prime,t,\zeta)Z(0,x-\xi^\prime-\zeta).
        \end{align}
   This helps us to analyse the first term in RHS of \eqref{eq:TimeDiffDecom} in the following way.
   \begin{align}
   \left\Vert\sum_{\xi\in \Xi(t)}p(0,t,\xi)\right.&\left.Z(t_1,x-\xi)-\sum_{\xi^\prime\in \Xi(t^\prime)}p(0,t^\prime,\xi^\prime)Z(0,x-\xi^\prime)\right\Vert_{2k}\nonumber\\
   &= \sum_{\xi^\prime\in \Xi(t^\prime)}p_\epsilon(0,t^\prime,\xi^\prime)\left[\sum_{\zeta\in \Xi(t^\prime,t)} p_\epsilon(t^\prime,t,\zeta)\left\Vert Z(0,x-\xi^\prime-\zeta)-Z(0,x-\xi^\prime)\right\Vert_{2k}\right]\nonumber\\
   &\leq C\sum_{\xi^\prime\in \Xi(t^\prime)}p_\epsilon(0,t^\prime, \xi^\prime)\sum_{\zeta\in \Xi(t^\prime,t)}\epsilon|\zeta|^v p_\epsilon(t^\prime,t,\zeta)\exp\left(\tau\epsilon(|\zeta|+2|x|+2|\xi^\prime|)\right)\nonumber\\
   &\leq C(\epsilon|t-t^\prime|)^{v/2} \exp(2\tau \epsilon|x|).
\end{align}
 Note that we have used \eqref{eq:KernelEstimate4} in Proposition~\ref{KernelEstimate} in the last inequality. One can bound $\|Z_{mg}(t^\prime,t,x)\|_{2k}$ using \eqref{eq:NormEstimate1}. Moreover, integrand in the RHS of \eqref{eq:NormEstimate1} is bounded above by $$C\sqrt{\epsilon}1_{t^\prime \leq s\leq t}\min\{1,(t^\prime-s)^{-1/2}\}\exp(2\tau\epsilon|x|)$$ thanks to \eqref{eq:MomentEstimate1}. This enforces $\|Z_{mg}(t^\prime,t,x)\|^2_{2k}$ to be bounded by $C(\epsilon|t-t^\prime|)^{1/2}\exp(2\tau\epsilon|x|)$. Hence, the claim follows.

 \textit{Proof of Proposition~\ref{MomentEstimateProp}:}  First, we show \eqref{eq:MomemtEstimateFromStep1}. In Lemma~\ref{ChaosSeriesLemma}, we multiply both side with $\epsilon^{-2}(1-\exp(-\lambda_\epsilon\nu_\epsilon))^2$. Consequently, we get
  \begin{align}\label{eq:MomentEstimatePropChaosSeries}
  \|\tilde{Z}(t,\zeta)\|^2_{2k}\leq & 2\left( [p_\epsilon(0,t)*\tilde{Z}(0)](\zeta)\right)^2\\
  &+2\sum_{n=1}^{\infty} \int_{\vec{s}\in \Delta_{n+1}(t)} h(s_1,s_2)\ldots h(s_n,s_{n+1})\tilde{\Psi}_n(s_{n+1}, \zeta) ds_1  \ldots ds_{n+1}
\end{align}
where 
\begin{align}\label{eq:PsiTildeDef}
\tilde{\Psi}_n(s_{n+1}, \zeta):=\left[(\delta_0+P_\epsilon)^n*([p_\epsilon(0,s_{n+1})*\tilde{Z}(0)])^2\right](\zeta).
\end{align}

 Note that $\tilde{Z}(0,.)$ is now deterministic. Using the relation $\epsilon \sum_{\zeta\in \Xi(0)} \tilde{Z}(0)=1$ and the upper bound on $p_\epsilon(0,t,.)$ from \eqref{eq:KernelEstimate3}, we further have
\begin{align}\label{eq:ConvolutionBoundAt0}
\left|[p_\epsilon(0,t)*\tilde{Z}(0)](x)\right|\leq C\min\{\epsilon^{-1/2}, (\epsilon t)^{-1/2}\}.
\end{align}
Similarly, the term $[p_\epsilon(0,s_{n+1})*\tilde{Z}(0)]$ inside $\tilde{\Psi}_n(s_{n+1},\zeta)$ in \eqref{eq:PsiTildeDef} is also bounded above by $\min\{\epsilon^{-1/2}, (\epsilon s_{n+1})^{-1/2}\}$.  To that effect, one can bound the rest of the integral as
\begin{align}\label{eq:ItereatedIntegralRevisited}
\int_{\vec{s}\in \Delta_{n+1}(t)} h(s_1,s_2) \ldots h(s_n,s_{n+1}) (\epsilon s_{n+1})^{-1/2} ds_1\ldots d_{s_{n+1}}\leq  \frac{\left((C\Gamma(1/2))^2\epsilon t\right)^{n/2}}{\Gamma(n/2)}
\end{align}
in the same way as in \eqref{eq:IteratedIntEst}. Thus, for all $t\in (0,\epsilon^{-1}T]$, summing \eqref{eq:ItereatedIntegralRevisited} over $n\in \ZZ_{+}$, we get
\begin{align}\label{eq:MomentEstimate1EndCalcul}
\|\tilde{Z}(t,x)\|^2_{2k}\leq C([p_\epsilon(0,t)*\tilde{Z}_\epsilon(0)](\zeta))^2+\exp((\epsilon t)^{1/2})\leq C^\prime\min\{\epsilon^{-1}, (\epsilon t)^{-1} \}
\end{align}
where the constants $C,C^\prime>0$ depend only on $T$. Thus the claim follows.

 Now, we turn to show \eqref{eq:MomemtEstimateFromStep2}. To this end, multiplying $Z(t,\zeta)-Z(t,\zeta^\prime)$ by $\epsilon^{-1}(1-\exp(-\lambda_\epsilon\nu_\epsilon))$, we get
 \begin{align}\label{eq:ZDiffEq}
 \|\tilde{Z}(t,\zeta)-\tilde{Z}(t,\zeta^\prime)\|^2_{2k}&\leq 2\left(\sum_{\xi \in \Xi(0)}|p_\epsilon(0,t,\zeta -\xi)-p_\epsilon(0,t,\zeta-\xi)|\tilde{Z}(0,\xi)\right)^2\\
 &+ c\|\tilde{Z}_{\nabla, mg}(0,t,\zeta_1,\zeta^\prime_1)\|^2_{2k}
 \end{align}
  for some constant $c$. Using the bound on $|p_\epsilon(0,t,\zeta -\xi)-p_\epsilon(0,t,\zeta-\xi)|$ from \eqref{eq:KernelEstimate4} and the relation $\epsilon\sum_{\zeta\in \Xi(0)}\tilde{Z}(0,\zeta)=1$, one can get the following
 \begin{align}\label{eq:FirstTerBoundOfZDiff}
 \sum_{\xi \in \Xi(0)}|p_\epsilon(0,t,\zeta -\xi)-p_\epsilon(0,t,\zeta-\xi)|\tilde{Z}(0,\xi)\leq C(\epsilon|\zeta-\zeta^\prime|)^v\min\{\epsilon^{-\frac{1+v}{2}}, (\epsilon t)^{-\frac{1+v}{2}}\}.
\end{align}
 Moreover, substituting $2v$ in place of $v$ and taking $(t_1,t_2)=(0,t)$ in \eqref{eq:NormEstimate2}, we get
 \begin{align}
 \|\tilde{Z}_{\nabla,mg}(0,t,\zeta,\zeta^\prime)\|^2_{2k}\leq &\epsilon^{\frac{1+2v}{2}}|\zeta-\zeta|^{2v}\int \min\{1, (t-s)^{-(1+2v)/2}\}\left[\tilde{\Psi}(s,t,\zeta)+\tilde{\Psi}(s,t,\zeta^\prime)\right]ds\\
 \label{eq:NormEstimateSurgery}
 \end{align}
 where $\tilde{\Psi}(s,t,\zeta):=\left[(\delta_0+P_\epsilon)*[p_\epsilon(s,t)*\|\tilde{Z}(s)\|^2_{2k}]\right](\zeta)$. 
 One can further bound the integrand on the right side in \eqref{eq:NormEstimateSurgery} using bound on $\|\tilde{Z}(s,\zeta)\|^2_{2k}$ from \eqref{eq:MomentEstimate1EndCalcul}. Henceforth, this implies
 \begin{align}\label{eq:FinalBoundOnNablaDiff}
 \|\tilde{Z}_{\nabla,mg}(0,t,\zeta,\zeta^\prime)\|^2_{2k}\leq C(\epsilon|\zeta-\zeta^\prime|)^{2v} (\epsilon t)^{-v+(1/2)}.
\end{align}
 But, we know $\epsilon t\leq T$ for all $t\in (0,\epsilon^{-1}T]$. Thus, we can improve the bound on the right side of \eqref{eq:FinalBoundOnNablaDiff} to $C^\prime (\epsilon|\zeta-\zeta|)^{2v}(\epsilon t)^{-(1+v)}$ where $C^\prime$ depends only on $T$.

\subsection{Proof of Proposition~\ref{SubseqArguements}: The Martingale Problem}
   Before going into the details of the proof, we present a brief exposure on the martingale problem from \cite{Bertini1997}.
 \bd
 Let $\mathcal{Z}$ be a $C([0,\infty)\times (-\infty,\infty))$ valued process such that for any given $T>0$, there exists $u<\infty$ such that
 \begin{equation}\label{eq:UniformMoment}
  \sup_{t \in [0,T]}\sup_{x\in \mathbb{R}}e^{-u|x|}\mathbb{E}(\mathcal{Z}(t,x))<\infty.
\end{equation}
For such $\mathcal{Z}$ and $\psi\in C^\infty_c(\RR)$, let $\langle \mathcal{Z}(t),\psi \rangle :=\int_{\RR}\mathcal{Z}(t,x)\psi(x)dx$. We say $\mathcal{Z}$ solves the martingale problem with initial condition $\mathcal{Z}^{ic}\in C(\mathbb{R})$, if $\mathcal{Z}(0,.)=\mathcal{Z}^{ic}(.)$ in distribution and
\begin{align*}
t\mapsto N_\psi(t)&:= \langle \mathcal{Z}(t),\psi\rangle -\langle \mathcal{Z}(0),\psi\rangle -\frac{1}{2}\int_{0}^t \left\langle \mathcal{Z}(\tau),\frac{\partial^2 \psi}{\partial x^2}\right\rangle d\tau\\
t\mapsto N^\prime_\psi(t)&:= (N_\psi(t))^2-\int_0^t \left\langle \mathcal{Z}^2(\tau), \psi^2\right\rangle d\tau
\end{align*}
are local martingales for any $\psi\in C^{\infty}_c(\mathbb{R})$.
 \ed
 It has been stated in \citep{Bertini1997}, that for any initial condition $\mathcal{Z}^{ic}$ satisfying
 \begin{equation}\label{eq:SubExponentialDecay}
 \|\mathcal{Z}^{ic}(x)\|_{2k}\leq Ce^{a|x|}\quad \text{for some }a>0,
 \end{equation}
 the solution of the martingale problem stated above coincides with the solution of the SHE with initial condition $\mathcal{Z}^{ic}$. In our case, if we can show any limit point of $\{Z_\epsilon\}_\epsilon$ started from $\mathcal{Z}^{ic}$ solves the martingale problem, then it follows from \citep{Bertini1997} that $\{Z_\epsilon\}_\epsilon$ converges to a unique process which solves SHE with initial condition being $\mathcal{Z}^{ic}$. For solving the martingale problem, we have to have a discrete analogue of $\langle \mathcal{Z}(t),\psi\rangle$ for any $Z_\epsilon$. Note that we can take
 \[\langle Z(s),\psi \rangle_\epsilon=\epsilon \sum_{\xi \in \ZZ}Z(t,\xi)\psi(\epsilon \xi)\]
 as the discrete analogue for $\langle Z_\epsilon(t),\psi\rangle$. Let us divide $(0,\infty)$ into a number of disjoint sub-intervals $\cup_{i=1}^\infty (s_{i},s_{i+1})$ where $0=s_1<s_2<\ldots$ and each has length $\epsilon$. Let us assume that $s_m\leq t<s_{m+1}$. On the basis of decomposition of $Z(t,x)$ in \eqref{eq:DiscreteSHE}, one can say the following:
 \begin{align}\label{eq:MartingaleProblemDecom}
 \langle Z(t),\psi\rangle_\epsilon -\langle Z(0),\psi\rangle_\epsilon &= \epsilon \sum_{i=1}^{m-1}\sum_{\zeta\in \Xi(s_i)}\left[\sum_{\xi\in \Xi(s_i,s_{i+1})+\zeta}p_\epsilon(s_i,s_{i+1}, \xi-\zeta)\psi(\epsilon\xi)-\psi(\epsilon\zeta)\right]Z(s_i,\zeta)\nonumber\\
 &+\epsilon\sum_{\zeta\in \Xi(s_i)}\left[\sum_{\xi\in \Xi(s_m,t)+\zeta}p_\epsilon(s_m,t, \xi-\zeta)\psi(\epsilon\xi)-\psi(\epsilon\zeta)\right]Z(s_m,\zeta)\nonumber\\
 &+ \epsilon\sum_{i=1}^{m-1}\int_{s_i}^{s_{i+1}} \sum_{\zeta\in \Xi(s)}\left(\sum_{\xi\in \Xi(s,s_{i+1})+\zeta}p_\epsilon(s,s_{i+1},\xi-\zeta)\psi(\epsilon\xi)\right)Z(s,\zeta)dM(s,\zeta)\nonumber\\
 &+\epsilon\int_{s_m}^{t}\sum_{\zeta\in \Xi(s)}\left(\sum_{\xi\in \Xi(s,t)+\zeta}p_\epsilon(s,s_{i+1},\xi-\zeta)\psi(\epsilon\xi)\right)Z(s,\zeta)dM(s,\zeta).
 \end{align}

 Let us define $R_\epsilon(t)$ as
 \begin{align}\label{eq:R1Def}
 R_\epsilon(t):=& \epsilon \sum_{i=1}^{m-1}\sum_{\zeta\in \Xi(s_{i})}\left[\sum_{\xi\in \Xi(s_i,s_{i+1})+\zeta}p_\epsilon(s_{i},s_{i+1}, \xi-\zeta)\psi(\epsilon\xi)-\psi(\epsilon\zeta)\right]Z(s,\zeta)ds\\
 &+\epsilon\sum_{\zeta\in \Xi(s_i)}\left[\sum_{\xi\in \Xi(s_m,t)+\zeta}p_\epsilon(s_m,t, \xi-\zeta)\psi(\epsilon\xi)-\psi(\epsilon\zeta)\right]Z(s_m,\zeta).
 \end{align}
 Also, define
 \begin{align}\label{eq:NepDef}
 N_\epsilon(t):=&\epsilon\sum_{i=1}^m\int_{s_i}^{s_{i+1}} \sum_{\zeta\in \Xi(s)}\left(\sum_{\xi\in \Xi(s,s_{i+1})}p_\epsilon(s,s_{i+1},\xi-\zeta)\psi(\epsilon\xi)\right)Z(s,\zeta)dM(s,\zeta)\\&+\epsilon\int_{s_m}^{t}\sum_{\zeta\in \Xi(s)}\left(\sum_{\xi\in \Xi(s,t)+\zeta}p_\epsilon(s,s_{i+1},\xi-\zeta)\psi(\epsilon\xi)\right)Z(s,\zeta)dM(s,\zeta).
 \end{align}

Clearly, $N_\epsilon(t)$ is a local martingale. Let us denote the quadratic variation of $N_\epsilon(t)$ by $\hat{N}_\epsilon(t)$. So, for showing that indeed any limit point of $\{Z_\epsilon\}_\epsilon$ solves the martingale problem, it is enough to show the following:

  \begin{align}
  \mathbb{E}\left|\langle Z(\epsilon^{-1}t),\psi\rangle_\epsilon - \langle Z_\epsilon(t),\psi\rangle\right|&\to 0,\label{eq:MartingaleCondOne}\\ \mathbb{E}\left|R_\epsilon(\epsilon^{-1}t)-\int_0^t \left\langle Z_\epsilon(s), \frac{\partial^2 \psi}{\partial x^2}\right\rangle ds\right|&\to 0,
  \label{eq:MartingaleCondTwo}\\
  \mathbb{E}\left|\hat{N}_\epsilon(\epsilon^{-1}t)-\int_0^t\langle Z^2_\epsilon(s),\psi^2\rangle ds\right|&\to 0.\label{eq:MartingaleCondThree}
  \end{align}
   To prove \eqref{eq:MartingaleCondOne} - \eqref{eq:MartingaleCondThree}, we need the following lemma.

  \textit{Proof of \eqref{eq:MartingaleCondOne}}:
     Here, we have to show that
     \[\langle Z_\epsilon(t),\psi\rangle =\int_{\RR} Z_\epsilon(t, x)\psi(x)dx\]
    and $\langle Z(\epsilon^{-1}t),\psi \rangle_\epsilon$ are asymptotically same in $L^{1}$. Notice that for $x$ such that $|x-\zeta|\leq 1$ where $\zeta\in \Xi(\epsilon^{-1}t)$, using smoothness of $\psi$ and moment estimates from \eqref{eq:MomentEstimate1}-\eqref{eq:MomentEstimate3}, one can say that the sequences $\epsilon^{-v}\left(Z(\epsilon^{-1}t,\zeta)\psi(\epsilon\zeta)-Z(\epsilon^{-1}t,x)\psi(\epsilon x)\right)$ are $L^{2k}$ bounded for all $k\in \NN$. This forces $\langle Z(\epsilon^{-1}t),\psi\rangle_\epsilon - \langle Z_\epsilon(t),\psi\rangle$ to go to $0$ in $L^1$ as $\epsilon\to 0$.

   \textit{Proof of \eqref{eq:MartingaleCondTwo}}:
     Recall that subintervals $\{[s_i,s_{i+1})\}^\infty_{i=1}$ form a partition of $(0,\infty)$ such that length of each of the interval is $\epsilon$. Assume here that $\epsilon s_m\leq t<\epsilon s_{m+1}$. Using the smoothness of $\psi$, one can get
     \begin{align}\label{eq:UseOfSmoothness}
     \sum_{\mathclap{\xi\in \Xi(s_{i},s_{i+1})+\zeta}}  p_\epsilon(s_{i},&s_{i+1}, \xi-\zeta)\psi(\epsilon\xi)-\psi(\epsilon\zeta)
     =\epsilon\psi^\prime (\epsilon\zeta)\sum_{\mathclap{\xi\in \Xi(s_{i},s_{i+1})}}p_\epsilon(s_{i},s_{i+1},\xi-\zeta)(\xi-\zeta)\nonumber\\&+\sum_{\mathclap{\xi\in \Xi(s_{i},s_{i+1})}}\epsilon^2\frac{1}{2}\left(\frac{\partial^2 \psi}{\partial x^2}(\epsilon\zeta)+o(\epsilon)\mathcal{B}(\zeta)\right)p_\epsilon(s_{i},s_{i+1},\xi-\zeta)(\xi-\zeta)^2
     \end{align}
     where $\mathcal{B}$ belongs to $L^2$. Note that $\sum_{\xi \in \Xi(s,t)}p_\epsilon(s_{i},s_{i+1},\xi-\zeta)(\xi-\zeta)=0$. This goes back to the fact that mean under $p_\epsilon(s_{i},s_{i+1},.)$ is $0$. Moreover, $\sum_{\xi\in \Xi(s_i,s_{i+1})}p_\epsilon(s_i,s_{i+1},\xi-\zeta)(\xi-\zeta)^2=(s_{i+1}-s_{i})\sigma^2_\epsilon=(s_{i+1}-s_{i})\frac{b^{-1}-1}{1-b^{\nu_\epsilon}}\frac{1}{(1-b^{\nu_\epsilon})^2}=(s_{i+1}-s_{i})(\nu^{-1}_\epsilon+O(\epsilon))(\nu_\epsilon\lambda_\epsilon)^{-2}\epsilon^{-1}$. Recall from Definition~\ref{GartnerTransform} that $\lambda_\epsilon=\nu^{-3/2}_\epsilon$. To this effect, we have
     \begin{align}\label{eq:Repestimates}
     R^{1}_\epsilon(\epsilon^{-1}t)&=\epsilon\frac{1}{2}\sum_{i=1}^m(s_{i+1}-s_i)\epsilon\sum_{\zeta\in \Xi(s_{i})}Z(s_{i},\zeta)\left(\frac{\partial^2 \psi}{\partial x^2}(\epsilon\zeta)+o(\epsilon)\right)\\
     &=\frac{1}{2}\epsilon\int_0^{t\epsilon^{-1}} \left\langle Z(s),\frac{\partial^2 \psi}{\partial x^2}\right\rangle_\epsilon ds + o(\epsilon)\int^{t\epsilon^{-1}}_0\langle Z(s),\mathcal{B}\rangle_\epsilon ds +\mathcal{D}.
     \end{align}
     Simple change of variable reduces first term in \eqref{eq:Repestimates} to $2^{-1}\int_0^{t} \left\langle Z(s),\frac{\partial^2 \psi}{\partial x^2}\right\rangle_\epsilon ds$ at a cost of an error $\mathcal{D}$. It follows from the moment estimates in \eqref{eq:MomentEstimate1}-\eqref{eq:MomentEstimate3} that $\mathcal{D}$ in \eqref{eq:Repestimates} indeed converges to $0$ in $L^1$ as $\epsilon\to 0$. Further, we know the fact that $\mathcal{B}$ is in $L^2$. Thus, second term also converegs to $0$ in $L^1$. Thereafter, one can use DCT along with \eqref{eq:MartingaleCondOne} to conclude $L^1$ distance between $\int_0^t \langle Z(\epsilon^{-1}s),\partial^2\psi/\partial x^2\rangle_\epsilon ds$ and $\int_0^t \langle Z_\epsilon(s),\partial^2\psi/\partial x^2\rangle ds$ converges to $0$, thus shows the claim.

     \textit{Proof of \eqref{eq:MartingaleCondThree}}:
   One can write
 \begin{align}\label{eq:TotQuadVar}
\hat{N}(\epsilon^{-1}t) &= \epsilon^2\sum_{i=1}^{m-1} \int_{s_i}^{s_{i+1}} \sum_{\zeta_1,\zeta_2}\psi^{\zeta_1,\zeta_2}(s,s_{i+1}) Z(s,\zeta_1)Z(s,\zeta_2)d\langle M(s,\zeta_1), M(s,\zeta_2)\rangle\\
&+\epsilon^2\int_{s_m}^{t} \sum_{\zeta_1,\zeta_2}\psi^{\zeta_1,\zeta_2}(s,t) Z(s,\zeta_1)Z(s,\zeta_2)d\langle M(s,\zeta_1), M(s,\zeta_2)\rangle
\end{align}
         where \[\psi^{\zeta_1,\zeta_2}(s_1,s_2)=\left(\sum_{\xi\in \Xi(s_1,s_2)}p_\epsilon(s_1,s_2,\xi-\zeta_1)\psi(\epsilon\xi)\right)\left(\sum_{\xi\in \Xi(s_1,s_2)}p_\epsilon(s_1,s_2,\xi-\zeta_2)\psi(\epsilon\xi)\right).\]
        Recall from \eqref{eq:TotVarDiff} that
        \begin{align}
(b^{-1}&-1)^{-2}Z(t,\zeta_1)Z(t,\zeta_2)d\langle M(t,\zeta_1),M(t,\zeta_2)\rangle \\&= \frac{b^{\nu_\epsilon|\zeta_1-\zeta_2|}}{1-b^{\nu_\epsilon}}Z(t,\zeta_1\wedge \zeta_2)\left([p_\epsilon(t,t+dt)*Z(t)](\zeta_1\wedge \zeta_2)-\frac{(1-b^{\nu_\epsilon})b^{-1}}{b^{-1}-1}Z(t,\zeta_1\wedge \zeta_2)dt\right).
\end{align}
   Furthermore, one can write
   \[[p_\epsilon(t,t+dt)*Z(t)](x)-Z(t,x )dt= \sum_{k=0}^\infty (1-b^{\nu_\epsilon})b^{k\nu_\epsilon}\left(Z(t,x-k)-Z(t,x)\right)dt.\]
   Using \eqref{eq:MomentEstimate2}, we have
 \[\left|\sum_{k=0}^\infty (1-b^{\nu_\epsilon})b^{k\nu_\epsilon}\left(Z(t,x-k)-Z(t,x)\right)\right|\leq \sum_{k=0}^\infty \sqrt{\epsilon}|\epsilon k|^{v}\exp\left(-\lambda_\epsilon\nu_\epsilon k\sqrt{\epsilon}+k\epsilon\right)=O(\epsilon^{v/2})\]
 for any $t>0$.
       Moreover, one can deduce from \eqref{eq:UseOfSmoothness} that
      \[\sum_{\xi\in \Xi(s,s_{i+1})}p_\epsilon(s,s_{i+1},\xi-\zeta)\psi(\epsilon\xi)=\psi(\epsilon \zeta)+\epsilon\frac{1}{2}\frac{\partial^2 \psi}{\partial x^2}(\epsilon\zeta)(s_{i+1}-s)+o(\epsilon^2)\mathcal{B},\]
 where $\mathcal{B}$ is in $L^2(\RR)$. Thus, we get $\hat{N}(\epsilon^{-1}t)$ equal to 
 \begin{align}\label{eq:ModTotVAr}
   \sum_{i=1}^{m-1}\int_{s_i}^{s_{i+1}} \epsilon^2\sum_{\zeta_1,\zeta_2} \psi^{\zeta_1,\zeta_2}(s,s_{i+1})&b^{\nu_\epsilon|\zeta_1-\zeta_2|}(b^{-1}-1)\left(\frac{b^{\nu_\epsilon-1}-1}{1-b^{\nu_\epsilon}}Z^2(s,\zeta_1\wedge \zeta_2)+O(\epsilon^{v/2})\right)ds.
 \end{align}
   where $\psi^{\zeta_1,\zeta_2}(s,s_{i+1})=\psi(\epsilon\zeta_1)\psi(\epsilon\zeta_2)+O(\epsilon^2) \mathcal{B}_1$. Further, the random variable $\mathcal{B}_1$ belongs to $L^2$. Note that $\nu_\epsilon>0$ (see Definition~\ref{GartnerTransform}) is chosen in way so that it satisfies $1-\nu_\epsilon=\nu^{2}_\epsilon$. For this specific choice of $\nu_\epsilon$, we have
\[\sum_{\zeta_2\geq \zeta_1}b^{\nu_\epsilon|\zeta_1-\zeta_2|}(b^{-1}-1)\frac{b^{\nu_\epsilon-1}-1}{1-b^{\nu_\epsilon}}=1+O(\sqrt{\epsilon}).\]
            Thereafter, using continuity of $\psi$, one can say further $\hat{N}_\epsilon(\epsilon^{-1}t)-\int_0^t \langle Z^2(\epsilon^{-1}),\psi^2\rangle_\epsilon ds$ converges to $0$ in $L^1$ as $\epsilon\to 0$. To this end, \eqref{eq:MartingaleCondOne} shows $\int_0^t \langle Z^2(\epsilon^{-1}s),\psi^2\rangle_\epsilon ds$ and $\int_0^t \langle Z^2_\epsilon(s),\psi^2\rangle ds$ both has same $L^1$ limit and hence, completes the proof.

\appendix
\section{Fredholm Determinant}\label{App:ApendixA}

We start with the definition of Fredholm determinant.
\bd\label{FredholmDeterminant}
Let $K(x,y)$ be the meromorphic function of two complex variable, invariably referred as kernel in the literature. Let $\Gamma\subset \mathbb{C}$ be a curve. Assume that $K$ has no singularities on $\Gamma\times \Gamma$. Then the Fredholm determinant $\det(I+K)$ of the kernel $K$ is defined as the sum of the series of the complex integrals
\begin{align}\label{eq:FredDef}
\det(I+K)=\sum_{n\geq 0}\frac{1}{n!(2\pi i)^n}\int_\Gamma\ldots \int_\Gamma \det\left(K(z_i,z_j)\right)_{i,j=1}^n dz_1,\ldots dz_n.
\end{align}
\ed
\begin{remark}
In the usual definition of the kernel of Fredholm determinant (see, \cite[Chapter 5]{BS2005}), there is no such pre-factor of $1/(2\pi i)^n$ as in \eqref{eq:FredDef}. Furthermore, to show the absolute convergence of the series in \eqref{eq:FredDef}, one need to have suitable control on $\sup_{x\in \Gamma}|K(x,z)|$. To that end, following inequality (see, \cite[Exercise 1.1.3]{RajBhat97}) proves the absolute convergence.
\end{remark}

\bl[Hadamard Inequality]\label{HadamardInequality}
 If $D$ is a $N\times N$ complex matrix, then
 \begin{align}\label{eq:Hadamard}
 |\mbox{det}(D)|\leq \prod_{i=1}^N\sqrt{\sum_{j=1}^n D^2_{ij}}.
\end{align}
\el

\bc\label{DominatedConvCor}
Let $\Gamma$ be a curve. Kernel $K$ is defined over $\Gamma\times \Gamma$. Furthermore, $\sup_{x\in \Gamma}|K(x,z)|\leq \mathbf{K}(z)$, where $\mathbf{K}(z)$ satisfies $\int_\Gamma \mathbf{K}(z)dz< \infty $. Then series in \eqref{eq:FredDef} converges absolutely.
\ec


In the next lemma, we show that Fredholm determinant of the limiting kernel in \eqref{eq:FinalLimitOfKernelBelowCriticality} is the Gaussian distribution function. Although, we couldn't find a direct proof that result, we believe that it is quite well known. For the sake of clarity, we present here a short proof.
\bl\label{KernelConversion}
Consider a piecewise linear infinite curve $\Gamma$ which extends linearly from $-\delta +\infty e^{-\pi/6}$ to $-\delta$ and from $-\delta$ to $-\delta +\infty e^{\pi/6}$. Let us define the kernel $K_G$ for any two points $w,w^\prime\in \Gamma$  as
\begin{align}\label{eq:GaussianKernel}
K_G:=\frac{1}{2\pi i}\int_{-2\delta+\RR } \frac{\exp\left(\frac{v^2}{2}-\frac{w^2}{2}+s(v-w)\right)}{(v-w)(v-w^\prime)}\frac{v}{w}dv.
\end{align}
Then,
\begin{align}\label{eq:KernelConvergeToGaussian}
\det(1+K_G)_{L^2(\Gamma)}=\int_{-\infty}^s \frac{1}{\sqrt{2\pi}} e^{-z^2/2}dz.
\end{align}
\el
\begin{proof}
To begin with, notice that for all $v$ on the line $-2\delta+i\RR$, $\mbox{Re}(v-w)<0$ whenever $w\in \Gamma$. Thus, we can write $\exp(s(v-w))=(v-w)\int_s^\infty\exp(z(v-w))dz$. Fix any $\sigma\in \mathfrak{S}(n)$. To this end, we have
\begin{align}\label{eq:MultipleIntegral}
\int_\Gamma &\ldots \int_\Gamma K(w_1,w_{\sigma(1)})K(w_2,w_{\sigma(2)})\ldots K(w_n,w_{\sigma(n)}) dw_1,\ldots dw_n\\
&=\frac{1}{(2\pi i)^n}\int_{\Gamma^n} \int_{(-2\delta+i\RR)^n }\frac{\exp\left(\sum_{i=1}^n\left(\frac{v_i^2}{2}-\frac{w_i^2}{2}+s(v_i-w_i)\right)\right)}{\prod_{i=1}^n(v_i-w_i)(v_i-w_{\sigma(i)})}\frac{v_i}{w_i}\prod_{i=1}^n dv_i\prod_{i=1}dw_i\\
&=\frac{1}{(2\pi i)^n}\int_{\Gamma^n} \int_{(-2\delta+i\RR)^n }\int_{[-s,\infty)^n} \frac{\exp\left(\sum_{i=1}^n\left(\frac{v_i^2}{2}-\frac{w_i^2}{2}+z_i(v_i-w_i)\right)\right)}{\prod_{i=1}^n(v_i-w_{\sigma(i)})}\frac{v_i}{w_i}\prod_{i=1}^n dz_i\prod_{i=1}^n dv_i \prod_{i=1}^n dw_i\\
&=\frac{1}{(2\pi i)^n}\int_{[s,\infty)^n}\int_{(-2\delta+i\RR)^n}\int_{\Gamma^n} \frac{\exp\left(\sum_{i=1}^n\left(\frac{v_i^2}{2}-\frac{w_i^2}{2}+z_i(v_i-w_i)\right)\right)}{\prod_{i=1}^n(v_i-w_{\sigma(i)})}\frac{v_i}{w_i}\prod_{i=1}^n dw_i\prod_{i=1}^n dv_i \prod_{i=1}^n dz_i.
\end{align}
Define a piecewise linear curve $\Gamma_r$ which extends linearly from $-\delta +re^{-i\pi/6}$ to $-\delta$ and from $-\delta$ to $-\delta +re^{i\pi/6}$. Notice that $\Gamma_r$ tends to $\Gamma$ as $r\to \infty$. Denote the closed contour formed by $\Gamma$ and the line joining $-\delta +re^{-i\pi/6}$ to $-\delta +re^{i\pi/6}$ as $\tilde{\Gamma}_r$. As $r$ increases, integrand inside the paranthesis of the last line in \eqref{eq:MultipleIntegral} exponentially decays to zero uniformly for all $w\in \tilde{\Gamma}_r\backslash \Gamma_r$. Consequently, the integral over $\tilde{\Gamma}_r$ will converges to the integral over $\Gamma$. Thus, computing the integral over $\Gamma^n$ boils down to finding out the residue of the integrand at $z_1=\ldots =z_n=0$ because these are the only poles sitting inside the contour $\tilde{\Gamma}^n_r$ for all large value of $r$. Thereafter,
\begin{align}\label{eq:InnerIntSimplification}
\int_{\Gamma^n} &\frac{\exp\left(\sum_{i=1}^n\left(\frac{v_i^2}{2}-\frac{w_i^2}{2}+z_i(v_i-w_i)\right)\right)}{\prod_{i=1}^n(v_i-w_{\sigma(i)})}\frac{v_i}{w_i}\prod_{i=1}^n dw_i= \exp\left(\sum_{i=1}^n\left(\frac{v_i^2}{2}+z_iv_i\right)\right).
\end{align}
Using the analytic behavior of the function $\exp(2^{-1}v^2+vz)$, one can argue further
\begin{align}\label{eq:MiddleIntegral}
\int_{(-2\delta+i\RR)^n}\exp\left(\sum_{i=1}^n\left(\frac{v_i^2}{2}+z_iv_i\right)\right)\prod_{i=1}^n dv_i&= \int_{(-i\infty, i\infty)^n }\exp\left(\sum_{i=1}^n\left(\frac{v_i^2}{2}+z_iv_i\right)\right)\prod_{i=1}^n dv_i\\
&= (\sqrt{2\pi} i)^n \exp\left(-\sum_{i=1}^n z_i^2\right).
\end{align}
To this end, multiple integral in \eqref{eq:MultipleIntegral} becomes equal to $(1-\Phi(s))^n$ where $\Phi(.)$ is the distribution function of standard normal distribution. This implies for all $n\geq 2$
\begin{align}\label{eq:HigherKerDetZero}
\int_\Gamma\ldots \int_\Gamma \det(K(w_i,w_j))_{i,j=1}^n dw_1\ldots dw_n =0.
\end{align}
Thus, we have $\det(1+K_G)_{L^2(\Gamma)}= 1-(1-\Phi(s))=\Phi(s)$. Hence, this completes the proof.
\end{proof}

 \section{eigenfunction of the Transition Matrix in HL-PushTASEP}\label{suppA}

\textsc{Proof of Proposition~\ref{EigenFunctionLemma}:} Here, we sketch a proof of the present proposition extending the arguments in \citep[Theorem 3.4]{BCG2016} in the continuous case. Let us define $N$ functions $g_i(x_{i},x_{i+1},\ldots ,x_N;z_{i},z_{i+1},\ldots ,z_N)$ as follows
\begin{align}\label{eq:TestFunctionDef}
g_i(x_{i},x_{i+1},\ldots ,x_N; z_1,z_2,\ldots ,z_N)= \sum_{(y_i,\ldots ,y_N)}\mathcal{T}^{(t,t+dt)}_i \left((x_i,\ldots ,x_N)\to (y_i,\ldots,y_N)\right) z^{y_i}_i\ldots z^{y_N}_N
\end{align}
where $\mathcal{T}^{t,t+dt}_i(.)$ denotes the probability of the transition from the configuration $(x_i,\ldots,x_N)$ to $(y_i,\ldots,y_N)$ once the particle at  position $x_i$ gets excited in the time interval $(t,t+dt)$. For convenience, let us consider
\begin{align}
g_i(x_{i},x_{i+1},\ldots ,x_N; z_{i},z_{i+1},\ldots ,z_N)&:= g^{(1)}(x_{i},x_{i+1},\ldots ,x_N; z_{i},z_{i+1},\ldots ,z_N)\\&+g^{(2)}_i(x_{i},x_{i+1},\ldots ,x_N; z_{i},z_{i+1},\ldots ,z_N).\label{eq:TransferMatrixBreak}
\end{align}
Subsequently, we describe two new functions introduced in \eqref{eq:TransferMatrixBreak} as follows. In the expansion of $g_i$, one can note the contribution of the event that the particle at position $x_i$  moves by $j$-steps where $1\leq j<x_{i+1}-x_{i}$ in the time gap $(t,t+dt)$ is $\frac{z_i}{1-bz_i}(z_{i}^{x_i}-b^{x_{i+1}-x_i-1}z_i^{x_{i+1}-1})z^{x_{i+1}}_{i+1}\ldots z^{x_N}_N$. Thus, define
\begin{align}\label{eq:g1Def2}
g^{(1)}_i(x_i,\ldots ,x_N; z_i,\ldots , z_N )= \frac{z_i(1-b)}{1-bz_i}(z_{i}^{x_i}-b^{x_{i+1}-x_i-1}z_i^{x_{i+1}-1})z^{x_{i+1}}_{i+1}\ldots z^{x_N}_N .
\end{align}
Lastly, $g^{(2)}$ takes account of the event that particle at $x_i$ jumps to the position $x_{i+1}$ and consequently, particle which was previously there in $x_{i+1}$ starts to hop to the right in the same fashion. To this end, one can write
\begin{align}\label{eq:g3Def3}
g^{(2)}_i(x_i,\ldots ,x_N; z_i,\ldots , z_N )= b^{x_{i+1}-x_{i}-1}z_i^{x_{i+1}}g_{i+1}(x_{i+1},\ldots ,x_N; z_i,\ldots , z_N ).
\end{align}
Let us define a function $g$ as
\begin{align}\label{eq:gDef}
g(x_1,\ldots ,x_N;z_1,\ldots ,z_N)&:=\sum_{(y_1,\ldots ,y_N)}\mathcal{T}^{(t,t+dt)} \left((x_1,\ldots ,x_N)\to (y_1,\ldots,y_N)\right) z^{y_i}_i\ldots z^{y_N}_N
\end{align}
where $\mathcal{T}^{(t,t+dt)}$ denotes the transfer matrix in the interval $(t,t+dt)$. Naively, $g$ acts as a probability generating function in the infinitesimal time interval $(t,t+dt)$. In what follows, it is easy to see how $g$ and $g_i$'s are connected
\begin{align}\label{eq:G-GiConnection}
g(x_1,\ldots ,x_N;z_1,\ldots ,z_N)=\sum_{i=1}^N z_1^{x_1}\ldots z_{i-1}^{x_{i-1}}g_i(x_{i},x_{i+1},\ldots ,x_N; z_{i},z_{i+1},\ldots ,z_N)dt.
\end{align}
If we take $N=2$, then we get
\[g(x_1,x_2;z_1,z_2)=\frac{z_1(1-b)}{1-bz_1}(z_{1}^{x_1}-b^{x_{2}-x_1-1}z_1^{x_{2}-1})z^{x_{2}}_{2}dt+\frac{z_2^{x_2+1}(1-b)}{1-bz_2}(b^{x_{2}-x_{1}-1}z^{x_{2}}_1+z^{x_{1}}_1)dt.\]
Let us consider $\mathcal{S}_{12}g(x_1,x_2;z_1,z_2)+\mathcal{S}_{21}g(x_1,x_2;z_2,z_1)$ for some functions $\mathcal{S}_{12}$ and $\mathcal{S}_{21}$ of $z_1,z_2$. One can write
\begin{align}\label{eq:EigenfunctionHunt}
\mathcal{S}_{12}g(x_1,x_2&;z_1,z_2)+\mathcal{S}_{21}g(x_1,x_2;z_2,z_1)  = \left(\frac{z_1(1-b)}{1-bz_1}+\frac{z_2(1-b)}{1-bz_2}\right)(\mathcal{S}_{12}z^{x_1}_1z^{x_2}_2+\mathcal{S}_{21}z^{x_1}_2z^{x_2}_1)dt\\
&+ (1-b)(z_1z_2)^{x_2}b^{x_2-x_1-1}\left(\mathcal{S}_{12}\left[\frac{z_2}{1-bz_2}-\frac{1}{1-bz_1}\right]+\mathcal{S}_{21}\left[\frac{z_1}{1-bz_1}-\frac{1}{1-bz_2}\right]\right)dt.
\end{align}
Thus, on taking $\mathcal{S}_{12}=1-(1+b)z_1+bz_1z_2$ and $\mathcal{S}_{12}=-1+(1+b)z_2-bz_1z_2$, one can note that $\mathcal{S}_{12}z^{x_1}_1z^{x_2}_2+\mathcal{S}_{21}z^{x_2}_1z^{x_1}_2$ is an eigenfunction of the transfer matrix $\mathcal{T}$ with the eigenvalue $\frac{z_1(1-b)}{1-bz_1}+\frac{z_2(1-b)}{1-bz_2}$ in the interval $(t,t+dt)$. One can also scale $\mathcal{S}_{12}$ and $\mathcal{S}_{21}$ by $1-(1+b)z_1+bz_1z_2$ without any major consequences. To obtain the eigenfunction and eigenvalues for any $N$, we use induction. Assume $N=n-1$ and
\begin{align}\label{eq:InductionStep}
\sum_{\sigma\in \mathfrak{S}(n-1)}A_{\sigma}&g(x_1,x_2,\ldots ,x_{n-1};z_{\sigma(1)},z_{\sigma(2)},\ldots ,z_{\sigma(n-1)})\\&=(1-b)\left(\sum_{i=1}^{n-1}\frac{z_i}{1-bz_i}\right)\left(\sum_{\sigma\in \mathfrak{S}(n-1)}A_{\sigma}z^{x_1}_{\sigma(1)}\ldots z^{x_{n-1}}_{\sigma(n-1)}\right)dt
\end{align}
where
\[A_\sigma:=\prod_{1\leq i<j\leq n-1}(-1)^{\sigma}\frac{1-(1+b)z_{\sigma(i)}+bz_{\sigma(i)}z_{\sigma(j)}}{1-(1+b)z_{i}+bz_{i}z_{j}}.\]
We need  to show that \eqref{eq:InductionStep} holds even when $N=n$. Consider the following
\begin{align}
\sum_{\sigma\in \mathfrak{S}(n)}A_{\sigma}&g(x_1,x_2,\ldots ,x_{n};z_{\sigma(1)},z_{\sigma(2)},\ldots ,z_{\sigma(n)})=(I)+(II)+(III)
\end{align}
where
\begin{align}\label{eq:DefI}
(I):&= \sum_{\stackrel{(y_1,y_2,\ldots ,y_n)}{y_1=x_1}}\mathcal{T}^{(t,t+dt)}\left((x_1,\ldots ,x_n)\to(y_1,\ldots ,y_n)\right)\left(\sum_{j=1}^n\sum_{\stackrel{\sigma\in \mathfrak{S}(n)}{\sigma(1)=j}}A_\sigma\prod z^{y_i}_{\sigma(i)}\right),\\
&=\sum_{j=1}^n  (1-b)\left(\sum_{\stackrel{i=1}{i\neq j}}^{n}\frac{z_i}{1-bz_i}\right)\left(\sum_{\sigma\in \mathfrak{S}(n)}A_{\sigma}z^{x_1}_{\sigma(1)}z^{x_2}_{\sigma(2)}\ldots z^{x_n}_{\sigma(n)}\right)dt,
\end{align}

\begin{align}\label{eq:DefII}
(II):&= \sum_{\stackrel{(y_1,y_2,\ldots ,y_n)}{x_1<y_1<x_{2}}}\mathcal{T}^{(t,t+dt)}\left((x_1,\ldots ,x_n)\to(y_1,\ldots ,y_n)\right)\left(\sum_{\sigma\in \mathfrak{S}(n)}A_\sigma\prod z^{y_i}_{\sigma(i)}\right)\\
&=(1-b)\sum_{\sigma\in \mathfrak{S}(n)}A_\sigma\left(\frac{z_{\sigma(1)}}{1-bz_{\sigma(1)}}-\frac{b^{x_2-x_1-1}z^{x_2-x_1}_{\sigma(1)}}{1-bz_{\sigma(1)}}\right)z^{x_1}_{\sigma(1)}z^{x_2}_{\sigma(2)}\ldots z^{x_n}_{\sigma(n)}dt,
\end{align}
and
\begin{align}\label{eq:DefIII}
(III):&=\sum_{\stackrel{(y_1,y_2,\ldots ,y_n)}{y_1=x_{2}}}\mathcal{T}^{(t,t+dt)}\left((x_1,\ldots ,x_n)\to(y_1,\ldots ,y_n)\right)\left(\sum_{\sigma\in \mathfrak{S}(n)}A_\sigma\prod z^{y_i}_{\sigma(i)}\right)\\
&=(1-b)b^{x_2-x_1-1}\sum_{\sigma\in \mathfrak{S}(n)}A_\sigma\left(\frac{z_{\sigma(2)}}{1-bz_{\sigma(2)}}-\frac{b^{x_3-x_2-1}z^{x_3-x_2}_{\sigma(2)}}{1-bz_{\sigma(2)}}\right)z^{x_1}_{\sigma(1)}z^{x_2}_{\sigma(2)}\ldots z^{x_n}_{\sigma(n)}dt\\
&+\sum_{\stackrel{(y_1,y_2,\ldots ,y_n)}{y_1=x_{2},y_2=x_{3}}}\mathcal{T}^{(t,t+dt)}\left((x_1,\ldots ,x_n)\to(y_1,\ldots ,y_n)\right)\left(\sum_{\sigma\in \mathfrak{S}(n)}A_\sigma\prod z^{y_i}_{\sigma(i)}\right).
\end{align}
For any $\sigma\in\mathfrak{S}(n)$, one can note that
\[A_\sigma\left[\frac{z^{x_2-x_1}_{\sigma(1)}}{1-bz_{\sigma(1)}}-\frac{z_{\sigma(2)}}{1-bz_{\sigma(2)}}\right](z_{\sigma(1)}z_{\sigma(2)})^{x_2}=-A_{\sigma^\prime}\left[\frac{z^{x_2-x_1}_{\sigma^\prime(1)}}{1-bz_{\sigma^\prime(1)}}-\frac{z_{\sigma^\prime(2)}}{1-bz_{\sigma^\prime(2)}}\right](z_{\sigma^\prime(1)}z_{\sigma^\prime(2)})^{x_2}\]
where $\sigma^\prime=(1,2)*\sigma$. Thus, there will be telescopic cancellations of large number of terms in the sum $(II)+(III)$. To that effect, we would have
\begin{align}\label{eq:IIandIIISum}
(II)+(III)= (1-b)\sum_{\sigma\in \mathfrak{S}(n)}A_\sigma\frac{z_{\sigma(1)}}{1-bz_{\sigma(1)}}z^{x_1}_{\sigma(1)}z^{x_2}_{\sigma(2)}\ldots z^{x_n}_{\sigma(n)}dt
\end{align}
and consequently, adding \eqref{eq:DefI} and \eqref{eq:IIandIIISum} shows \eqref{eq:InductionStep} when $N=n$. To this end, denoting the transfer matrix for exactly $m$-many clock counts in the interval $(0,t)$ by $\mathcal{T}^{(0,t)}_{m}(.)$, one can deduce
\begin{align}\label{eq:PoissonInt}
&\sum_{(y_1,\ldots ,y_N)}\mathcal{T}^{(0,t)}_{m}\left((x_1,\ldots ,x_N)\to (y_1,\ldots ,y_N)\right)\left(\sum_{\sigma\in  \mathfrak{S}(n)}A_\sigma z^{y_1}_{\sigma(1)}\ldots z^{y_N}_{\sigma(N)}\right)\\
&=\int_{(t_0,\ldots,t_{m+1})\in \Delta(t)}\left(\sum_{i=1}^N\frac{(1-b)z_i}{(1-bz_i)}\right)^m\left(\sum_{\sigma\in  \mathfrak{S}(n)}A_\sigma z^{x_1}_{\sigma(1)}\ldots z^{x_N}_{\sigma(N)}\right)\prod_{k=1}^{m} \exp(-N(t_k-t_{k-1}))dt_k\\
&=\exp(-Nt)\frac{1}{m!}\left(t\sum_{i=1}^N\frac{(1-b)z_i}{(1-bz_i)}\right)^m \left(\sum_{\sigma\in  \mathfrak{S}(n)}A_\sigma z^{x_1}_{\sigma(1)}\ldots z^{x_N}_{\sigma(N)}\right)
\end{align}
where $\Delta(t_1,\ldots ,t_{m+1})=\left\{0=t_0 <t_1<t_2<\ldots <t_{m+1}=t\right\}$. Note that $\exp(-N(t_k-t_{k-1}))$ term inside the integral above accounts for the probability that $N$ independent clocks associated with $N$ particles do not ring in the interval $(t_{k-1},t_{k})$. Once any of the $N$ clock rings, it results in a factor of $\sum_{i=1}^N\frac{(1-b)z_i}{(1-bz_i)}$ at the front. Thus, $m$ many clock counts add upto $m$ such factors. Using \eqref{eq:PoissonInt}, one can further deduce
\begin{align}\label{eq:CollectingAll}
\sum_{(y_1,\ldots ,y_n)}&\mathcal{T}^{(N)}_t\left((x_1,\ldots ,x_N)\to (y_1,\ldots, y_N)\right)\left(\sum_{\sigma\in \mathfrak{S}(N)}A_\sigma z^{y_1}_{\sigma(1)}\ldots z^{y_N}_{\sigma(N)}\right)\\
&=\sum_{m=0}^\infty\sum_{(y_1,\ldots ,y_n)}\mathcal{T}^{(0,t)}_m\left((x_1,\ldots ,x_N)\to (y_1,\ldots, y_N)\right)\left(\sum_{\sigma\in \mathfrak{S}(N)}A_\sigma z^{y_1}_{\sigma(1)}\ldots z^{y_N}_{\sigma(N)}\right)\\
&=\sum_{m=0}^\infty\exp(-Nt)\frac{1}{m!}\left(t\sum_{i=1}^N\frac{(1-b)z_i}{(1-bz_i)}\right)^m \left(\sum_{\sigma\in  \mathfrak{S}(n)}A_\sigma z^{x_1}_{\sigma(1)}\ldots z^{x_N}_{\sigma(N)}\right)\\
&=\exp\left(-\sum_{i=1}^N\frac{1-z_i}{1-bz_i}\right)\left(\sum_{\sigma\in  \mathfrak{S}(n)}A_\sigma z^{x_1}_{\sigma(1)}\ldots z^{x_N}_{\sigma(N)}\right).
\end{align}
This completes the proof.

\section{Moment Formula}\label{suppC}

\textsc{Proof of Proposition~\ref{MGF}:}  
In what follows, we first find out the moment formula with the initially $N$ many particles. First, we claim that under step Bernoulli initial condition the probability of transition of $m$-th particle to $x$ by time $t$ in HL-PushTASEP is given as 
\begin{align}\label{eq:OneDimProb}
\mathbb{P}_{stepb}&(x_m(t)=x;t)= (-1)^{m-1}b^{m(m-1)/2}\sum_{m\leq k\leq N}\sum_{|S|=k}b^{\kappa(S,\ZZ_{>0})-mk-k(k-1)/2}\rho^{S_k}\binom{k-1}{m-1}_b\nonumber\\
&\times \frac{1}{(2\pi i)^k}\oint\ldots \oint \prod_{1\leq i<j\leq k}\frac{z_j-z_i}{1-(1+b^{-1})z_i+b^{-1}z_iz_j)}\prod_{i=1}^k \left(\frac{1}{1-(1-\rho)z^{-1}_i\cdots z^{-1}_k}\right)^{S_i-S_{i-1}}\nonumber\\
&\times \frac{1-\prod_{i=1}^k z_i}{\prod_{i=1}^k (1-z_i)}\prod_{i=1}^k z_i^{x-S_i-1}\exp\left(-t\frac{1-z^{-1}_i}{1-bz_i^{-1}}\right) dz_i,
\end{align}
where $S=\{S_1,\ldots ,S_k\}$ with $1\leq S_1<S_2<\ldots <S_k\leq N$  and $\kappa(S,\ZZ_{>0})$ counts the sum of the elements in $S$.
One  can note that when HL-PushTASEP is started with particles at $N$ ordered places $(y_1,\ldots ,y_n)$ (call it $\mathbf{Y}$) in $\ZZ_{\geq 0}$, we can have an explicit expression of $
\mathbb{P}_{\vec{y}}(x_m(t)=x;t)$ given in Theorem~\ref{TransitProb}. To put all the ingredients together, let us just recall the formula of $\mathbb{P}_{\vec{y}}(x_m(t)=x;t)$ given there. For a fixed $N\geq 0$, when the initial data is $\vec{y}=(y_1,\ldots ,y_N)$, then
\begin{align}\label{eq:OneDimProbRepeat}
\mathbb{P}_{\vec{y}}(x_m(t)=x;t)&= (-1)^{m-1}b^{m(m-1)/2}\sum_{m\leq k\leq N}\sum_{|S|=k}b^{\kappa(S,\ZZ_{> 0})-mk-k(k-1)/2}\binom{k-1}{m-1}_b\nonumber\\
&\times \frac{1}{(2\pi i)^k}\oint\ldots \oint \prod_{i,j\in S, i<j}\frac{z_j-z_i}{1-(1+b^{-1})z_i+b^{-1}z_iz_j)}\nonumber\\
&\times \frac{1-\prod_{i\in S} z_i}{\prod_{i\in S} (1-z_i)}\prod_{i\in S} z_i^{x-y_i-1}\exp\left(-t\frac{1-z^{-1}_i}{1-bz_i^{-1}}\right) dz_i,
\end{align}
where $\kappa(S,\ZZ_{>0})$ is the sum of the element in $S$. The term which relates initial data with $
\mathbb{P}_{\vec{y}}(x_m(t)=x;t)$ in \eqref{eq:OneDimProbRepeat}, is $z^{-y_i}_i$ for $i\in S$. To obtain the result in case of step initial condition, we have to sum over all $\vec{y}$ by multiplying $
\mathbb{P}_{\vec{y}}(x_m(t)=x;t)$ with appropriate probability. Let us introduce $k$-many new variables $t_1,\ldots ,t_k$ which are defined as
\[t_1:=S_1-1,t_2:=S_1-S_2-1,\ldots , t_k:=S_k-S_{k-1}-1.\]
To this end, we can write down $S_i$ in terms of $t_i$ by $S_i=\sum_{j=1}^i t_i +i$. So, initially, we would have the configuration $y_{S_1}=S_1+i_i,\ldots , y_{S_k}= S_k+\sum_{l=1}^k i_l$ with probability
\[\rho^{S_k}\prod_{l=1}^k \binom{t_l+i_l}{t_l}(1-\rho)^{i_l}.\]
 One can note that
\begin{align}\label{eq:3}\sum_{i_1,\ldots ,i_k\geq 0}\prod_{j=1}^k \binom{t_j+i_j}{t_j}(1-\rho)^{i_j}z_{S_j}^{-S_j+\sum_{l=1}^j i_l}&=\prod_{j=1}^k z_{S_j}^{-S_j}\left(\frac{1}{1-(1-\rho)z^{-1}_{S_j}\cdots z^{-1}_{S_k}}\right)^{t_j+1}.
\end{align}
Now, $\prod_{j=1}^k z_{S_j}^{-S_j}$ in \eqref{eq:3} counts for the replacement of $\prod_{i\in S} z_i^{x-y_i-1}$ in \eqref{eq:OneDimProbRepeat} with $\prod_{i=1}^k z_i^{x-S_i-1}$ in \eqref{eq:OneDimProb}, whereas extra factor in second line of \eqref{eq:OneDimProb} comes from the second term at the end of the right side in \eqref{eq:3}. Thus, the claim is proved. 

 Now, we turn to showing the moment formula. Given any $L\in \NN$, we multiply \eqref{eq:OneDimProb} by $b^{mL}$ and sum over all $1\leq m\leq N$. Let us note down the following result (see, \cite[Theorem 10.2.1]{andrews1999}) known as $q$- binomial Theorem. For any $q,t\in \mathbb{C}$ such that $|q|\neq 1$ and any positive integer $n$, we have
\begin{align}\label{eq:QBinomTheo}
\prod_{i=0}^{n-1}(1+q^it)=\sum_{k=0}^n q^{k(k-1)/2}\binom{n}{k}_q t^k.
\end{align}
Thus, using \eqref{eq:QBinomTheo}, one can write
\[\sum_{m=1}^k(-1)^{m-1}b^{m(m-1)/-mk+mL}\binom{k-1}{m-1}_b=b^{L-k}\prod_{i=0}^{k-2}(1-b^{1-k+L}.b).\]
 To that effect, we get
 \begin{align}\label{eq:ArbInitDataMomFol}
 \mathbb{E}_{stepb}&\left(b^{LN_x(t)}\eta_x(t)\right)=\sum_{k=1}^{\min(L,N)}\left(\prod_{i=0}^{k-2}(1-b^{1-k+L}.b)\right)\sum_{|S|=k}b^{\kappa(S,\ZZ_{>0})-k(k-1)/2+L-k}\rho^{S_k}\nonumber\\
 &\times \frac{1}{(2\pi i)^k}\oint\ldots \oint \prod_{1\leq i<j\leq k}\frac{z_j-z_i}{1-(1+b^{-1})z_i+b^{-1}z_iz_j}\prod_{i=1}^k \left(\frac{1}{1-(1-\rho)z^{-1}_i\cdots z^{-1}_k}\right)^{S_i-S_{i-1}}\nonumber\\
&\times \frac{1-\prod_{i=1}^k z_i}{\prod_{i=1}^k (1-z_i)}\prod_{i=1}^k z_i^{x-S_i-1}\exp\left(-t\frac{1-z^{-1}_i}{1-bz_i^{-1}}\right) dz_i,
\end{align}
where $S=\{S_1,\ldots ,S_k\}$ with $1\leq S_1<S_2<\ldots <S_k\leq N$  and $\kappa(S,\ZZ_{>0})$ counts the sum of the elements in $S$. Letting $N\to\infty$, we sum \eqref{eq:ArbInitDataMomFol} over all possible $k$ ordered subsets $1\leq S_1<S_2<\ldots <S_k$ of various sizes. 
This infinite sum can be simplified in the following way 
\begin{align}\label{eq:SumOfAllSubset}
\sum_{1\leq S_1<S_2<\ldots <S_k}b^{S_1+\ldots +S_k}\rho^{S_k}\prod_{i=1}^k z^{-S_i}_i\left(\frac{1}{1-(1-\rho)z^{-1}_i\cdots z^{-1}_k}\right)^{S_i-S_{i-1}}&=\prod_{i=1}^k \frac{\frac{\rho (b/z_i)\cdots(b/z_k)}{1-(1-\rho)z^{-1}_i\cdots z^{-1}_k}}{1-\frac{\rho (b/z_i)\cdots(b/z_k)}{1-(1-\rho)z^{-1}_i\cdots z^{-1}_k}}\\
=b^{k(k+1)/2}\rho^k \frac{1}{z_i\cdots z_k-(1-\rho)-b^{k-i+1}\rho}.
\end{align}
In what  follows, we use a symmetrization identity from \cite[Section III]{Tracy2009b} to conclude that
\begin{align}\label{eq:IntermediateDStep}
\mathbb{E}_{stepb}\left(b^{LN_x(t)}\eta_x(t)\right)&=\sum_{k=1}^L\left(\prod_{i=1}^{k-2}(1-b^{1-k+L}b^i)\right)\frac{b^{-k(k-1)/2+L-k}\rho^k}{(b^{-1};b^{-1})_k}\nonumber\\
 &\times \frac{1}{(2\pi i)^k}\oint\ldots \oint \prod_{1\leq i<j\leq k}\frac{z_j-z_i}{1-(1+b^{-1})z_i+b^{-1}z_iz_j}\left(1-\prod_{i=1}^k z_i\right)\nonumber\\
 &\times \prod_{i=1}^k\frac{z_i^x(b^{-1}-1)}{(\rho+(1-\rho)b^{-1}-z_ib^{-1})(1-z_i)}\exp\left(-t\frac{1-z_i^{-1}}{1-bz_i^{-1}}\right)dz_i.
\end{align}
Moreover, following recursion
\begin{align}\label{eq:Recursion}
b^{LN_x(t)}=b^{LN_{x-1}(t)}+b^{LN_x(t)}\eta_t(x)\left(1-b^{-L}\right)=\ldots &=(1-b^{-L})\sum_{y=-\infty}^x\eta_t(x)b^{LN_x(t)} +\lim_{y\to \infty}b^{LN_y(t)}\\
&=(1-b^{-L})\sum_{y=-\infty}^x\eta_t(x)b^{LN_x(t)} +1
\end{align}
 cancels out the term $(1-\prod_{i=1}^k z_i)$ from \eqref{eq:IntermediateDStep}, thus yielding $\mathbb{E}(b^{LN_x(t)})$ in \eqref{eq:5} from  $\mathbb{E}(b^{LN_x(t)}\eta_x(t))$.
\end{proof}

\bibliographystyle{alpha}
\bibliography{Reference}

\end{document}